\documentclass[11pt]{amsart}
\usepackage{pb-diagram,amssymb,epic,eepic,verbatim,graphicx}
\usepackage{graphics,epsfig,psfrag}

\newtheorem{theorem}{Theorem}[subsection]
\newtheorem{corollary}[theorem]{Corollary}

\newtheorem{conjecture}[theorem]{Conjecture}
\newtheorem{proposition}[theorem]{Proposition}
\newtheorem{lemma}[theorem]{Lemma}
\newtheorem{lem}[theorem]{}
\theoremstyle{definition}
\newtheorem{definition}[theorem]{Definition}
\theoremstyle{remark}
\newtheorem{remark}[theorem]{Remark}
\newtheorem{example}[theorem]{Example}
\newcommand{\blem}{\begin{lem} \rm}
\newcommand{\elem}{\end{lem}}
\newcommand\TD{\mathcal{T} }
\newcommand\TDN{\mathcal{T}^N }

\newcommand\E{\mathcal{E}}
\newcommand\M{\mathcal{M}}
\newcommand\D{\mathcal{D}}

\newcommand\cA{\mathcal{A}}

\renewcommand\M{\mathcal{M}}
\renewcommand\S{\mathcal{S}}

\newcommand{\J}{\mathcal{J}}

\newcommand{\U}{\mathcal{U}}

\newcommand{\N}{\mathbb{N}}
\newcommand{\R}{\mathbb{R}}

\newcommand{\C}{\mathbb{C}}
\newcommand{\RP}{\mathbb{RP}}
\newcommand{\CP}{\mathbb{CP}}
\newcommand{\cC}{\mathcal{C}}

\newcommand{\Z}{\mathbb{Z}}

\renewcommand{\P}{\mathbb{P}}

\newcommand\lie[1]{\mathfrak{#1}}

\newcommand{\g}{\lie{g}}

\newcommand{\on}{\operatorname}

\newcommand{\Ham}{\on{Ham}}

\renewcommand{\top}{{\on{top}}}

\newcommand{\graph}{\on{graph}}
\newcommand{\Symp}{\on{Symp}}

\newcommand{\Lag}{\on{Lag}}
\newcommand{\Fr}{\on{Fr}}
\newcommand{\Sp}{\on{Sp}}

\newcommand{\End}{\on{End}}

\newcommand{\Ind}{ \on{Ind}}
\renewcommand{\ker}{ \on{ker}}

\newcommand{\im}{ \on{im}}

\newcommand{\codim}{\on{codim}}

\newcommand\dirac{/\kern-1.2ex\partial} % Dirac operator
\newcommand\qu{/\kern-.7ex/} % Categorical quotients
\newcommand\lqu{\backslash \kern-.7ex \backslash} % Categorical
\newcommand\dr{r_+ \kern-.7ex - \kern-.7ex r_-}
 
\def\pd{\partial}

\renewcommand{\d}{{\mbox{d}}}

\newcommand\Phinv{\Phi^{-1}}
\newcommand\muinv{\mu^{-1}}

\newcommand\eps{\epsilon}

\newcommand{\f}{\frac}
\newcommand{\lan}{\langle}
\newcommand{\ran}{\rangle}
\newcommand{\hh}{{\f{1}{2}}}

\newcommand{\ti}{\tilde}

\newcommand\pt{\on{pt}}
\newcommand\cE{\mathcal{E}}
\newcommand\cL{\mathcal{L}}
\newcommand\cK{\mathcal{K}}

\newcommand\curv{\on{curv}}

\newcommand\Map{\on{Map}}

\newcommand\Vect{\on{Vect}}
\newcommand\ul{\underline}

\renewcommand\Im{\on{Im}}

\newcommand\reg{{\on{reg}}}
\newcommand\bra[1]{ \lan {#1} \ran} 
\newcommand\bdefn{\begin{definition}}
\newcommand\edefn{\end{definition}}
\newcommand\bea{\begin{eqnarray*}}
\newcommand\eea{\end{eqnarray*}}
\newcommand\bcv{\left[ \begin{array}{r} }
\newcommand\ecv{\end{array} \right] }

\newcommand\bma{\left[ \begin{array} }
\newcommand\ema{\end{array} \right]}
\newcommand\ben{\begin{enumerate}}
\newcommand\een{\end{enumerate}}
\newcommand\beq{\begin{equation}}
\newcommand\eeq{\end{equation}}
\newcommand\bex{\begin{example}}
\newcommand\bsj{\left\{ \begin{array}{rrr} }
\newcommand\esj{\end{array} \right\}}

\newcommand\Id{\on{Id}}
\newcommand\cI{\mathcal{I}}
\newcommand\cP{\mathcal{P}}

\newcommand\eex{\end{example}}

\newcommand\sx{*\kern-.5ex_X}

\def\mathunderaccent#1{\let\theaccent#1\mathpalette\putaccentunder}
\def\putaccentunder#1#2{\oalign{$#1#2$\crcr\hidewidth \vbox
to.2ex{\hbox{$#1\theaccent{}$}\vss}\hidewidth}}

\newcommand{\odd}{{\on{odd}}}
\newcommand{\even}{{\on{even}}}

\begin{document}

\title{Quilted Floer cohomology}

\author{Katrin Wehrheim and Chris T. Woodward}

\address{Department of Mathematics,
Massachusetts Institute of Technology,
Cambridge, MA 02139.
{\em E-mail address: katrin@math.mit.edu}}

\address{Department of Mathematics, 
Rutgers University,
Piscataway, NJ 08854.
{\em E-mail address: ctw@math.rutgers.edu}}

\begin{abstract}  
We generalize Lagrangian Floer cohomology to sequences of Lagrangian
correspondences. For sequences related by the geometric composition of Lagrangian correspondences we establish an isomorphism of the Floer cohomologies.
This provides the foundation for the construction of a symplectic $2$-category as well as for the definition of topological invariants via decomposition and representation in the symplectic category.
Here we give some first direct symplectic applications: Calculations of Floer cohomology, displaceability of Lagrangian correspondences, and transfer of displaceability under geometric composition.

\vspace{-10mm}
\end{abstract} 

\maketitle

\tableofcontents

\section{Introduction}  

Lagrangian Floer cohomology associates to a pair of compact Lagrangian manifolds
a chain complex whose differential counts pseudoholomorphic strips 
with boundary values in the given Lagrangians.  
In this paper we generalize Floer cohomology to include compact Lagrangian correspondences.
Recall that if $(M_0,\omega_0)$ and $(M_1,\omega_1)$ are symplectic manifolds, then a {\em Lagrangian correspondence} $L_{01}$ from $M_0$ to $M_1$ is a Lagrangian submanifold  $L_{01}\subset M_0^- \times M_1$, where $M_0^-:= (M_0,-\omega_0)$. 
These were introduced by Weinstein \cite{we:sc} in an attempt to create a symplectic category with morphisms between not necessarily symplectomorphic manifolds.
So we also denote a Lagrangian correspondence by $M_0 \overset{L_{01}}{\longrightarrow} M_1$.
With this notation we can view a pair of Lagrangian submanifolds $L,L'\subset M$ as sequence of Lagrangian correspondences 
$\pt \overset{L}{\longrightarrow} M \overset{L'}{\longrightarrow} \pt $
from the point via $M$ back to the point.
This is a special case of a {\em cyclic sequence of Lagrangian correspondences}
$$
M_0 \overset{L_{01}}{\longrightarrow} M_1\overset{L_{12}}{\longrightarrow} M_2 \; \ldots  M_r
\overset{L_{r(r+1)}}{\longrightarrow} M_{r+1}=M_0 ,
$$
for which we will define a {\em quilted Floer cohomology} 
\begin{equation} \label{eq:HF}
HF(L_{01},L_{12}, \ldots, L_{r(r+1)}).
\end{equation}
The quilted differential counts tuples of pseudoholomorphic strips 
$(u_j:\R\times[0,1]\to M_j)_{j=0,\ldots,r}$ whose boundaries match up via the Lagrangian correspondences, $(u_{j}(s,1),u_{j+1}(s,0))\in L_{j(j+1)}$ for $j=0,\ldots, r$.
These tuples are examples of pseudoholomorphic quilts with the strips thought of as patches and the boundary matching conditions thought of as seams. The theory of quilts is developed in higher generality in \cite{quilts}.
In this paper, we next investigate the effect of geometric composition on Floer cohomology.
The {\em geometric composition} of two Lagrangian correspondences $L_{01}\subset M_0^-\times M_1$, $L_{12}\subset M_1^-\times M_2$ is
\begin{equation}\label{alg comp}
L_{01} \circ L_{12} := \bigl\{ (x_0,x_2)\in M_0\times M_2 \,\big|\, \exists x_1 :
(x_0,x_1)\in L_{01}, (x_1,x_2)\in L_{12} \bigr\} .
\end{equation}
In general, this will be a singular subset of $M_0^-\times M_2$. However, if we assume transversality of the intersection
$ L_{01} \times_{M_1} L_{12} :=
\bigl( L_{01} \times L_{12}\bigr) \cap \bigl(M_0^- \times \Delta_{M_1}\times M_2 \bigr)$,
then the restriction of the projection $ \pi_{02}: M_0^- \times M_1 \times
M_1^- \times M_2 \to M_0^- \times M_2$ to $L_{01} \times_{M_1} L_{12}$
is automatically a Lagrangian immersion.
We will study the class of {\em embedded} geometric compositions, for which in addition $\pi_{02}$ is injective, and hence $L_{01} \circ L_{12}$ is a smooth Lagrangian correspondence.
If the composition $L_{(\ell-1)\ell}\circ L_{\ell(\ell+1)}$ is embedded, then we obtain under suitable monotonicity assumptions a canonical isomorphism
\begin{equation}\label{eq:iso}
HF(\ldots, L_{(\ell-1)\ell}, L_{\ell(\ell+1)}, \ldots)
\cong
HF(\ldots, L_{(\ell-1)\ell}\circ L_{\ell(\ell+1)}, \ldots) .
\end{equation}
For the precise monotonicity and admissibility conditions
see Section \ref{shrink}. The proof proceeds in two steps.
First, we allow for varying widths $(\delta_j>0)_{j=0,\ldots,r}$ of the pseudoholomorphic strips $(u_j:\R\times[0,\delta_j]\to M_j)_{j=0,\ldots,r}$ defining the differential.
Section~\ref{sec:invariance} of this paper shows that Floer cohomology is independent of the choice of widths.
(These domains are not conformally equivalent due to the identification
between boundary components that is implicit in the seam conditions.) 
The second (hard analytic) part is to prove that with the width $\delta_\ell>0$ sufficiently close to zero, the $r+1$-tuples of holomorphic strips with seam conditions in
$(\ldots, L_{(\ell-1)\ell}, L_{\ell(\ell+1)}, \ldots)$ are in one-to-one correspondence with the $r+1$-tuples of holomorphic strips with seam conditions in $(\ldots, L_{(\ell-1)\ell}\circ L_{\ell(\ell+1)}, \ldots)$.
This analysis is completely analogous to \cite{isom}, where we establish the bijection for the Floer trajectories of the special cyclic sequence
$\pt \overset{L_{0}}{\longrightarrow}  M_0 \overset{L_{01}}{\longrightarrow} M_1\overset{L_{12}}{\longrightarrow} M_2 \overset{L_{2}}{\longrightarrow} \pt $
when $\delta_1\to 0$.
The monotonicity assumptions are crucial for this part since the exclusion of a novel ``figure eight bubble" in \cite{isom} hinges on a strict energy-index proportionality.

In section~\ref{sec:app} we provide a number of new tools for the calculation of Floer cohomology (and hence detection of non-displaceability), arising as direct consequences of \eqref{eq:iso} or from a conjectural generalization of Perutz' long exact Gysin sequence \cite{per:gys}. These are meant to exemplify the wide applicabilty and reach of the basic isomorphism \eqref{eq:iso}. We believe that it should have much more dramatic consequences once systematically employed. 
As first specific example of direct consequences of \eqref{eq:iso}
we confirm the calculation $HF(T^n_{\rm Cl},T^n_{\rm Cl})\cong H_*(T^n)$  of Cho \cite{cho:hol} for the Clifford torus in $\CP^n$, and we calculate some further Floer cohomologies in $\CP^n$ using reduction at pairs of transverse level sets.
Next, we prove Hamiltonian non-displaceability of the Lagrangian $3$-sphere $\Sigma \subset (\CP^1)^-\times \CP^2$ arising from reduction at the level set of an $S^1$-action on $\CP^2$ containing $T_{\rm Cl}$. The latter will be deduced from the nontriviality of $HF(T^2_{\rm Cl},T^2_{\rm Cl})$ together with our isomorphism
$$
HF(S^1 \times T_{\rm Cl}, \Sigma ) \cong HF(T_{\rm Cl},T_{\rm Cl}) ,
$$
and the fact that $HF(L',L)\neq 0 \Rightarrow HF(L,L)\neq 0$ (since $HF(L,L)$ contains an element that acts as the identity on $HF(L',L)$).
Finally, we generalize this non-displaceability result to the Lagrangian embedding
$\Sigma \subset (\CP^{k-1})^-\times \CP^n$ of $(S^1)^{n-k} \times S^{2k-1}$, 
which arises from the monotone level set of an $S^{n-k+1}$-action on $\CP^n$
for $2\leq k \leq n$.

Another consequence of our results is a general prescription for defining topological invariants by
decomposing into elementary pieces. For example, let $Y$ be a compact
manifold and $f: Y \to \R$ a Morse function, which induces a decomposition $Y
= Y_{01} \cup \ldots \cup Y_{(k-1)k}$ into elementary cobordisms by
cutting along non-critical level sets $X_1,\ldots, X_{k-1}$.  First
one associates to each $X_j$ a monotone symplectic manifold $M(X_j)$,
and to each $Y_{(j-1)j}$ with $\partial Y_{(j-1)j}=X_{j-1}^-\sqcup
X_j$ a smooth monotone Lagrangian correspondence $L(Y_{(j-1)j})\subset
M(X_{j-1})^-\times M(X_j)$ (taking $M(X_0)$ and $M(X_k)$ to be
points.)  Second, one checks that the basic moves described by Cerf
theory (cancellation or change of order of critical points) 
change the sequence of Lagrangian correspondences by
replacing adjacent correspondences with an embedded composition, or
vice-versa.  In other words, the equivalence class of sequences of
Lagrangian correspondences by embedded compositions
$[L(Y_{01}),\ldots,L(Y_{(k-1)k})]$ does not depend on the choice of
the Morse function $f$.  
Then the results of this paper provide a
group-valued invariant of $Y$, by taking the Floer homology of the
sequence of Lagrangian correspondences.
For example, in \cite{field, fieldb} we investigate the theory which uses as
symplectic manifolds the moduli spaces of flat bundles with compact
structure group on three-dimensional cobordisms containing tangles.
Conjecturally this provides a symplectic construction of Donaldson type gauge theoretic invariants: $SO(3)$-instanton Floer homology, it's higher rank version (though not strictly defined in the gauge theoretic setting), and the $SU(n)$-tangle invariants defined by Kronheimer-Mrowka \cite{KMknot} from singular instantons.
The same setup is used to give alternative constructions of Heegard Floer homology \cite{yanki} and Seidel-Smith invariants \cite{reza}.

Even more generally, in \cite{cat} the quilted Floer cohomology groups provide the $2$-morphism spaces in our construction of a symplectic $2$-category which contains all Lagrangian correspondences as morphisms, and where the composition is equivalent to geometric composition, if the latter is embedded.
This setup in turn is used by Abouzaid-Smith \cite{as:hms} to prove mirror symmetry for the $4$-torus and deduce various further symplectic consequences.
\\

\noindent
{\bf Notation and Organization:}
We will frequently refer to the assumptions (M1-2), (L1-3),
and (G1-2) that can be found on pages \pageref{M ass} -
%, \pageref{L ass}, 
\pageref{G ass}.

Section~\ref{LC} is a detailed introduction to Lagrangian correspondences, geometric composition, and sequences of correspondences, which also provides the basic framework for the sequels \cite{cat,quilts}  to this paper.
In Section~\ref{gradings} we generalize gradings to Lagrangian correspondences and establish their behaviour under geometric composition, so that the isomorphism \eqref{eq:iso} becomes an isomorphism of graded groups.
Section~\ref{sec:HF} provides a review of monotonicity and Floer cohomology and gives a first definition of the Floer cohomology \eqref{eq:HF} by building a pair of Lagrangians in the product $M_0\times M_1 \times \ldots M_{k-1}$. 
The latter is however unsatisfactory since it does not provide an approach to the isomorphism \eqref{eq:iso}.
Section~\ref{quilted FH} gives the general definition of quilted Floer cohomology \eqref{eq:HF} and finalizes the proof of the isomorphism \eqref{eq:iso}. 
Finally, Section~\ref{sec:app} gives a number of direct symplectic applications of the isomorphism \eqref{eq:iso}.

\medskip

{\em We thank Paul Seidel and Ivan Smith for encouragement and helpful
discussions, and we are very grateful to the meticulous referee for help with cleaning up the details.}  

\section{Lagrangian correspondences}
\label{LC}

Let $M$ be a smooth manifold.  A {\em symplectic form} on $M$ is a
closed, non-degenerate two-form $\omega$.  A {\em symplectic manifold}
is a smooth manifold equipped with a symplectic form.  If
$(M_1,\omega_1)$ and $(M_2,\omega_2)$ are symplectic manifolds, then a
diffeomorphism ${\varphi: M_1 \to M_2}$ is a {\em symplectomorphism} if
$\varphi^* \omega_2 = \omega_1$.  Let $\Symp$ denote the category
whose objects are symplectic manifolds and whose morphisms are
symplectomorphisms.  We have the following natural operations on $\Symp$.

\begin{enumerate}
\item (Duals)  If $M = (M,\omega)$ is a symplectic manifold, then $M^- =
(M,-\omega)$ is a symplectic manifold, called the {\em dual} of $M$.
\item (Sums)  If $M_j = (M_j,\omega_j), j = 1,2$ are symplectic
  manifolds, then the disjoint union $M_1 \cup M_2$ equipped with the symplectic 
structure $\omega_1$ on $M_1$ and $\omega_2$ on $M_2$, is a symplectic manifold. 
The empty set $\emptyset$ is a unit for the disjoint union.

\item (Products) Let $M_j = (M_j,\omega_j), j = 1,2$ be symplectic manifolds, 
then the Cartesian product
$(M_1 \times M_2, \pi_1^*\omega_1 + \pi_2^* \omega_2)$ is a symplectic manifold.
(Here $\pi_j: M_1 \times M_2 \to M_j$ denotes the projections.)
The symplectic manifold $\pt$, consisting of a single point, is a unit for the 
Cartesian product.
\end{enumerate}

Clearly the notion of symplectomorphism is very restrictive; in
particular, the symplectic manifolds must be of the same dimension.  A
more flexible notion of morphism is that of Lagrangian correspondence,
defined as follows \cite{we:le,we:sc,gu:rev}.  Let $M = (M,\omega)$ be
a symplectic manifold.  A submanifold 
$L \subset M$ is {\em isotropic}, resp.\ {\em coisotropic}, resp.\ {\em Lagrangian} if the
$\omega$-orthogonal complement $TL^\omega$ satisfies 
$TL^\omega \supseteq TL$ resp.\ $TL^\omega \subseteq TL$ resp.\
$TL^\omega = TL$.

\begin{definition}
Let $M_1, M_2$ be symplectic manifolds.  
A {\em Lagrangian correspondence from $M_1$ to $M_2$}
is a Lagrangian submanifold $L_{12}\subset M_1^- \times M_2$.  
\end{definition}

\begin{example} \label{ex:corr}
The following are examples of Lagrangian correspondences:
\begin{enumerate}
\item (Trivial correspondence) The one and only Lagrangian correspondence between
$M_1=\emptyset$ and any other $M_2$ is $L_{12}=\emptyset$.
\item (Lagrangians) Any Lagrangian submanifold $L\subset M$ can be viewed both as correspondence $L\subset \pt^-\times M$ from the point to $M$ and as correspondence $L\subset M^-\times \pt$ from $M$ to the point.
\item (Graphs) If $\varphi_{12}:M_1 \to M_2$ is a symplectomorphism then its
graph
$$ \graph(\varphi_{12}) = \{ (m_1,\varphi_{12}(m_1)) \,|\, \ m_1 \in M_1 \}
\subset M_1^- \times M_2 $$
is a Lagrangian correspondence.  
\item (Fibered coisotropics) Suppose that $\iota: C \to M$ is a coisotropic submanifold.
Then the null distribution $TC^\omega$ is integrable, see e.g.\ \cite[Lemma 5.30]{ms:jh}.
Suppose that $TC^\omega$ is in fact fibrating, that is, there exists
a symplectic manifold $(B,\omega_B)$ and a fibration $\pi: C \to B$ such
that $\iota^* \omega$ is the pull-back $\pi^* \omega_B$.  Then
$$ (\iota \times \pi): C \to M^- \times B $$
maps $C$ to a Lagrangian correspondence.
\item (Level sets of moment maps) Let $G$ be a Lie group with Lie
algebra $\g$.  Suppose that $G$ acts on $M$ 
by Hamiltonian symplectomorphisms generated by a moment map
$\mu: M \to \g^*$. (That is $\mu$ is equivariant and the
generating vector fields $\g \to \Vect(M), \xi \mapsto \xi_M$
satisfy $\iota(\xi_M) \omega = - d(\mu,\xi)$.)
If $G$ acts freely on $\muinv(0)$, then $\muinv(0)$ is a smooth
coisotropic fibered over the symplectic quotient
$ M \qu G = \muinv(0)/G$, which is a symplectic manifold.  
Hence we have a Lagrangian correspondence
$$ (\iota \times \pi): \ \muinv(0) \to M^- \times (M \qu G) .$$
The symplectic two-form $\omega_{M \qu G}$ on $M \qu G$ is the unique
form on $M \qu G$ satisfying $ \pi^* \omega_{M \qu G} = \iota^* \omega $.
\een
\end{example}

\begin{definition}
Let $M_0,M_1,M_2$ be symplectic manifolds and
${L_{01} \subset M_0^- \times M_1}$,  $L_{12} \subset M_1^- \times M_2$
Lagrangian correspondences.  
\begin{enumerate}
\item The {\em dual} Lagrangian correspondence of $L_{01}$ is
$$(L_{01})^t := \{ (m_1,m_0) \,|\, (m_0,m_1) \in L_{01} \} 
 \subset M_1^- \times M_0.$$
\item
The {\em geometric composition} of $L_{01}$ and $L_{12}$ is 
$$ 
L_{01} \circ L_{12} :=
\left\{ (m_0,m_2) \in M_0^- \times M_2  \,\left|\, \exists m_1
\in M_1 : \begin{aligned}
&(m_0,m_1) \in L_{01}\\
&(m_1,m_2) \in L_{12} \end{aligned}
\right.\right\} \subset M_0^-\times M_2 .
$$
\end{enumerate}
\end{definition}

Geometric composition and duals of Lagrangian correspondences satisfy the following:
\begin{enumerate}
\item (Graphs) 
If $\varphi_{01}: M_0 \to M_1$ and
$\varphi_{12}: M_1 \to M_2$ are symplectomorphisms, then
$$ \graph(\varphi_{01}) \circ \graph(\varphi_{12}) =
\graph(\varphi_{12} \circ \varphi_{01}), $$
$$\graph(\varphi_{01})^t = \graph(\varphi_{01}^{-1}) .$$
\item(Zero) Composition with $\emptyset$ always yields $\emptyset$, that is for any
Lagrangian correspondence  $L_{01} \subset M_0^- \times M_1$ we have
$$  \emptyset \circ L_{01} =  \emptyset, \qquad L_{01} \circ \emptyset = \emptyset .$$
\item (Identity)  If $L_{01} \subset M_0^- \times M_1$ is 
a Lagrangian correspondence and $\Delta_j \subset M_j^- \times M_j, j =0,1$
are the diagonals, then 
$$ L_{01} = \Delta_0 \circ L_{01} = L_{01} \circ \Delta_1 .$$
\item (Associativity) If $ L_{01} \subset M_0^- \times M_1, L_{12}
\subset M_1^- \times M_2, L_{23} \subset M_2^- \times M_3$ are
Lagrangian correspondences, then
$$ (L_{01} \circ L_{12}) \circ L_{23} = L_{01} \circ (L_{12} \circ
L_{23}) , $$
$$ (L_{01} \circ L_{12})^t = (L_{12})^t \circ (L_{01})^t . $$
\een

The geometric composition can equivalently be defined as 
$ L_{01} \circ L_{12} = \pi_{02}( L_{01} \times_{M_1}
L_{12}) $, the image under the projection $ \pi_{02}: M_0^- \times
M_1 \times M_1^- \times M_2 \to M_0^- \times M_2 $ of
$$ 
L_{12} \times_{M_1} L_{01} := (L_{01} \times L_{12})
 \cap (M_0^- \times \Delta_1 \times M_2) .
$$ Here $\Delta_1 \subset M_1^- \times M_1$ denotes the diagonal.
$L_{01} \circ L_{12}\subset M_0^-\times M_2$ is an immersed Lagrangian
submanifold if $L_{01}\times L_{12}$ intersects $M_0^- \times \Delta_1
\times M_2$ transversally.  In general, the geometric composition of smooth
Lagrangian submanifolds may not even be immersed. We will be working
with the following class of compositions, for which the resulting
Lagrangian correspondence is in fact a smooth submanifold, as will be seen in Lemma~\ref{immersion} below.

\begin{definition} \label{embedded}
We say that the composition $ L_{01} \circ L_{12}$ is {\em embedded}
if $L_{12} \times_{M_1} L_{01}$ is cut out transversally (i.e.\  
$(L_{01} \times L_{12}) \pitchfork (M_0^- \times \Delta_1 \times M_2)$)
and the projection $\pi_{02}:L_{12} \times_{M_1} L_{01}\to L_{01}\circ L_{12}\subset
M_0^-\times M_2$ is an embedding.
(For compact Lagrangians it suffices to assume that $\pi_{02}$ is injective, by Lemma~\ref{immersion} below.)
\end{definition}

By some authors (e.g.\  \cite{gu:rev}) geometric composition of Lagrangian correspondences is more generally defined under clean intersection hypotheses.  This extension is not needed in the present paper, because the quilted Floer cohomology is invariant under Hamiltonian isotopy, and after such an isotopy transversality may always be achieved. However, transverse intersection only yields an immersed\footnote{One can not necessarily remove all self-intersections of the immersed composition by Hamiltonian isotopy on one correspondence. A basic example is the
composition of transverse Lagrangian submanifolds $L,L'\subset M$. Identifying 
$M\cong M\times\{\pt\}\cong \{\pt\}\times M$ the projection 
$L\times_M L'\to L\circ L'\subset \{\pt\}\times\{\pt\}$ maps the (finite) intersection 
$L\pitchfork L'$ to a point.
}
Lagrangian correspondence, as the following Lemma shows.

\begin{lemma} \label{immersion}
Let ${L_{01} \subset M_0^- \times M_1}$,  $L_{12} \subset M_1^- \times M_2$ be Lagrangian correspondences such that the intersection
$(L_{01} \times L_{12}) \pitchfork (M_0^- \times \Delta_1 \times M_2)$ is transverse.
Then the projection $\pi_{02}:L_{12} \times_{M_1} L_{01}\to L_{01}\circ L_{12}\subset
M_0^-\times M_2$ is an immersion.

In particular, if the Lagrangians are compact, the intersection is transverse, and the projection is injective, then the composition $L_{01} \circ L_{12}=:L_{02}$ is embedded.
\end{lemma}
\begin{proof}
The proof essentially is a special case of the fact that the geometric composition of linear Lagrangian correspondences is always well defined (i.e.\ yields another linear Lagrangian correspondence), see e.g.\ \cite[Section 4.1]{gu:rev}.

Fix a point $\ul{x}=(x_0,x_1,x_1,x_2)\in  L_{01} \times_{M_1} L_{12}$ then we need to check that
$\ker {\rm d}_{\ul{x}} \pi_{02}=\{0\}$ for the projection restricted to 
$L_{12} \times_{M_1} L_{01}$. In fact, we will show that
\beq \label{imm claim}
\ker {\rm d}_{\ul{x}} \pi_{02} \cong \frac{ T_{\ul{x}}(M_0\times M_1\times M_1\times M_2) }
{( T_{(x_0,x_1)} L_{01} \times T_{(x_1,x_2)} L_{12} ) + ( T_{x_0} M_0 \times T_{(x_1,x_1)} \Delta_1 \times T_{x_2} M_2 ) },
\eeq
which is zero by transversality.
To simplify notation we abbreviate $\Lambda_{01}:=T_{(x_0,x_1)} L_{01}$, $\Lambda_{12}:=T_{(x_1,x_2)} L_{12}$, and $V_i:= T_{x_i} M_i$. Now \eqref{imm claim} follows as in \cite[Section 4.1]{gu:rev}. For completeness we recall the precise argument:
We identify
\begin{align}
\ker {\rm d}_{\ul{x}} \pi_{02} \;
&= ( \Lambda_{01}\times_{V_1} \Lambda_{12} ) \cap ( \{0\} \times V_1 \times V_1 \times \{0\} ) 
\nonumber \\
&\cong \bigl\{ v_1 \in V_1 \,\big|\, (0,v_1)\in\Lambda_{01}, (v_1,0)\in\Lambda_{12} \bigr\} 
\;=\; \ker\Lambda_{01}^t \cap \ker\Lambda_{12} ,
\label{first} 
\end{align}
where $\ker\Lambda_{12}:=\{ v_1 \in V_1 \,|\, (v_1,0)\in\Lambda_{12} \}\subset V_1$ 
and similarly $\ker\Lambda_{01}^t\subset V_1$.
On the other hand, we use the symplectic complements with respect to 
$\omega_{0112}:=(- \omega_0) \oplus \omega_1 \oplus (- \omega_1) \oplus \omega_2$
on $V_0\times V_1\times V_1 \times V_2$ to identify
\begin{align}
&(V_0\times V_1\times V_1\times V_2) \; / \; 
(\Lambda_{01} \times \Lambda_{12} ) + ( V_0 \times \Delta_{V_1} \times V_2 ) \nonumber \\
&\quad\cong 
(\Lambda_{01} \times \Lambda_{12} )^{\omega_{0112}} \cap ( V_0 \times \Delta_{V_1} \times V_2 )^{\omega_{0112}}  \nonumber \\
&\quad= \bigl\{(0,v_1,v_1, 0) \,\big|\, (0,v_1)\in\Lambda_{01}^{(-\omega_0)\oplus\omega_1}, (v_1,0)\in\Lambda_{12}^{(-\omega_1)\oplus\omega_2} \bigr\} \nonumber \\
&\quad \cong ( \im\Lambda_{01} )^{\omega_1} \cap ( \im\Lambda_{12}^t )^{\omega_1} , 
\label{second} 
\end{align}
where $\im\Lambda_{01}:=\pi_{V_1}(\Lambda_{01})\subset V_1$, similarly $\im\Lambda_{12}^t\subset V_1$, and we used the equivalence 
$$
(-\omega_0 \oplus \omega_1) ((0,v_1), (\lambda_0,\lambda_1) ) = 0 
\;\;\forall  (\lambda_0,\lambda_1)\in\Lambda_{01}
\quad\Leftrightarrow\quad
\omega_1(v_1,\lambda_1 ) = 0 \;\;\forall  \lambda_1\in\pi_{V_1}(\Lambda_{01}) .
$$
Now the two vector spaces in \eqref{first} and \eqref{second} are identified by the dualities 
$\ker\Lambda_{01}^t  = (\im\Lambda_{01})^{\omega_1}$ and
$ \ker\Lambda_{12} = (\im\Lambda_{12}^t )^{\omega_1} $,
which follow from the Lagrangian property of $\Lambda_{01}$ resp.\ $\Lambda_{12}^t$,
$$
(0,v_1)\in \Lambda_{01}
\quad\Leftrightarrow\quad
(0,v_1)\in \Lambda_{01}^{-\omega_0\oplus\omega_1} 
\quad\Leftrightarrow\quad
v_1\in (\im\Lambda_{01})^{\omega_1} .
$$
This proves \eqref{imm claim} and hence finishes the proof that $\pi_{02}$ is an immersion.
Finally, note that an injective immersion of a compact set is automatically an embedding.
\end{proof}

\begin{remark} \label{rmk:embedded}
Suppose that the composition $L_{01} \circ L_{12}=:L_{02}$ is embedded.
\begin{enumerate}
\item
By the injectivity, for every $(x_0,x_2)\in L_{02}$ there is a
unique solution $x_1\in M_1$ to $(x_0,x_1,x_1,x_2)\in L_{01} \times
L_{12}$.  Due to the transversality assumption, this solution is given
by a smooth map $\ell_1 : L_{02} \to M_1$.
\item
If $L_{01}$ and $L_{12}$ are compact, oriented, and equipped with a
relative spin structure, then $L_{02}$ is also compact and inherits an
orientation and relative spin structure, see \cite{orient}.
The orientation is induced from the canonical orientation of the diagonal, see Remark~\ref{Maslov diagonal}(b), and the splitting
$$
T(M_0\times M_2\times M_1\times M_1) 
= \bigl( TL_{02}\times\{0\}\bigr) \oplus \bigl(\{0\} \times (T\Delta_1)^\perp \bigr)
\oplus T(L_{01} \times L_{12})^\perp .
$$

\item 
If both fundamental groups $\pi_1(L_{01})$ and $\pi_1(L_{12})$ have torsion
image in the respective ambient space, then
$\pi_1(L_{02})$ has torsion image in $\pi_1(M_0\times M_2)$.  
(Any loop $\gamma:S^1\to L_{02}$ lifts to a loop 
$\ti\gamma:S^1\to L_{01}\times L_{12}$ with $\gamma=\pi_{02}\circ\ti\gamma$. By assumption, some multiple cover $\ti\gamma^N$ is the boundary of a map 
$\ti u:D^2\to M_0\times M_1\times M_1\times M_2$. 
Hence the same cover $\gamma^N = (\pi_{02}\circ\ti u)|_{\partial D^2}$ is contractible.)
\end{enumerate}
\end{remark}

\subsection{Generalized Lagrangian correspondences}
\label{symp cat}

A simple resolution of the composition problem (that geometric composition is not always defined) is
given by passing to sequences of Lagrangian correspondences and composing them by concatenation.
In \cite{cat} we employ these to define a symplectic category containing all smooth Lagrangian correspondences as composable morphisms, yet retaining geometric composition in cases where it is well defined.

\begin{definition} \label{Lag cor}
Let $M,M'$ be symplectic manifolds.  A {\em generalized Lagrangian
correspondence} $\ul{L}$ from $M$ to $M'$ consists of
\begin{enumerate}
\item a sequence $N_0,\ldots,N_r$ of any length $r+1\geq 2$ of
symplectic manifolds with $N_0 = M$ and $N_r = M'$ ,
\item a sequence $L_{01},\ldots, L_{(r-1)r}$ of compact Lagrangian
correspondences with $L_{(j-1)j} \subset N_{j-1}^-\times N_{j}$ for
$j=1,\ldots,r$.
\end{enumerate}
\end{definition}
\begin{definition}
Let $\ul{L}$ from $M$ to $M'$ and $\ul{L}'$ from $M'$ to $M''$ 
be two generalized Lagrangian correspondences. Then we define composition
$$
(\ul{L},\ul{L}') :=
\bigl(L_{01},\ldots,L_{(r-1)r},L'_{01},\ldots,L'_{(r'-1)r'}\bigr) 
$$
as a generalized Lagrangian correspondence from $M$ to $M''$.
Moreover, we define the dual
$$
\ul{L}^t := \bigl(L_{(r-1)r}^t,\ldots,L_{01}^t\bigr) .
$$
as a generalized Lagrangian correspondence from $M'$ to $M$.
\end{definition}

We conclude this subsection by mentioning special cases of
generalized Lagrangian correspondences. The first is the case $M=M'$, 
which we will want to view separately as a cyclic correspondence, without
fixing the ``base point'' $M$.

\begin{definition} \label{Lag cycle} 
A {\em cyclic generalized Lagrangian correspondence} $\ul{L}$ consists of
\begin{enumerate}
\item a cyclic sequence $M_0,M_1,\ldots,M_r,M_{r+1}=M_0$ of symplectic manifolds 
of any length $r+1\geq 1$, 
\item a sequence $L_{01},\ldots, L_{r(r+1)}$ of compact Lagrangian correspondences 
with $L_{j(j+1)} \subset M_{j}^-\times M_{j+1}$ for $j=0,\ldots,r$.
\end{enumerate}
\end{definition}

The second special case is $M=\{pt\}$, which generalizes the concept of Lagrangian
submanifolds. Namely, note that any Lagrangian submanifold $L\subset M'$ can be
viewed as correspondence $L\subset\{pt\}^-\times M'$.

\begin{definition} \label{genLag}
Let $M'$ be a symplectic manifold.
A {\em generalized Lagrangian submanifold} $\ul{L}$ of $M'$ 
is a generalized Lagrangian correspondence from a point $M=\{pt\}$ to $M'$.
That is, $\ul{L}$ consists of 
\begin{enumerate}
\item a sequence $N_{-r},\ldots,N_0$ of any length $r\geq 0$ 
of symplectic manifolds with $N_{-r}=\{\pt\}$ a point and $N_0 = M'$,
\item a sequence $L_{(-r)(-r+1)},\ldots,L_{(-1)0}$ of compact 
Lagrangian correspondences $L_{(i-1)i}\subset N_{i-1}^- \times N_i$.
\end{enumerate}
\end{definition}

\section{Gradings} 
\label{FH pairs}
\label{gradings}

The purpose of this section is to review the theory of graded
Lagrangians and extend it to generalized Lagrangian
correspondences. It can be skipped at first reading.  
% CW: out 
%When considering
%ungraded Floer homologies only the index identity in Lemma~\ref{index}
%is relevant.

Following Kontsevich and Seidel \cite{se:gr} one can define graded Lagrangian
subspaces as follows.  Let $V$ be a symplectic vector space and 
let $\Lag(V)$ be the Lagrangian Grassmannian of $V$.  An {\em $N$-fold Maslov
covering for $V$} is a $\Z_N$-covering $\Lag^N(V) \to \Lag(V)$
associated to the Maslov class in ${\rm Hom}(\pi_1(\Lag(V)),\Z)$.  (More precisely, the mod $N$ reduction of the Maslov class defines a representation $\pi_1(\Lag(V))\to\Z_N$. 
Let $\widetilde{\Lag}(V)$ be the universal cover of $ \Lag(V)$, then the $N$-fold Maslov
covering associated to the given representation is the associated bundle
$\widetilde{\Lag}(V)\times_{\pi_1(\Lag(V))} \Z_N \to  \Lag(V)$.)
A {\em grading} of a Lagrangian subspace $\Lambda\in\Lag(V)$ is a lift to $\ti{\Lambda}\in\Lag^N(V)$.

\begin{remark} \label{Maslov diagonal} 
\begin{enumerate}
\item
For any basepoint $\Lambda_0\in\Lag(V)$ we obtain an $N$-fold Maslov
cover $\Lag^N(V,\Lambda_0)$ given as the homotopy classes of paths
$\ti{\Lambda}:[0,1]\to\Lag(V)$ with base point
$\ti{\Lambda}(0)=\Lambda_0$, modulo loops whose Maslov index is a multiple of $N$. The
covering is $\ti{\Lambda}\mapsto\ti{\Lambda}(1)$.  The base point has
a canonical grading given by the constant path
$\tilde\Lambda_0\equiv\Lambda_0$.  Any path between basepoints
$\Lambda_0,\Lambda_0'$ induces an identification $\Lag^N(V,\Lambda_0)
\to \Lag^N(V,\Lambda_0')$.
\item
For the diagonal $\Delta \subset V^-\times V$ we fix a canonical
grading and orientation as follows. 
(This choice is made in order to obtain the degree identity in Lemma~\ref{degprop}~(d).)
 We identify the Maslov coverings
$\Lag^N(V^- \times V,\Lambda^- \times \Lambda)$ and $\Lag^N(V^- \times
V, \Delta)$ by concatenation of the paths\footnote{
The first path arises from the canonical path between ${\rm Id}$ and $J$ in ${\rm Symp}(V)$. The second path can be understood as the graphs of the symplectomorphisms $(t^{-1}{\rm Id}) \times (t{\rm Id})$ on $V\cong \Lambda \times J\Lambda$. For $t\to 0$ this graph converges to the split Lagrangian $J\Lambda \times \Lambda \subset V^-\times V$; for $t=1$ the graph is the diagonal.
}
\begin{equation} \label{diag grading} ( e^{Jt}\Lambda^-\times \Lambda)_{t\in[0,\pi/2]}, \ \ \ 
(\{ (tx+Jy,x+tJy) | x,y\in \Lambda \})_{t\in[0,1]} , \end{equation}
where $J\in\End(V)$ is an $\omega$-compatible complex structure on $V$
(i.e.\ $J^2=-\Id$ and $\omega(\cdot,J\cdot)$ is symmetric and positive
definite).  In particular, this induces the {\em canonical grading} on
the diagonal $\Delta$ with respect to any Maslov covering $\Lag^N(V^-
\times V, \Lambda^- \times \Lambda)$, by continuation.  Any
identification $\Lag^N(V^- \times V, \Lambda^-_0 \times \Lambda_0) \to
\Lag^N(V^- \times V, \Lambda^-_1 \times \Lambda_1)$ induced by a path
in $\Lag^N(V)$ maps the graded diagonal to the graded diagonal, since
the product $\gamma^-\times\gamma$ of any loop $\gamma:S^1\to\Lag(V)$
has Maslov index $0$.  Similarly, we define a {\em canonical
orientation} on $\Delta$ by choosing any orientation on $\Lambda$,
giving the product $\Lambda^- \times \Lambda$ the product orientation
(which is well defined), and extending the orientation over the path
\eqref{diag grading}.  This is related to the orientation induced by
projection of the diagonal on the second factor by a sign
$(-1)^{n(n-1)/2}$, where $\dim(M) = 2n$.
\end{enumerate}
\end{remark}

Let $M$ be a symplectic manifold and let $\Lag(M) \to M$ be the fiber
bundle whose fiber over $m \in M$ is the space $\Lag(T_mM)$ of
Lagrangian subspaces of $T_m M$.  An {\em $N$-fold Maslov covering} of
$M$ is an $N$-fold cover $\Lag^N(M) \to \Lag(M)$ whose restriction to
each fiber is an $N$-fold Maslov covering $\Lag^N(T_mM) \to \Lag(T_m
M)$.  Any choice of Maslov cover for $\R^{2n}$ induces a one-to-one
correspondence between $N$-fold Maslov covers of $M$ and
$\Sp^N(2n)$-structures on $M$.  Here $2n=\dim M$ and $\Sp^N(2n)$ is
the $N$-fold covering group of $\Sp(2n)$ associated to the mod $N$ reduction of the Maslov class in ${\rm Hom}(\pi_1(\Sp(2n)),\Z_N)$.  (Explicitly, this is realized by using the
identity as base point.)  An {\em $\Sp^N(2n)$-structure} on $M$ is an
$\Sp^N(2n)$-bundle $\Fr^N(M)\to M$ together with an isomorphism
$\Fr^N(M)\times_{\Sp^N(2n)} \Sp(2n) \simeq \Fr(M)$ to the symplectic
frame bundle of~$M$.  It induces the $N$-fold Maslov covering
$$ \Lag^N(M) = \Fr^N(M)\times_{\Sp^N(2n)} \Lag^N(\R^{2n}). $$
The notions of duals, disjoint union, and Cartesian product extend
naturally to the graded setting as follows.
The dual $\Lag^N(M^-)$ of a Maslov covering $\Lag^N(M) \to \Lag(M)$ is the 
same space with the inverted $\Z_N$-action. 
We denote this identification by 
\begin{equation}\label{dual g}
\Lag^N(M) \to \Lag^N(M^-), \qquad
\ti{\Lambda}\mapsto\ti{\Lambda}^-.
\end{equation}
For $\Sp^N$-structures $\Fr^N(M_0)$ and $\Fr^N(M_1)$ the embedding
$$ \Sp^N(2n_0) \times_{\Z_N} \Sp^N(2n_1) \to \Sp^N(2n_0 + 2n_1) $$
induces an $\Sp^N(2n_0 + 2n_1)$-structure $\Fr^N(M_0 \times M_1)$ on the product
and an equivariant map
\begin{equation} \label{fiberfr}
\Fr^N(M_0) \times \Fr^N(M_1) \to \Fr^N(M_0 \times M_1)
\end{equation}
covering the inclusion $\Fr(M_0) \times \Fr(M_1) \to \Fr(M_0 \times M_1)$.
The corresponding product of $N$-fold Maslov covers on $M_0 \times M_1$ is the 
$N$-fold Maslov covering
$$
\Lag^N(M_0 \times M_1):=
\bigl( \Fr^N(M_0)\times\Fr^N(M_1) \bigr) \times_{\Sp^N(2n_0)\times\Sp^N(2n_1)}
\Lag^N(\R^{2n_0}\times\R^{2n_1}) .
$$
Combining this product with the dual yields a Maslov covering for $M_0^-\times M_1$
which we can identify with
$$
\Lag^N(M_0^- \times M_1) =
\bigl( \Fr^N(M_0)\times\Fr^N(M_1) \bigr) \times_{\Sp^N(2n_0)\times\Sp^N(2n_1)}
\Lag^N(\R^{2n_0,-}\times\R^{2n_1}) .
$$
Finally, the inclusion $\Lag(M_0) \times \Lag(M_1) \to \Lag(M_0 \times M_1)$ lifts to a map
\begin{equation} \label{fiberzn} 
\Lag^N(M_0) \times \Lag^N(M_1) \to \Lag^N(M_0 \times M_1), \ \ \ 
(\ti{L}_0,\ti{L}_1) \mapsto \ti{L}_0 \times^N \ti{L}_1 
\end{equation}
with fiber $\Z_N$.  It is defined by combining the product
\eqref{fiberfr} with the basic product of the linear Maslov cover
$\Lag^N(\R^{2n_0}) \times \Lag^N(\R^{2n_1}) \to \Lag^N(\R^{2n_0}\times
\R^{2n_1})$.

\begin{definition} \label{def:graded symp, Lag}
\begin{enumerate}
\item
Let $M_0$, $M_1$ be two symplectic manifolds equipped with $N$-fold Maslov covers and 
let $\phi: M_0 \to M_1$ be a symplectomorphisms.  A {\em grading} of $\phi$ is a
lift of the canonical isomorphism $\Lag(M_0) \to \Lag(M_1)$ to an
isomorphism $\phi^N:\Lag^N(M_0) \to \Lag^N(M_1)$, or equivalently, a
lift of the canonical isomorphism $\Fr(M_0) \to \Fr(M_1)$ of symplectic frame bundles to 
an isomorphism $\Fr^N(M_0) \to \Fr^N(M_1)$. 
\item
Let $L \subset M$ be a Lagrangian submanifold and 
$M$ be equipped with an $N$-fold Maslov cover.  
A {\em grading} of $L$ is a lift $\sigma_L^N : L \to \Lag^N(M)$
of the canonical section $\sigma_L: L \to \Lag(M)$. 
\end{enumerate}
\end{definition}

\begin{remark} \label{rmk:gradings}
\begin{enumerate}
\item
The set of graded symplectomorphisms forms a group under composition.
In particular, the identity on $M$ has a canonical grading, given by the
identity on $\Lag^N(M)$. 
\item 
Given a one-parameter family $\phi_t$ of symplectomorphisms 
with $\phi_0 = \Id_M$, we obtain a grading of $\phi_t$ by continuity.  
\item Any choice of grading on the diagonal
$\ti{\Delta}\in\Lag^N(\R^{2n,-}\times\R^{2n})$ induces a 
bijection between gradings of a symplectomorphism $\phi:M_0 \to M_1$
and gradings of its graph $\graph(\phi) \subset M_0^- \times M_1$ with
respect to the induced Maslov cover $\Lag^N(M_0^-\times M_1)$.
Indeed, the graph of the grading, $\graph(\phi^N) \subset ( \Fr^N(M_0)
\times \Fr^N(M_1)) |_{\graph(\phi)}$ is a principal bundle over
$\graph(\phi)$ with structure group $\Sp^N(2n)$, $2n=\dim M_0=\dim
M_1$.  The graded diagonal descends under the associated fiber bundle
construction with $\graph(\phi^N)$ to a section of ${\Lag^N(M_0^-
\times M_1)|_{\graph(\phi)}}$ lifting $\graph(\phi)$.  Moreover, this
construction is equivariant for the transitive action of
$H^0(M_0,\Z_N)$ on both the set of gradings of $\phi$ and the set of
gradings of $\graph(\phi)$.

We will refer to this as the {\em canonical bijection} when using the canonical grading
$\ti{\Delta}\in\Lag^N(\R^{2n,-}\times\R^{2n})$ in Remark~\ref{Maslov diagonal}.  
In particular, the diagonal in $M^- \times M$ has a
canonical grading induced by the canonical bijection from the
canonical grading of the identity on $M$.
\item
Any grading $\sigma_L^N$ of a Lagrangian submanifold $L \subset M$
induces a grading of $L \subset M^-$ via the diffeomorphism $\Lag^N(M^-) \to \Lag^N(M)$. 
\item
Given graded Lagrangian submanifolds $L_0 \subset M_0, L_1 \subset M_1$, the 
product $L_0 \times L_1 \subset M_0 \times M_1$ inherits a grading from \eqref{fiberzn}.  
\item
Given a graded symplectomorphism $\phi: M_0 \to M_1$ and a graded Lagrangian
submanifold $L \subset M_0$, the image $\phi(L) \subset M_1$ inherits a grading
by composition $\sigma_{\phi(L)}^N = \phi^N \circ \sigma_L^N$.    
\end{enumerate}
\end{remark}

\begin{example}  \label{maslovexamples}
\begin{enumerate}
\item
Let $\Lag^2(M)$ be the bundle 
whose fiber over $m$ is the space of oriented Lagrangian subspaces of
$T_m M$. Then $\Lag^2(M) \to \Lag(M)$ is a $2$-fold Maslov covering.  
A $\Lag^2(M)$-grading of a Lagrangian $L \subset M$ is equivalent to an orientation on $L$.
\item
By \cite[Section 2]{se:gr}, any symplectic manifold $M$ with
$H^1(M)=0$ and minimal Chern number $N_M$ admits an $N$-fold Maslov
covering $\Lag^N(M)$ iff $N$ divides $2N_M$.  Any Lagrangian with
minimal Maslov number $N_L$ admits a $\Lag^N(M)$-grading iff $N$
divides $N_L$.  In particular, if $H^1(M)=0$ and $L$ is simply connected, 
then $N_L = 2N_M$ and $L$ admits a $\Lag^{2N_M}(M)$ grading.
\item
Suppose that $[\omega]$ is integral, $[\omega] = (1/l) c_1(TM) $, and
$\cL$ is a line bundle with connection $\nabla$ and curvature
$\curv(\nabla) = (2\pi /i) \omega$.
This induces a $2l$-fold Maslov cover $\Lag^{2l}(M) \to \Lag(M)$, see
\cite[Section 2b]{se:gr}.  Let $L \subset M$ be a Bohr-Sommerfeld
monotone Lagrangian as in Remark \ref{BS}.  
A grading of $L$ is equivalent to a choice of
(not necessarily horizontal) section of $\cL | L$ whose $l$-th tensor
power is $\phi_L^{\cK}$; that is, a choice of the section $\exp(2 \pi
i \psi) \phi_L^{\cL}$ in \eqref{stupid}.
\end{enumerate}
\end{example}

\begin{definition}
Let $\Lambda_0,\Lambda_1\subset V$ be a transverse pair of Lagrangian
subspaces in a symplectic vector space $V$ and let $\tilde\Lambda_0,
\tilde\Lambda_1\in\Lag^N(V)$ be gradings. 
The {\em degree} $d(\tilde\Lambda_0,\tilde\Lambda_1)\in\Z_N$ is
defined as follows.  Let
$\tilde\gamma_0,\tilde\gamma_1:[0,1] \to \Lag^N(V)$
be paths with common starting point
$\tilde\gamma_0(0) = \tilde\gamma_1(0)$ and end points
$\tilde\gamma_j(1) =\tilde\Lambda_j$.
Let $\gamma_j:[0,1] \to \Lag(V)$ denote their image under the
projection $\Lag^N(V) \to \Lag(V)$ and define
\begin{equation}\label{the degree}
 d(\tilde\Lambda_0,\tilde\Lambda_1) := \tfrac12\dim(\Lambda_0) +
I(\gamma_0,\gamma_1 ) \quad \text{mod} \ N ,
\end{equation}
where $I(\gamma_0,\gamma_1)$ denotes the Maslov index for the pair of
paths as in \cite{vi:in,rs:maslov}.
\end{definition}

Let us recall from \cite{rs:maslov} that the Maslov index for a pair of
paths with regular crossings (in particular with a finite set of crossings
$\cC:=\{s\in[0,1]\,|\, \gamma_0(s) \cap \gamma_1(s) \neq \{ 0 \} \}$)
is given by the sum of crossing numbers with the endpoints weighted by $1/2$,
$$ I(\gamma_0,\gamma_1) = \hh \sum_{s \in \cC\cap\{ 0,1 \}} {\rm sign}(\Gamma(\gamma_0,\gamma_1,s)) +
\sum_{s \in \cC\cap(0,1)} {\rm sign}(\Gamma(\gamma_0,\gamma_1,s)) . $$
Each crossing operator $\Gamma(\gamma_0,\gamma_1,s)$ is defined 
on $v\in\gamma_0(s)\cap\gamma_1(s)$ by
fixing Lagrangian complements $\gamma_0(s)^c$, $\gamma_{1}(s)^c$ of
$\gamma_0(s),\gamma_1(s)$ and setting
\begin{equation} \label{crossing}  
\Gamma(\gamma_0,\gamma_1,s)v = \tfrac d{dt}\bigr|_{t=0}
\omega(v,w(t) - w'(t)) \end{equation}
where $w(t) \in \gamma_0(s)^c$ such that $v+w(t)\in \gamma_0(s + t)$
and $w'(t) \in \gamma_1(s)^c$ such that $v+w'(s + t)\in \gamma_1(s)$.

\begin{remark} \label{alt deg}
The degree can alternatively be defined by fixing $\ti{\gamma}_0\equiv\ti{\Lambda}_0$
and choosing a path $\tilde\gamma:[0,1] \to \Lag^N(V)$ from $\tilde\gamma(0) = \tilde\Lambda_0$ to
$\tilde\gamma(1) =\tilde\Lambda_1$ such that the crossing form
$\Gamma(\gamma,\Lambda_0,0)$ of the underlying path $\gamma:[0,1] \to \Lag(V)$ 
is positive definite at $s=0$.  
Then the degree
$$
d(\tilde\Lambda_0,\tilde\Lambda_1) 
= \tfrac {\dim\Lambda_0} 2 + I(\Lambda_0,\gamma)
= - \sum_{s\in(0,1)} {\rm sign}(\Gamma(\gamma,\Lambda_0,s))
= - I'(\gamma,\Lambda_0) \qquad \text{mod}\;N
$$
is given by the Maslov index $I'$ of $\gamma|_{(0,1)}$ (not counting the endpoints)
relative to $\Lambda_0$.
Equivalently, we have
$$
d(\tilde\Lambda_0,\tilde\Lambda_1) 
= I'(\gamma^{-1},\Lambda_0) \qquad \text{mod}\;N
$$
for the reversed path 
$\gamma^{-1}:[0,1] \to \Lag(V)$ from
 $\gamma^{-1}(0) = \Lambda_1$ to
$\gamma^{-1}(1) = \Lambda_0$ such that the crossing form
$\Gamma(\gamma^{-1},\Lambda_0,1)$ is negative definite at $s=1$.  
\end{remark}

\begin{lemma} (Index theorem for once-punctured disks) 
\label{grFr}
Let $\Lambda_0,\Lambda_1\subset V$ be a transverse pair of Lagrangian
subspaces with gradings $\tilde\Lambda_0,\tilde\Lambda_1\in\Lag^N(V)$.
Then for any smooth path of graded Lagrangian subspaces
$\tilde\Lambda:[0,1]\to\Lag^N(V)$ with endpoints
$\tilde\Lambda(j)=\tilde\Lambda_j$, $j=0,1$ we have
$$
d(\tilde\Lambda_0,\tilde\Lambda_1) = {\rm Ind}(D_{V,\Lambda})   \quad \text{\em mod} \ N  .
$$ 
Here $D_{V,\Lambda}$ is any Cauchy-Riemann operator in $V$ on the
disk $D$ with one outgoing strip-like end $(0,\infty)\times
[0,1]\hookrightarrow D$ and with boundary conditions given by
$\Lambda$ (the projection of $\tilde\Lambda$ to $\Lag(V)$) such that
$\Lambda(j)=\Lambda_j$ is the boundary condition over the boundary
components $(0,\infty)\times\{j\}$, $j=0,1$ of the end.
\end{lemma}
\begin{proof}
It suffices to prove the index identity for a fixed path
$\tilde\Lambda$.  Indeed, if $\tilde\Lambda'$ is any other path with
the same endpoints then we have ${\rm Ind}(D_{V,\Lambda})-{\rm
Ind}(D_{V,\Lambda'}) = {\rm Ind}(D_{V,\Lambda})+{\rm
Ind}(D_{V,-\Lambda'}) = {\rm Ind}(D_{V,\Lambda\#(-\Lambda')})$ by
gluing.  Here the last Cauchy-Riemann operator is defined on the disk
with no punctures and with boundary conditions given by the loop
$\Lambda\#(-\Lambda')$.  Since the loop lifts to a loop
$\tilde\Lambda\#(-\tilde\Lambda')$ in $\Lag^N(V)$, its Maslov index
(and thus index) is $0$ modulo $N$.

By Remark \ref{alt deg}, the degree can be defined by a path
$\ti{\Lambda}$ from $\ti{\Lambda}_1$ to $\ti{\Lambda}_0$ whose
projection $\Lambda$ has negative definite crossing form at $s=1$. 
The sum of crossing numbers in
$d(\tilde\Lambda_0,\tilde\Lambda_1) = 
\sum_{s\in(0,1)} {\rm sign}(\Gamma(\Lambda,\Lambda_0,s))$ 
is the Maslov index $I_H(\Lambda)$ in \cite[Lemma 11.11]{se:bo} 
and hence equals to the Fredholm index
${\rm Ind}(D_{V,\Lambda})$ over the half space, or the conformally
equivalent disk with strip-like end.
This conformal isomorphism takes the boundary ends $(-\infty,-1)$
resp.\ $(1,\infty)$ in the half space $\{\Im z \ge 0\}$ 
(over which $\Lambda$ equals to $\Lambda_1$ resp.\ $\Lambda_0$)
to $\{1\}\times(1,\infty)$ resp.\ $\{0\}\times(1,\infty)$ in the strip-like end. 
\end{proof}

\begin{lemma} \label{degprop}
The degree map satisfies the following properties.
\begin{enumerate}
\item (Additivity)
If $V = V' \times V''$ then
$$ d( \ti{\Lambda}_0' \times^N \ti{\Lambda}_0'',
\ti{\Lambda}_1' \times^N \ti{\Lambda}_1'') = d( \ti{\Lambda}_0',\ti{\Lambda}_1') +
d( \ti{\Lambda}_0'',\ti{\Lambda}_1'') $$
for $\ti{\Lambda}_j',\ti{\Lambda}_j''$ graded Lagrangian subspaces in $V',V''$ respectively, $j = 0,1$.  
\item (Multiplicativity)
For $\ti{\Lambda}_0,\ti{\Lambda}_1$ graded Lagrangian subspaces and any $c \in \Z_N$
$$ d( \ti{\Lambda}_0 ,  c \cdot \ti{\Lambda}_1) = c + d(  \ti{\Lambda}_0 , \ti{\Lambda}_1) .
$$
\item (Skewsymmetry)
For $\ti{\Lambda}_0,\ti{\Lambda}_1$ graded Lagrangian subspaces
$$ d( \ti{\Lambda}_0 , \ti{\Lambda}_1) + d(  \ti{\Lambda}_1 , \ti{\Lambda}_0)
= \dim\Lambda_0 
= d( \ti{\Lambda}_0 , \ti{\Lambda}_1) + d(  \ti{\Lambda}_0^- , \ti{\Lambda}_1^-).
$$
\item (Diagonal)
For a transverse pair $\ti{\Lambda}_0,\ti{\Lambda}_1$ of graded 
Lagrangian subspaces in $V$ 
and $\ti\Delta$ the canonically graded diagonal in $V^-\times V$
$$ d( \ti{\Delta} , \ti{\Lambda}_0^- \times^N \ti{\Lambda}_1 ) 
= d(  \ti{\Lambda}_0 , \ti{\Lambda}_1) .
$$
\end{enumerate}
\end{lemma}

\begin{proof}
The first three properties are standard, see \cite[Section 2d]{se:gr}.
We prove the diagonal property to make sure all our sign conventions match up.
For that purpose we fix $\ti{L}\in\Lag^N(V)$ and choose the following paths $\ti\gamma_{..}$ 
of graded Lagrangian subspaces (with underlying paths $\gamma_{..}$ of Lagrangian subspaces):
\begin{itemize}
\item
$\ti{\gamma}_{0}:[-1,1]\to\Lag^N(V)$ from $\ti{\gamma}_{0}(-1)=\ti{L}$
to $\ti{\gamma}_{0}(1)=\ti{\Lambda}_{0}$ such that
${{\ti{\gamma}_{0}}|_{[-1,0]}\equiv\ti{L}}$,
\item
$\ti{\gamma}_{1}:[-1,1]\to\Lag^N(V)$ from $\ti{\gamma}_{1}(-1)=\ti{L}$ to 
$\ti{\gamma}_{1}(1)=\ti{\Lambda}_{1}$, 
such that ${\gamma}_{1}|_{[-1/2,0]}\equiv JL\pitchfork L$ and 
${\gamma}_{1}|_{[-1,-1/2]}$ is a smoothing of $t\mapsto e^{\pi(1+t)J}L$.
\item
$\ti{\gamma}:[-1,1]\to\Lag^N(V^-\times V)$ 
starting with ${\ti{\gamma}}|_{[-1,-1/2]}=(\ti{\gamma}_1^-\times^N\ti\gamma_0)|_{[-1,-1/2]}$,
ending at $\ti{\gamma}|_{[0,1]}\equiv\ti{\Delta}$, and such that
$\gamma|_{[-\frac 12,0]}$ is a smoothing of
$t\mapsto\{((2t+1)x+Jy,x+(2t+1)Jy) | x,y\in L \}$.
(The lift to graded subspaces matches up since $\gamma|_{[-1,0]}$
is exactly the path of \eqref{diag grading}
which defines $\ti{\Delta}$ by connecting it to  $\ti{L}^-\times\ti{L}$.)
\end{itemize}

Note that we have $I(\gamma_0,\gamma_1)|_{[-1,0]}=-\frac
12\dim\Lambda_0$ and $I(\gamma,\gamma_0^-\times\gamma_1)|_{[-1,0]} =
I(\gamma_1^-,\gamma_0^-)|_{[-1,0]} + I(\gamma_0,\gamma_1)|_{[-1,0]} =
-\dim\Lambda_0$ since $\gamma|_{[-\frac 12,0]}$ is transverse to
$L^-\times JL$.  With these preparations we can calculate
\begin{align*}
 d(\ti{\Lambda}_0,\ti{\Lambda}_{1}) 
&\,=\; \tfrac12\dim\Lambda_0 + I(\gamma_{0},\gamma_{1}) 
\;=\;  I(\gamma_0,\gamma_1)\bigr|_{[0,1]} \\
&\,=\;  I(\Delta,\gamma_0^-\times\gamma_1)\bigr|_{[0,1]} \\
&\,=\; \dim\Lambda_0 + I(\gamma,\gamma_0^-\times\gamma_1) 
\;=\; d(\ti{\Delta},\ti{\Lambda}_0^-\times^N\ti{\Lambda}_1).
\end{align*}
Here the identity of the Maslov indices over the interval $[0,1]$
follows from identifying the intersections
$K(s):=\gamma_0\cap\gamma_1 \cong \Delta\cap (\gamma_0^-\times\gamma_1)$ 
and the crossing forms $\Gamma(s),\hat\Gamma(s) : K(s) \to \R$
at regular crossings $s\in[0,1]$ (after a homotopy of the paths to regular crossings).
Fix Lagrangian complements $\gamma_0(s)^c$ and $\gamma_1(s)^c$, 
then for $v\in K(s)$ pick $w_i(t)\in \gamma_i(s)^c$ such that $v+w_i(t)\in \gamma_i(s+t)$.
For the corresponding vector $\hat{v}=(v,v)\in \Delta\cap(\gamma_0^-\times\gamma_1)$
we can pick $\hat{w}(t)=(0,0)\in \Delta^c$ satisfying $\hat{v}+\hat{w}(t)\in\Delta$
and $\hat{w}'(t)=(w_0,w_1)\in \gamma_0(s)^c\times\gamma_1(s)^c$ 
satisfying $\hat{v}+\hat{w}'(t)\in (\gamma_0\times\gamma_1)(s+t)$
to identify the crossing forms
\begin{align*}
\hat{\Gamma}(s)\hat{v} 
&= \tfrac d{dt}\bigr|_{t=0} (-\omega\oplus\omega)(\hat{v},\hat{w}(t) - \hat{w}'(t)) \\
&= \tfrac d{dt}\bigr|_{t=0} \bigl( - \omega(v,-w_0(t)) + \omega(v, - w_1(t) ) \bigr) \\
&= \tfrac d{dt}\bigr|_{t=0}  \omega(v, w_0(t)-w_1(t)) = \Gamma(s)v .
\end{align*}

\vspace{-5mm}
\end{proof}

If $L_0,L_1\subset M$ are $\Lag^N(M)$-graded Lagrangians and intersect
transversally then one obtains a {\em degree map}
$$  \cI(L_0,L_1):=L_0\cap L_1 \to \Z_N ,
\qquad x\mapsto |x|:=d(\sigma_{L_0}^N(x),\sigma_{L_1}^N(x)) . $$ 
More generally, if $L_0,L_1$ do not necessarily intersect
transversally, then we can pick a Hamiltonian perturbation
$H:[0,1]\times M\to \R$ such that its time $1$ flow $\phi_1: M \to M$
achieves transversality $\phi_1(L_0)\pitchfork L_1$.  Then the
Hamiltonian isotopy and the grading on $L_0$ induce a grading on
$\phi_1(L_0)$, which is transverse to $L_1$.  The degree map is then
defined on the perturbed intersection points, $d:
\cI(L_0,L_1):=\phi_1(L_0)\cap L_1 \to \Z_N$.

\subsection{Graded generalized Lagrangian correspondences} 
\label{graded cor}

In this section we extend the grading and degree constructions to generalized Lagrangian
correspondences and discuss their behaviour under geometric composition and insertion of
the diagonal.

\begin{definition} \label{gen grad s}
Let $M$ and $M'$ be symplectic manifolds equipped with $N$-fold Maslov coverings.
Let $\ul{L}=(L_{01},\ldots, L_{(r-1)r})$ be a generalized Lagrangian correspondence
from $M$ to $M'$ (i.e.\ $L_{(j-1)j} \subset M_{j-1}^-\times M_{j}$ for a sequence 
$M=M_1, \ldots, M_{r}=M'$ of symplectic manifolds).
A {\em grading} on $\ul{L}$ consists of a collection of $N$-fold
Maslov covers $\Lag^N(M_j)\to M_j$ and gradings of the Lagrangian correspondences 
$L_{(j-1)j}$ with respect to $\Lag^N(M_{j-1}^-\times M_j)$, where
the Maslov covers on $M_1=M$ and $M_{r}=M'$ are the fixed ones.
\end{definition}

A pair of graded generalized Lagrangian correspondences $\ul{L}_1$ and $\ul{L}_2$ 
from $M$ to $M'$ (with fixed Maslov coverings)
defines a cyclic Lagrangian correspondence $\ul{L}_1\#(\ul{L}_2)^t$,
which is graded in the following sense.

\begin{definition} \label{gen grad c}
Let $\ul{L}=(L_{01},\ldots, L_{r(r+1)})$ be a cyclic generalized Lagrangian correspondence
(i.e.\ $L_{j(j+1)} \subset M_{j}^-\times M_{j+1}$ for a cyclic sequence 
$M_0, M_1, \ldots, M_{r+1}=M_0$ of symplectic manifolds).
An {\em $N$-grading} on $\ul{L}$ consists of a collection of $N$-fold
Maslov covers $\Lag^N(M_j)\to M_j$ and gradings of the Lagrangian correspondences 
$L_{j(j+1)}$ with respect to $\Lag^N(M_{j}^-\times M_{j+1})$.  
\end{definition}

In the following, we will consider a cyclic generalized Lagrangian correspondence $\ul{L}$
and assume that it intersects the generalized diagonal transversally, i.e.\
\begin{equation}\label{diag trans}
\bigl(L_{01}\times L_{12}\times\ldots\times L_{r(r+1)}\bigr)\pitchfork
\bigl(\Delta_{M_0}^-\times\Delta_{M_1}^-\times\ldots\times\Delta_{M_r}^-\bigr)^T,
\end{equation}
where $\Delta^-_{M}\subset M\times M^-$ denotes the (dual of the) diagonal and 
$M_0\times M_0^-\times M_1\times \ldots \times M_{r}^-\to 
M_0^- \times M_1\times \ldots\times M_{r}^-\times M_0$, $Z\mapsto Z^T$ 
is the transposition of the first to the last factor.
In section \ref{sequences} this transversality will be achieved by a suitable Hamiltonian isotopy. 
It ensures that the above transverse intersection cuts out a finite set, which we identify
with the {\em generalized intersection points}
\begin{align*}
\cI(\ul{L})&:=
\times_{\Delta_{M_0}}\bigl( L_{01} \times_{\Delta_{M_1}} L_{12} \ldots \times_{\Delta_{M_r}} L_{r(r+1)}\bigr)\\
&=
\bigl\{ \ul{x}=(x_0,\ldots,x_r) \in M_0 \times \ldots \times M_r \,\big|\, 
(x_0,x_1) \in L_{01}, \ldots ,(x_r,x_0) \in L_{r(r+1)} \bigr\} .
\end{align*}

\begin{remark} \label{rmk:canonical ident}
Consider two cyclic generalized Lagrangian correspondences 
\begin{align*}
\ul{L}&=(L_{01},\ldots,L_{(j-1)j},L_{j(j+1)}, \ldots, L_{r(r+1)}), \\
\ul{L}'&=(L_{01},\ldots,L_{(j-1)j}\circ L_{j(j+1)}, \ldots, L_{r(r+1)})
\end{align*}
such that the composition $L_{(j-1)j}\circ L_{j(j+1)}$ is embedded
in the sense of Definition~\ref{embedded}.
Then the generalized intersection points
\begin{align*}
\cI(\ul{L})&=
\bigl\{(\ldots,x_{j-1},x_j,x_{j+1},\ldots) \in \ldots \times M_{j-1}\times M_j \times M_{j+1} \ldots  
\,\big|\,  \\
&\qquad\qquad\qquad\qquad
\ldots,(x_{j-1},x_j)\in L_{(j-1)j}, (x_{j},x_{j+1})\in L_{j(j+1)}, \ldots  \bigr\}\\
&=
\bigl\{ (\ldots,x_{j-1},x_{j+1},\ldots) \in \ldots \times M_{j-1} \times M_{j+1} \ldots  \,\big|\, \\
&\qquad\qquad\qquad\qquad
 \ldots , (x_{j-1},x_{j+1})\in L_{(j-1)j}\circ L_{j(j+1)}, \ldots  \bigr\}
=\cI(\ul{L}')
\end{align*}
are canonically identified, since the intermediate point 
$x_{j}\in M_j$ with $(x_{j-1},x_j)\in L_{(j-1)j}$ and $(x_{j},x_{j+1})\in L_{j(j+1)}$ 
is uniquely determined by the pair $(x_{j-1},x_{j+1})\in L_{(j-1)j}\circ L_{j(j+1)}$.
\end{remark}

Now an $N$-grading on $\ul{L}$ induces an $N$-fold Maslov covering on
${M}:=M_0^-\times M_1\times\ldots\times M_r\times M_r^-\times M_0$
and a grading of ${L}:=L_{01}\times L_{12}\times\ldots\times L_{r(r+1)}$.
In addition, we have a grading on
$\Delta^T:=(\Delta_{M_0}^-\times\Delta_{M_1}^-\times\ldots\times\Delta_{M_r}^-)^T$
from the canonical grading on each factor.
In order to define a degree we then identify generalized intersection points 
$\ul{x}=(x_0,x_1,\ldots,x_r)$ with the actual intersection points 
${x}=(x_0,x_1,x_1,\ldots,x_r,x_r,x_0)\in{L}\cap\Delta^T$.

\begin{definition}\label{def deg}
Let $\ul{L}$ be a graded cyclic generalized Lagrangian correspondence $\ul{L}$ 
that is transverse to the diagonal \eqref{diag trans}. Then the {\em degree} is 
$$
\cI(\ul{L}) \to \Z_N, \qquad 
\ul{x}\mapsto|\ul{x}|= d(\sigma^N_{{L}}({x}),\sigma^N_{\Delta^T}({x})) .
$$
\end{definition}

\begin{lemma} \label{alt seq deg}
Alternatively, the degree is defined as follows:
\begin{enumerate}
\item
Pick any tuple of Lagrangian subspaces 
$\Lambda_i'\in\Lag(T_{x_i}M_i)$, $\Lambda_i''\in\Lag(T_{x_i}M_i^-)$,
$i=0,\ldots,r$ whose product
is transverse to the diagonal, $\Lambda_i'\times\Lambda_i''\pitchfork\Delta_{T_{x_i}M_i}$.
Then there exists a path (unique up to homotopy) $\gamma:[0,1]\to \Lag(T_{{x}}{M})$
from $\gamma(0)=T_{{x}}{L}$ to
$\gamma(1)=\Lambda_0''\times\Lambda_1'\times\ldots\times \Lambda_r'\times\Lambda_r''\times\Lambda_0'$
that is transverse to the diagonal at all times,
$\gamma(t) \pitchfork T_{x}\Delta^T$.
We lift the grading $\sigma^N_{{L}}({x})\in\Lag^N(T_{{x}}{M})$
along this path and pick preimages under the graded product map \eqref{fiberzn} to define
$\tilde{\Lambda}_{i}'\in\Lag^N(T_{x_i}M_i)$ and
$\tilde{\Lambda}_{i}''\in\Lag^N(T_{x_i}M_i^-)$. Then
$$
|\ul{x}|= \sum_{i=0}^r d(\tilde{\Lambda}_{i}',\tilde{\Lambda}_{i}''^-) .
$$
\item
If $\ul{L}$ has even length $r+1\in 2\N$ then it defines an $N$-fold Maslov cover on
$\widetilde{M}:=M_0^-\times M_1 \times M_2^-\times\ldots \times M_{r}$ 
and a pair of graded Lagrangian submanifolds,
\begin{align*}
{L}_{(0)}&:= L_{01}\times L_{23} \times \ldots \times L_{(r-1)r}  \subset \widetilde{M} , \\
{L}_{(1)}&:= ( L_{12} \times L_{34} \times \ldots \times L_{r(r+1)})^T  \subset \widetilde{M}^- ,
\end{align*}
where we denote by $M_1^-\times \ldots \times M_{r}^-\times M_0 \to 
M_0\times M_1^-\times \ldots\times M_{r}^-$, $Z\mapsto Z^T$ 
the transposition of the last to the first factor.
If $\ul{L}$ has odd length $r+1\in 2\N +1$ we insert the diagonal  
$\Delta_{M_0}\subset M_{0}^-\times M_0 = M_{r+1}^-\times M_0$ (with its canonical grading) 
before defining a pair of graded Lagrangian submanifolds as above.
By \eqref{diag trans} the Lagrangians intersect transversally ${L}_{(0)}\pitchfork{L}_{(1)}^-$,
and this intersection is canonically identified with $\cI(\ul{L})$.
Then for $\ul{x}\in\cI(\ul{L})$ corresponding to $y\in{L}_{(0)}\cap{L}_{(1)}^-$ 
we have
$$  
|\ul{x}|=|{y}|=d(\sigma_{{L}_{(0)}}^N({y}),\sigma_{{L}_{(1)}}^N({y})^-) . 
$$ 
\end{enumerate}
\end{lemma}

\begin{proof}
In (a) we use the fact that the path $\gamma$ has zero Maslov index to rewrite
$$
d(\sigma^N_{{L}}({x}),\sigma^N_{\Delta^T}({x}))
=d(\tilde\Lambda_0'\times^N\tilde\Lambda_0''\times^N\ldots\times^N \tilde\Lambda_r'\times^N\tilde\Lambda_r'',
\tilde\Delta_{T_{x_0}M_0}^-\times^N \ldots \times^N\tilde\Delta_{T_{x_r}M_r}^-),
$$
where we moreover transposed the factors.
Now by Lemma \ref{degprop} the right-hand side can be written as the sum over
$d(\tilde\Lambda_i'\times^N\tilde\Lambda_i'',\tilde\Delta_{T_{x_i}M_i}^-)
=d(\tilde\Lambda_i',\tilde\Lambda_i''^-)$.

In (b) note that a reordering of the factors identifies 
the pair of graded Lagrangians $(L_{(0)}\times L_{(1)},\Delta_{\widetilde M}^-)$
with $(L,\Delta^T)$ for $r$ odd. 
So Lemma \ref{degprop} implies
$$
d(\sigma^N_{{L}}({x}),\sigma^N_{\Delta^T}({x}))
= d(\sigma^N_{L_{(0)}}({y})\times^N \sigma^N_{L_{(1)}}({y}),\tilde\Delta_{T_{(y,y)}\widetilde M}^-)
= d(\sigma^N_{L_{(0)}}({y}) , \sigma^N_{L_{(1)}}({y})^-) .
$$
For $r$ even the same argument proves
$$
d(\sigma^N_{L_{(0)}}({y}) , \sigma^N_{L_{(1)}}({y})^-) 
= d\bigl( \sigma^N_{L}({x})\times^N \ti\Delta_{T_{x_0}M_0} ,
(\ti\Delta_{T_{x_0}M_0}^-\times\ldots\times\ti\Delta_{T_{x_r}M_r}^-\times\ti\Delta_{T_{x_0}M_0}^-)^T ) ,
$$
which equals to $d(\sigma^N_{{L}}({x}),\sigma^N_{\Delta^T}({x}))$ by Lemma \ref{insertdiag} (b) below.
\end{proof}

The following Lemma describes the effect of inserting a diagonal on the 
grading of generalized Lagrangian correspondences.
Part (a) addresses noncyclic correspondences, whereas (b) applies to cyclic correspondences
with $\Lambda=T_{(x_0,x_1,\ldots,x_r,x_0)}( L_{01}\times L_{12}\times \ldots \times L_{r(r+1)})$,
$K=T_{(x_0,x_0,x_1,\ldots,x_r)}(\Delta_{M_0}^-\times\Delta_{M_1}^-\times\ldots\times\Delta_{M_r}^-)$,
$V_0=T_{x_0}M_0$, and $V_1=T_{(x_1,\ldots,x_r)}(M_1\times M_1^-\times\ldots\times M_r\times M_r^-)$.

\begin{lemma} \label{insertdiag}
Let $V_0, V_1, V_2$ be symplectic vector spaces. 
\begin{enumerate}
\item
Let
$\ti{\Lambda}_{0}\subset\Lag^N(V_0)$,
$\ti{\Lambda}_{01}\subset\Lag^N(V_0^-\times V_1)$,
$\ti{\Lambda}_{12}\subset\Lag^N(V_1^-\times V_2)$,
and $\ti{\Lambda}_{2}\subset\Lag^N(V_2^-)$ be graded Lagrangian subspaces. 
If the underlying Lagrangian subspaces are transverse then
$$
d(\ti{\Lambda}_0\times^N\ti{\Lambda}_{12},\ti{\Lambda}_{01}^-\times^N\ti{\Lambda}_2^-)
= d(\ti{\Lambda}_0\times^N\ti{\Delta}_1\times^N\ti{\Lambda}_2,\ti{\Lambda}_{01}^-\times^N\ti{\Lambda}_{12}^-) .
$$
\item
Let
$\ti{\Lambda}\subset\Lag^N(V_0^-\times V_1\times V_0)$ and
$\ti{K}\subset\Lag^N(V_0\times V_0^-\times V_1)$ be graded Lagrangian subspaces. 
If the underlying Lagrangian subspaces are transverse then
$$
d(\ti{\Lambda}\times^N\ti{\Delta}_0,(\ti{K}\times^N\ti{\Delta}_0^-)^T)
= d(\ti{\Lambda},\ti{K}^T) ,
$$
with the transposition $V_0\times W \to W \times V_0$, $Z\mapsto Z^T$.
\end{enumerate}
\end{lemma}
\begin{proof}
To prove (a) pick a path $\gamma_{0112}:[0,1]\to\Lag(V_0\times V_1^-\times V_1 \times V_2^-)$
from $\gamma_{0112}(0)={\Lambda}_{01}^-\times{\Lambda}_{12}^-$
to a split Lagrangian subspace 
$\gamma_{0112}(1)=\Lambda_0'\times\Lambda_1'\times\Lambda_1''\times\Lambda_2'$
that is transverse to ${\Lambda}_0\times{\Delta}_1\times{\Lambda}_2$ at all times
and hence has Maslov index $I(\gamma_{0112},{\Lambda}_0\times{\Delta}_1\times{\Lambda}_2)=0$.
We can homotope this path with fixed endpoints to 
$\gamma_{0112}=\gamma_{01}\times\gamma_{12}:[0,1]\to\Lag(V_0\times V_1^-)\times \Lag(V_1 \times V_2^-)$
that may intersect $\Lambda_0\times\Delta_1\times\Lambda_2$ but still has vanishing Maslov index.
We lift the grading along the paths $\gamma_{01}$ and $\gamma_{12}$ and pick preimages 
under the graded product map \eqref{fiberzn} to obtain gradings 
$\tilde{\Lambda}_{0}'\in\Lag^N(V_0)$,
$\tilde{\Lambda}_{1}'\in\Lag^N(V_1^-)$,
$\tilde{\Lambda}_{1}''\in\Lag^N(V_1)$,
$\tilde{\Lambda}_{2}'\in\Lag^N(V_2^-)$.
With these we calculate, using Lemma \ref{degprop}
\begin{align*}
d(\ti{\Lambda}_0\times^N\ti{\Lambda}_{12},\ti{\Lambda}_{01}^-\times^N\ti{\Lambda}_2^-)
&= d(\ti{\Lambda}_0\times^N\ti{\Lambda}_{1}''^-\times^N\ti{\Lambda}_{2}'^-,
\ti{\Lambda}_{0}'\times^N\ti{\Lambda}_{1}'\times^N\ti{\Lambda}_2^-) \\
&= d(\ti{\Lambda}_0,\ti{\Lambda}_{0}') 
+ d(\ti{\Lambda}_{1}''^-,\ti{\Lambda}_{1}') 
+ d(\ti{\Lambda}_{2}'^-,\ti{\Lambda}_2^-) \\
&= d(\ti{\Lambda}_0,\ti{\Lambda}_{0}') 
+ d(\ti{\Delta}_1,\ti{\Lambda}_{1}'\times^N\ti{\Lambda}_{1}'') 
+ d(\ti{\Lambda}_{2},\ti{\Lambda}_2') \\
&= d(\ti{\Lambda}_0\times^N\ti{\Delta}_1\times^N\ti{\Lambda}_2,
\ti{\Lambda}_{0}'\times^N\ti{\Lambda}_{1}'\times^N\ti{\Lambda}_{1}''\times^N\ti{\Lambda}_{2}') \\
&= d(\ti{\Lambda}_0\times^N\ti{\Delta}_1\times^N\ti{\Lambda}_2,\ti{\Lambda}_{01}^-\times^N\ti{\Lambda}_{12}^-).
\end{align*}
The first and last degree identity are due to the vanishing of the Maslov index
$$
0 \;=\; I(\Lambda_0\times\Delta_1\times\Lambda_2,\gamma_{01}\times\gamma_{12})
\;=\; I(\Lambda_0\times\gamma_{12}^-,\gamma_{01}\times\Lambda_2^-)=0 .
$$
The identity of these Maslov indices follows from identifying the intersections
$K(s):=(\Lambda_0\times\gamma_{12}^-(s)) \cap (\gamma_{01}(s)\times\Lambda_2^-)
\cong (\Lambda_0\times\Delta_1\times\Lambda_{2}) \cap
(\gamma_{01}\times\gamma_{12})$ and the crossing form $\Gamma(s),\hat{\Gamma}(s) : K(s) \to \R$
given by \eqref{crossing} at regular crossings $s\in[0,1]$.
Fix Lagrangian complements
$\gamma_{01}(s)^c\subset V_0\times V_1^-$ and $\gamma_{12}(s)^c\subset V_1\times V_2^-$, 
then for $v=(v_0,v_1,v_2)\in K(s)$ we can pick 
$(w_1,w_2)(t)\in \gamma_{12}(s)^c$ 
such that $v+(0,w_1,w_2)(t)\in \Lambda_0\times\gamma_{12}(s+t)$
and $(w_0',w_1')(t)\in \gamma_{01}(s)^c$ 
such that $v+(w_0',w_1',0)(t)\in \gamma_{01}(s+t)\times\Lambda_2$.
For the corresponding vector $\hat{v}=(v_0,v_1,v_1,v_2)\in
(\Lambda_0\times\Delta_1\times\Lambda_{2}^-) \cap (\gamma_{01}^-\times\gamma_{12}^-)$
we have $\hat{v}+(0,0,0,0)\in(\Lambda_0\times\Delta_1\times\Lambda_2)$
and $\hat{v}+(w_0',w_1',w_1,w_2)(t)\in (\gamma_{01}\times\gamma_{12})(s+t)$.
With this we identify the crossing forms 
\begin{align*}
\hat{\Gamma}(s)\hat{v} 
&= \tfrac d{dt}\bigr|_{t=0} (\omega_0\oplus -\omega_1\oplus \omega_1 \oplus -\omega_2)
\bigl(\hat{v},(0,0,0,0) - (w_0',w_1',w_1,w_2)(t)\bigr) \\
&= \tfrac d{dt}\bigr|_{t=0} \bigl( - \omega_0(v_0, w_0') - \omega_1(v_1, w_1 - w_1') 
+ \omega_2(v_2, w_2) \bigr) \\
&= \tfrac d{dt}\bigr|_{t=0} (\omega_0\oplus -\omega_1\oplus \omega_2)
\bigl(v, (0,w_1,w_2)(t) - (w_0',w_1',0)(t) \bigr)
= \Gamma(s)v .
\end{align*}
This proves (a). To prove (b) we pick a path 
$\gamma:[0,1]\to\Lag(V_0^-\times V_1\times V_0)$
from $\gamma(0)={\Lambda}$ to a split Lagrangian subspace 
$\gamma(1)=\Lambda_0^-\times\Lambda_1\times\Lambda_0'\in\Lag(V_0^-)\times\Lag(V_1)\times\Lag(V_0)$
that is transverse to $K^T$ at all times and hence has Maslov index 
$$
0 = I(\gamma,K^T) = I(\gamma\times\Delta_0,(K\times\Delta_0^-)^T) .
$$
Here the equality of Maslov follows directly from the identification of the
trivial intersections
$(\gamma\times\Delta_0)\cap(K\times\Delta_0^-)^T\cong \gamma\cap K^T = \{0\}$.
Now we can lift the grading along $\gamma$ to obtain gradings 
$\tilde{\Lambda}_{0}\in\Lag^N(V_0)$,
$\tilde{\Lambda}_{1}\in\Lag^N(V_1)$,
$\tilde{\Lambda}_{0}'\in\Lag^N(V_0)$.
With these we calculate, using part (a) and the fact
that gradings are invariant under simultaneous transposition of both factors
\begin{align*}
d(\ti{\Lambda}\times^N\ti{\Delta}_0,(\ti{K}\times^N\ti{\Delta}_0^-)^T)
&= d(\ti\Lambda_0^-\times^N\ti\Lambda_1\times^N\ti\Lambda_0'\times^N\ti{\Delta}_0,
(\ti{K}\times^N\ti{\Delta}_0^-)^T) \\
&= d(\ti\Lambda_0'\times^N\ti{\Delta}_0\times^N\ti\Lambda_0^-\times^N\ti\Lambda_1,
\ti{\Delta}_0^-\times^N\ti{K}) \\
&= d(\ti\Lambda_0'\times^N \ti{K}^-,
\ti{\Delta}_0^-\times^N(\ti\Lambda_0^-\times^N\ti\Lambda_1)^-) \\
&= d(\ti{\Delta}_0\times^N(\ti\Lambda_0^-\times^N\ti\Lambda_1),
\ti\Lambda_0'^-\times^N \ti{K}) \\
&= d(\ti{K}^-,\ti\Lambda_0'^-\times^N(\ti\Lambda_0^-\times^N\ti\Lambda_1)^-) \\
&= d(\ti\Lambda_0'\times^N\ti\Lambda_0^-\times^N\ti\Lambda_1, \ti{K}) \\
&= d(\ti\Lambda_0^-\times^N\ti\Lambda_1\times^N\ti\Lambda_0',\ti{K}^T) 
\;=\;  d(\ti{\Lambda},\ti{K}^T) 
\end{align*}
\vspace{-10mm}

\end{proof}

In the rest of this section we investigate the effect of geometric
composition on the grading of Lagrangian correspondences.  This
requires a generalization of Viterbo's index calculations \cite{vi:in}.

First, we lift the composition map to Maslov covers.  Let $M_0, M_1,
M_2$ be symplectic manifolds equipped with $N$-fold Maslov coverings
$\Lag^N(M_j), j = 0,1,2$. We equip the products $M_i^-\times M_j$ and
$M_0^-\times M_1 \times M_1^-\times M_2$ with the induced Maslov
coverings $\Lag^N(M_i^-\times M_j)$ resp.\ $\Lag^N(M_{0}^-\times
M_{1}\times M_{1}^-\times M_{2})$. We denote by
$$\TD(M_1)\subset \Lag(M_0^-\times M_1 \times M_1^-\times M_2)
\bigr|_{M_0\times \Delta_{M_1} \times M_2} $$
the subbundle whose fibre over $(m_0,m_1,m_1,m_2)$ consists of the Lagrangian subspaces
$\Lambda_{0112}\subset T_{(m_0,m_1,m_1,m_2)}(M_0^- \times M_1 \times M_1^- \times M_2)$ 
that are transverse to the diagonal 
$\Delta_{0112}:=T_{m_0}M_0\times \Delta_{T_{m_1}M_1} \times T_{m_2}M_2$.
The linear composition of Lagrangian subspaces extends a smooth map
$$
\circ: \TD(M_1) \to \Lag(M_{0}^-\times M_{2}) , \qquad  
\Lambda_{0112}\mapsto \pi_{M_0\times M_2}\bigl( \Lambda_{0112}\cap\Delta_{0112}\bigr) .
$$
The preimage of $\TD(M_1)$ in the Maslov cover will be denoted by
$$
\TDN(M_1) \subset \Lag^N(M_{0}^-\times M_{1}\times M_{1}^-\times M_{2})
\bigr|_{M_0\times \Delta_{M_1} \times M_2} .
$$
Finally, recall that we have a canonical grading of the diagonal
$\ti{\Delta}_{M_1}\in\Lag^N(M_1^-\times M_1)$ and its dual
$\ti{\Delta}^-_{M_1}\in\Lag^N(M_1\times M_1^-)$,
and let us denote another exchange of factors by
$\Lag^N(M_{0}^-\times M_{2}\times M_{1}\times M_{1}^-)\to
\Lag^N(M_{0}^-\times M_{1}\times M_{1}^-\times M_{2})$, $\ti{\Lambda}\mapsto\ti{\Lambda}^T$.

\begin{lemma} The linear composition $\circ: \TD(M_1) \to \Lag(M_{0}^-\times M_{2})$ 
lifts to a unique smooth map $\circ^N: \TDN(M_1) \to
\Lag^N(M_{0}^-\times M_{2})$ with the property that
\begin{equation} \label{property}
\circ^N \bigl(\bigl( \ti{\Lambda}_{02} \times^N \ti{\Lambda}_{11}\bigr)^T \bigr) =
d(\ti{\Lambda}_{11},\ti{\Delta}_{M_1}^-) \cdot \ti{\Lambda}_{02}.\end{equation}
for all graded Lagrangians $\ti{\Lambda}_{02} \in \Lag^N(M_0^-\times M_2)$
and $\ti{\Lambda}_{11}\in\Lag^N(M_1\times M_1^-)$,
such that the underlying Lagrangian $\Lambda_{11}\in\Lag(M_1\times M_1^-)$ 
is transverse to the diagonal.
\end{lemma}

\begin{proof}
We denote by $\Lag(\R^{2n})$ the Lagrangian Grassmannian in $\R^{2n}$,
write $\dim M_i=2n_i$, and abbreviate $\R_{0112}:=\R^{2n_0,-}\times\R^{2n_1}\times \R^{2n_1,-}\times\R^{2n_2}$.
Let $\TD \subset \Lag(\R_{0112})$ be the subset of Lagrangian subspaces meeting the diagonal
$\R^{2n_0}\times \Delta_{\R^{2n_1}} \times \R^{2n_2}$ transversally.  
The linear composition map
$$ \Lag(\R_{0112}) \supset \TD \to \Lag(\R^{2n_0,-}\times \R^{2n_2}), \qquad 
\Lambda\mapsto \pi_{\R^{2n_0}\times\R^{2n_2}}\bigl(\Lambda\cap(\R^{2n_0}\times \Delta_{\R^{2n_1}}\times \R^{2n_2})\bigr)
$$
is $\Sp(2n_0)\times\Sp(2n_1)\times\Sp(2n_2)$-equivariant, 
and lifts to a unique
$\Sp^N(2n_0)\times\Sp^N(2n_1)\times\Sp^N(2n_2)$-equivariant map
\begin{equation} \label{equivariant} 
\Lag^N(\R_{0112}) \supset \TDN \to \Lag^N(\R^{2n_0,-}\times \R^{2n_2}) \end{equation}
with the property \eqref{property}.
On the other hand, the
restriction of $ \Fr(M_0) \times \Fr(M_1) \times \Fr(M_1) \times \Fr(M_2) $ 
to $M_0 \times \Delta_{M_1} \times M_2$ 
admits a reduction of the structure group to
$\Sp(2n_0) \times \Sp(2n_1) \times \Sp(2n_2)$, and similarly the restriction 
$$
\Fr^N_{0112}:= \bigl( \Fr^N(M_0) \times \Fr^N(M_1) \times \Fr^N(M_1) \times \Fr^N(M_2)  \bigr)
\big|_{M_0 \times \Delta_{M_1} \times M_2}
$$ 
admits a reduction of the structure group to 
$\Sp^N(2n_0) \times \Sp^N(2n_1) \times \Sp^N(2n_2)$.  
This group acts on $\Lag^N(\R_{0112})$ by the diagonal action of 
$\Sp^N(2n_1)$ on $\R^{2n_1}\times\R^{2n_1,-}$.
Finally, we use the associated fiber bundle construction to identify
\begin{align*}
&\Lag^N(M_{0}^-\times M_{1}\times M_{1}^-\times M_{2}) 
 \bigr|_{M_0\times \Delta_{M_1} \times M_2}  \\
&\cong \Fr^N_{0112} \times_{\Sp^N(2n_0) \times \Sp^N(2n_1) \times \Sp^N(2n_1) \times \Sp^N(2n_2)}
\Lag^N(\R_{0112}) \\
&\cong
\bigl( \Fr^N(M_0) \times \Fr^N(M_1) \times \Fr^N(M_2)  \bigr) 
\times_{\Sp^N(2n_0) \times \Sp^N(2n_1) \times \Sp^N(2n_2)} \Lag^N(\R_{0112})
\end{align*}
and
$$
\Lag^N(M_0^- \times M_2)=
\bigl( \Fr^N(M_0^-)\times\Fr^N(M_2) \bigr) \times_{\Sp^N(2n_0)\times\Sp^N(2n_2)}
\Lag^N(\R^{2n_0,-}\times\R^{2n_2}) .
$$
Then the forgetful map on the first factor and the
equivariant map \eqref{equivariant} on the second factor 
define the unique lift $\circ^N$.
\end{proof}

Now consider two graded Lagrangian correspondences $L_{01}\subset
M_0^-\times M_1$ and $L_{12}\subset M_1^-\times M_2$ and suppose that
the composition $L_{01} \circ L_{12}=:L_{02}\subset M_0^-\times M_2$
is smooth and embedded. The canonical section
$\sigma_{L_{02}}:L_{02}\to\Lag(M_0^-\times M_2)$ is given by the
linear composition $\circ$ applied to $(\sigma_{L_{01}}\times
\sigma_{L_{12}})|_{L_{01}\times_{\Delta_{M_1}}L_{12}}$.  The gradings
$\sigma^N_{L_{01}}$, $\sigma^N_{L_{12}}$ induce a grading on $L_{02}$,
\begin{equation} \label{dfn comp grad}
\sigma^N_{L_{02}}:= \circ^N \bigl(
\sigma^N_{L_{01}}\times^N \sigma^N_{L_{12}}\bigr)\big|_{L_{01}\times_{\Delta_{M_1}}L_{12}} ,
\end{equation}
where the map $\times^N$ is defined in \eqref{fiberzn}
and we identify $L_{02}\cong L_{01}\times_{\Delta_{M_1}}L_{12}$.

\begin{proposition} \label{gradingcomp} 
Let $L_0\subset M_0$, $L_{01}\subset M_0^-\times M_1$, $L_{12}\subset
M_1^-\times M_2$, and $L_2\subset M_2^-$ be graded Lagrangians such
that the composition $L_{01} \circ L_{12}=:L_{02}$ is embedded.  Then,
with respect to the induced grading on $L_{02}$, the degree map
$\cI(L_0 \times L_2, L_{02}) \to \Z_N$ is the pull-back of the degree
map $\cI( L_0\times L_{12}, L_{01}\times L_2) \to \Z_N$ under the
canonical identification \footnote{ Here it suffices to allow for
Hamiltonian perturbation on $M_0$ and $M_2$, i.e.\ replacing $L_0,L_2$
with $L_0':=\phi_1^{H_0}(L_0)$, $L_2':=(\phi_{1}^{H_2})^{-1}(L_2)$.
Then for every $(m_0,m_2)\in (L_0' \times L_2')\cap L_{02}$ there is a
unique $m_1\in M_1$ such that $(m_0,m_1)\in L_{01}$, $(m_1,m_2)\in
L_{12}$, and hence $(m_0,m_1,m_2)\in (L_0'\times L_{12})\cap
(L_{01}\times L_2')$.  
% CW: here
%The same identification will be used in \eqref{ident int}.  
} of intersection points.

\end{proposition}

\begin{proof}
Suppose for simplicity that Hamiltonian perturbations have been
applied to the Lagrangians $L_0,L_2$ such that $\cI(L_0 \times L_2,
L_{02})$ (and hence also $\cI(L_0\times L_{12},L_{01}\times L_2)$) is
the intersection of transverse Lagrangians.  We need to consider
$(m_0,m_1,m_2) \in (L_0\times L_{12})\cap (L_{01}\times L_2)$, which
corresponds to $(m_0,m_2)\in (L_0 \times L_2) \cap L_{02}$.  We
abbreviate the tangent spaces of the Lagrangians by $\Lambda_j =
T_{m_j} L_j$, $\Lambda_{ij} = T_{(m_i,m_j)} L_{ij}$, and
$\Delta_1=\Delta_{T_{m_1}M_1}$ and their graded lifts by
$\ti{\Lambda}_j = \sigma^N_{L_j}(m_j)$, $\ti{\Lambda}_{ij} =
\sigma^N_{L_{ij}}(m_i,m_j)$, and
$\ti{\Delta}_1=\ti{\Delta}_{T_{m_1}M_1}$.  We claim that
\begin{align}
d(\ti{\Lambda}_0\times^N\ti{\Lambda}_{12},\ti{\Lambda}_{01}^-\times^N\ti{\Lambda}_2^-)
& =
d(\ti{\Lambda}_0\times^N\ti{\Delta}_1\times^N\ti{\Lambda}_2,\ti{\Lambda}_{01}^-\times^N\ti{\Lambda}_{12}^-)
\nonumber\\ & =
d(\ti{\Lambda}_0\times^N\ti{\Lambda}_2,\ti{\Lambda}_{01}^-\circ^N\ti{\Lambda}_{12}^-)
. \label{claim2}
\end{align}
The first identity is Lemma \ref{insertdiag}.  To prove \eqref{claim2}
we begin by noting the transverse intersection
$\Lambda_{02}\pitchfork\Lambda_{0}\times\Lambda_{2}$.  We denote
$\ti{\Lambda}_{02}:=\ti{\Lambda}_{01}\circ^N\ti{\Lambda}_{12}$ (hence
$\ti{\Lambda}_{02}^-=\ti{\Lambda}_{01}^-\circ^N\ti{\Lambda}_{12}^-$)
and pick a path $\ti{\gamma}_{02}:[0,1]\to\Lag^N(T_{m_0}M_0^-\times
T_{m_2}M_2)$ from
$\ti{\gamma}_{02}(0)=\ti{\Lambda}_0\times^N\ti{\Lambda_2}$ to
$\ti{\gamma}_{02}(1)=\ti{\Lambda}_{02}^-$ whose crossing form with
$\Lambda_0\times\Lambda_2$ at $s=0$ is positive definite and hence by
Remark \ref{alt deg}
$$ d(\ti{\Lambda}_0\times^N\ti{\Lambda_2},\ti{\Lambda}_{02}^-) = -
I'(\gamma_{02},\Lambda_0\times\Lambda_2) .
$$
Here $I'$ denotes the Maslov index of a pair of paths (the
second one is constant), not counting crossings at the endpoints.
Next, fix a complement $L_{11}\in\Lag(T_{(m_1,m_1)}M_1\times M_1^-)$
of the diagonal.  Then both $(\Lambda_{02}\times^N L_{11})^T$ and
$\Lambda_{01}\times\Lambda_{12}$ are transverse to $T_{m_0}M_0\times
\Delta_1\times T_{m_2}M_2$ and their composition is $\Lambda_{02}$. By
Lemma~\ref{the one} below we then find a path $\gamma_{0112}$ and lift
it to
$\ti{\gamma}_{0112}:[0,1]\to\Lag^N(T_{(m_0,m_1,m_1,m_2)}M_0^-\times
M_1\times M_1^-\times M_2)$ from
$\ti{\gamma}_{0112}(0)=[\ti{\Lambda}_{02}\times^N\ti{L}_{11}]^T$ to
$\ti{\gamma}_{0112}(1)=\ti{\Lambda}_{01}\times^N\ti{\Lambda}_{12}$
whose composition $\circ(\gamma_{0112})=\Lambda_{02}$ is constant and
that has no crossings with $\Lambda_0\times\Delta_1\times\Lambda_2$
(by the transversality
$\gamma_{0112}\cap(\Lambda_0\times\Delta_1\times\Lambda_2)
=\Lambda_{02}\cap(\Lambda_0\times\Lambda_2)=\{0\}$).  Here the grading
of $\ti{L}_{11}$ is determined by continuation along this path.  Since
the composition $\circ(\gamma_{0112})$ is constant this continuation
yields
$$
\ti{\Lambda}_{02}
=\circ^N(\ti{\gamma}_{0112})
=\circ^N((\ti{\Lambda}_{02}\times^N\ti{L}_{11})^T)
= d(\ti{L}_{11},\ti{\Delta}_1^-)\cdot\ti{\Lambda}_{02} .
$$
Here we also used \eqref{property}, and we deduce that $d(\ti{L}_{11},\ti{\Delta}^-_1)=0$ mod $N$.
Furthermore, we fix a path
$\ti{\gamma}_{11}:[0,1]\to\Lag^N(T_{(m_1,m_1)}M_1^-\times M_1)$ 
from $\ti{\gamma}_{11}(0)=\ti{\Delta}_1$ to $\ti{\gamma}_{11}(1)=\ti{L}_{11}^-$ 
whose crossing form with $\Delta_1$ at $s=0$ is positive definite, and thus
$$
- I'(\gamma_{11},\Delta_1) 
= d(\ti{\Delta}_1,\ti{L}_{11}^-)
= d(\ti{L}_{11},\ti{\Delta}_1^-) =0
\quad\text{mod}\;N .
$$ 
Now the concatenated path
$(\ti{\gamma}_{02}\times\ti{\gamma}_{11})^T\#\ti{\gamma}_{0112}^-$ connects
$\ti{\Lambda}_0\times^N\ti{\Delta}_1\times^N\ti{\Lambda}_2$
to $\ti{\Lambda}_{01}^-\times^N\ti{\Lambda}_{12}^-$
with positive definite crossing form at $s=0$,
and \eqref{claim2} can be verified,
\begin{align*}
 d(\ti{\Lambda}_0\times^N\ti{\Delta}_1\times^N\ti{\Lambda}_2,\ti{\Lambda}_{01}^-\times^N\ti{\Lambda}_{12}^-) 
&= - I'((\gamma_{02}\times\gamma_{11})^T\#\gamma_{0112}^-,\Lambda_0\times\Delta_1\times\Lambda_2) \\
&= - I'(\gamma_{02},\Lambda_0\times\Lambda_2)
- I'(\gamma_{11},\Delta_1)
- I'(\gamma_{0112}^-,\Lambda_0\times\Delta_1\times\Lambda_2) \\
&= - I'(\gamma_{02},\Lambda_0\times\Lambda_2)
= d(\ti{\Lambda}_0\times^N\ti{\Lambda_2},\ti{\Lambda}_{02}^-) .
\end{align*}

\vspace{-5mm}
\end{proof}

\begin{lemma} \label{the one}
Let $V_0, V_1, V_2$ be symplectic vector spaces,
$\Lambda_{02}\subset V_0^-\times V_2$ a Lagrangian subspace, and denote by
$$\TD_{\Lambda_{02}}\subset \Lag(V_0^-\times V_1 \times V_1^-\times V_2) $$
the subset of Lagrangian subspaces
$\Lambda\subset V_0^- \times V_1 \times V_1^- \times V_2$ with 
$\Lambda\pitchfork (V_0\times \Delta_{V_1} \times V_2)=:\hat\Lambda_{02}$ and
$\pi_{02}(\hat\Lambda_{02})=\Lambda_{02}$.
Then $\TD_{\Lambda_{02}}$ is contractible.
\end{lemma}
\begin{proof}
We fix metrics on $V_0$, $V_1$, and $V_2$. 
Then we will construct a contraction $(\rho_t)_{t\in[0,1]}$, 
$\rho_t:\TD_{\Lambda_{02}}\to \TD_{\Lambda_{02}}$
with $\rho_0={\rm Id}$ and $\rho_1\equiv\Psi(\Lambda_{02} \times (\Delta_1)^\perp)$,
where $\Psi:V_0^- \times V_2 \times V_1 \times V_1^-\to V_0^- \times V_1 \times V_1^- \times V_2$
exchanges the factors. 
To define $\rho_t(\Lambda)$ we write
$\Lambda=\hat\Lambda_{02}\oplus\hat\Lambda_{11}$, 
where $\hat\Lambda_{11}$ is the orthogonal complement of 
$\hat\Lambda_{02}$ in $\Lambda$.
Now $\hat\Lambda_{02}$ is the image of 
$({\rm Id}_{V_0},i_1,i_1,{\rm Id}_{V_2}): \Lambda_{02}\to V_0^- \times V_1 \times V_1^- \times V_2$
for a linear map $i_1:\Lambda_{02}\to V_1$
and $\hat\Lambda_{11}$ is the image of 
$(j_0,{\rm Id}_{V_1}+j_1,-{\rm Id}_{V_1}+j_1,j_2): V_1\to V_0^- \times V_1 \times V_1^- \times V_2$
for linear maps $j_i:V_1\to V_i$.
One can check that
$$
\rho_t(\Lambda):=
\im\bigl({\rm Id}_{V_0},t\cdot i_1,t\cdot i_1,{\rm Id}_{V_2} \bigr)
\oplus
\im\bigl(t\cdot j_0,{\rm Id}_{V_1}+t^2\cdot j_1,-{\rm Id}_{V_1}+t^2\cdot j_1,t\cdot j_2 \bigr)
$$
is an element of $\TD_{\Lambda_{02}}$ for all $t\in[0,1]$ and defines a smooth contraction.
\end{proof}

\section{Floer cohomology}
\label{sec:HF}

The main content of this section is a review of the construction of
graded Floer cohomology for pairs of Lagrangian submanifolds 
in monotone and exact cases by Floer, Oh, and Seidel. 
In \ref{sequences} we then extend Floer cohomology to generalized Lagrangian correspondences, which in Section \ref{quilted FH} will be reformulated in terms of pseudoholomorphic quilts.

\subsection{Monotonicity}
\label{sec:monotone}

Let $(M,\omega)$ be a symplectic manifold.  Let $\J(M,\omega)$
denote the space of compatible almost complex structures on
$(M,\omega)$.  Any $J \in \J(M,\omega)$ gives rise to a complex
structure on the tangent bundle $TM$; the first Chern class $c_1(TM)
\in H^2(M,\Z)$ is independent of the choice of $J$.  
Throughout, we will use the following standing assumptions on all
symplectic manifolds:

\label{M ass}
\begin{itemize}
\item[\bf(M1):] $(M,\omega)$ is {\em monotone}, that is for some
$\tau\geq 0$
\begin{equation*}
[\omega] = \tau c_1(TM) .
\end{equation*} 
\item[\bf(M2):] If $\tau>0$ then $M$ is compact. If $\tau=0$ then 
$M$ is (necessarily) noncompact but satisfies ``bounded geometry'' 
assumptions as in \cite{se:bo}.
\end{itemize}

Note here that we treat the exact case $[\omega]=0$ as special case of
monotonicity (with $\tau=0$).  Next, we denote the index map by
$$ c_1: \pi_2(M) \to \Z, \ \ u \mapsto (c_1,u_*[S^2]) .$$
The {\em minimal Chern number} $N_M \in \N$ is the non-negative generator
of its image. 

Associated to a Lagrangian submanifold $L\subset M$ 
are the Maslov index and action (i.e.\ symplectic area) maps
$$ I: \pi_2(M,L) \to \Z, \qquad A: \ \pi_2(M,L) \to \R . $$ 
Our standing assumptions on all Lagrangian submanifolds are the
following:

\begin{itemize}
\setlength{\itemsep}{1mm}
\item[\bf(L1):]
\label{L ass}
$L$ is {\em monotone}, that is
$$ 2 A(u) = \tau I(u) \quad \forall u \in \pi_2(M,L) $$
where the $\tau\geq 0$ is (necessarily) that from (M1).
\item[\bf(L2):]
$L$ is compact and oriented.
\end{itemize}

Any homotopy class $[u]\in\pi_2(M,L)$ that is represented by a nontrivial $J$-holomorphic disk $u:(D,\partial D) \to (M,L)$ 
has positive action $A([u])=\int u^*\omega>0$. 
Monotonicity with $\tau>0$ then implies that the index is also positive.
So, for practical purposes, we define the (effective) {\em minimal Maslov number} 
$N_{L}\in\N$ as the generator of $I(\{[u]\in\pi_2(M,L)|A([u])>0\}) \subset \N$.
If $M$ and $L$ are exact ($\tau=0$), then $A\equiv 0$, so we have $N_{L}=\infty$.

If the Lagrangian submanifold $L$ is oriented then $I(u)$ is always even since it is the Maslov index of a loop of oriented Lagrangian subspaces.
So the orientation and monotonicity assumption on $L$ imply $N_L\geq 2$, i.e.\ 
any nontrivial holomorphic disk must have $I(u)\geq 2$, which excludes disk bubbling in 
transverse moduli spaces of index $0$ and $1$. 
%KW removed this since we need orientation to have index of general surfaces fixed mod 2 by the ends
%
%For most purposes\footnote{
%For nonoriented Lagrangians or Lagrangian correspondences 
%all constructions and results in this paper extend directly
%to the ungraded Floer cohomologies with $\Z_2$-coefficients. 
%We do not discuss orientations of moduli spaces in this case.
%The grading also directly extends if we drop the assumption (G2)
%of compatibility with orientations.
%}
%we could drop the orientation assumption and replace (L2) by

%\begin{itemize}
%\item[\bf(L2'):]
%$L$ is compact and has minimal Maslov number $N_L\geq 2$.
%\end{itemize}
%
In order for the Floer cohomology groups to be well defined
we will also have to make the following additional assumption.

\begin{itemize}
\item[\bf(L3):] 
 $L$ has minimal Maslov number $N_L\geq 3$.
\end{itemize}

Alternatively, following \cite{oh:fl1,fooo},
we may replace (L3) by an assumption on the disk counts, which we briefly recall here.
Given any $p\in L$ and  $\omega$-compatible almost complex structure $J\in \J(M,\omega)$
let $\M_1^2(L,J,p)$ be the moduli space of $J$-holomorphic disks $u:(D,\partial D) \to (M,L)$ with Maslov number $2$ and one marked point satisfying $u(1)=p$, modulo automorphisms of the disk fixing $1 \in \partial D$.  
Oh proves that for any $p\in L$ there exists a dense subset $\J^{\rm reg}(p) \subset \J(M,\omega)$ such that $\M_1^2(L,J,p)$ is cut out transversely, and consequentially a finite set. Moreover, any relative spin structure on $L$ induces an orientation on $\M_1^2(L,J,p)$.
Letting $\eps: \M_1^2(L,J,p) \to \{ \pm 1 \} $ denote the map
comparing the given orientation to the canonical orientation of a
point, the disk number of $L$,
\begin{equation}  \label{diskcount}
w(L) := \sum_{u \in \M_1^2(L,J,p)} \eps(u) 
\end{equation}
is independent of $J \in J^{\rm reg}(p)$ and $p \in L$. If we work without orientations then $w(L)\in\Z_2$ is still well defined. 

As a third class of assumptions, 
we will restrict our considerations to Maslov coverings and gradings
that are compatible with orientations, that is we make
the following additional assumptions on the grading of the symplectic manifolds $M$
and Lagrangian submanifolds $L\subset M$.
(In the case $N=2$ these assumptions reduce to (L2).)
\label{G ass}
\begin{itemize}
\item[\bf(G1):]
$M$ is equipped with a Maslov covering $\Lag^N(M)$ for $N$ even, 
and the induced $2$-fold Maslov covering $\Lag^2(M)$ is the one 
described in Example \ref{maslovexamples} (i).
\item[\bf(G2):]
$L$ is equipped with a grading $\sigma_L^N:L\to\Lag^N(M)$, and
the induced $2$-grading $L\to\Lag^2(M)$ is the one given by the orientation of $L$.
\end{itemize}

In the following we discuss topological situations which ensure monotonicity.

\begin{lemma}  \label{monotone}
Suppose that $M$ is monotone and $L \subset M$
is a compact Lagrangian such that $\pi_1(L)$ is torsion; 
that is every element has finite order, and there is a finite maximal order $k$. 
Then $L$ is monotone and the minimal
Maslov number is at least $2N_M/k$.
\end{lemma}

\begin{proof}   Let $u:(D,\partial D) \to (M,L)$ and let $k(u)$ be
the order of the restriction of $u$ to the boundary in $\pi_1(L)$.
After passing to a $k(u)$-fold branched cover $\ti{u}$, we may assume that the
restriction of $\ti{u}$ to $\partial D$ is homotopically trivial in $L$. 
By adding the homotopy we obtain a sphere $v: S^2 \to M$ with
$ k(u) I(u) = I(\ti{u}) = 2 c_1(v) $ divisible by $2N_M$.  
For the relation between the first Chern class and the Maslov index
see e.g.\ \cite[Appendix~C]{ms:jh}.
The similar identity for the actions (due to $\omega|_L=0$) 
completes the proof.
\end{proof}

In practice, we will need the action-index relation not only for disks as in (L1) but also for other surfaces with several boundary components mapping to several Lagrangians. (This really only becomes relevant in \cite{quilts} for the definition of relative invariants from surfaces with strip-like ends.)
In particular, to define Floer cohomology for a pair of Lagrangians (and especially later to prove the isomorphism \eqref{eq:iso}) we need the action-index for annuli between the two Lagrangians. This provides the energy-index relation in Remark~\ref{rmk:monotone}. In fact, it also implies monotonicity (L1) for both Lagrangians as long as $M$ is connected.

\begin{definition}  \label{tuple monotone} \begin{enumerate}
\item
We say that a tuple $(L_e)_{e\in\cE}$ is {\em monotone} 
with monotonicity constant $\tau\ge 0$ if the
following holds: Let $\Sigma$ be any connected compact surface with 
nonempty boundary $\partial\Sigma=\sqcup_{e\in\cE} C_e$ 
(with $C_e$ possibly empty or disconnected).  
Then for every map $u :\Sigma \to M$ satisfying
$u(C_e)\subset L_e$ we have the action-index relation
$$ 2\int u^*\omega = \tau \cdot I(u^*TM,(u^*TL_e)_{e\in\E}),$$
where $I$ is the sum of the Maslov indices of the totally real subbundles
$(u|_{C_e})^*TL_e$ in some fixed trivialization of $u^*TM$.
\item
We say that a pair $(L_0,L_1)$ is {\em monotone for Floer theory} if (a) holds for the annulus $\Sigma=[0,1]\times S^1$ and every map $u$ with boundary values $u(\{j\}\times S^1)\subset L_j$ for $j=0,1$.
\end{enumerate}
\end{definition}

The following is a minor generalization of \cite[Proposition 2.7]{oh:fl1}.

\begin{lemma}  \label{monotone2}  
Suppose that $M$ is monotone.
\begin{enumerate}
\item
If each $L_e\subset M$ is monotone, and the image of each 
$\pi_1(L_e)$ in $\pi_1(M)$ is torsion, then the tuple $(L_e)_{e\in\E}$ 
is monotone.
\item
If both $L_0,L_1\subset M$ are monotone, and the image of 
$\pi_1(L_0)$ or $\pi_1(L_1)$ in $\pi_1(M)$ is torsion, then the pair $(L_0,L_1)$ is monotone for Floer theory. 
\end{enumerate}
\end{lemma}

\begin{proof} 
To check (a) consider $u :\Sigma \to M$ satisfying $u(C_e)\subset L_e$.  By
assumption we have integers $N_e\in\N$ such that $N_e u|_{C_e}$ is
contractible in $M$.  Let $N=\prod_{e\in \E} N_e$, so that $N
u|_{C_e}$ is contractible for all boundary components $C_e$ of
$\Sigma$.  Let $\ti{\Sigma} \to \Sigma$ be a ramified $N$-cover defined as follows:
Pick one ramification point $z_1,\ldots,z_k\in\Sigma$ in each connected component of $\Sigma$ with nonempty boundary. Then there exists a representation 
$\rho:\pi_1(\Sigma \setminus \{z_1,\ldots,z_k\}) \to \Z_N$ with $\rho([C_e])=[N/N_e]$. 
The induced ramified cover $\rho:\ti{\Sigma} \to \Sigma$ satisfies the following:
The inverse image $\ti{C}_e$ of $C_e$ consists of $N/N_e$ connected components, 
each of which is an $N_e$-fold cover of $C_e$.
Now the pull-back $\ti{u}:
\ti{\Sigma} \to M$ of $u:\Sigma\to M$ has restrictions to the boundary
$\ti{u} |_{\ti{C}_e}$ that are homotopically trivial in $M$.  Thus
$\ti{u}$ is homotopic to the connected sum of some maps $v_{e,j}: (D,\partial D)
\to (M,L_e)$ for $j=1,\ldots,N/N_e$ and a map $v:S\to M$ on a closed surface $S$.
We can now use the closedness of $\omega$ and the monotonicity of $M$
and each $L_e$ to deduce
\begin{align*}
2 N \int_\Sigma u^*\omega \;=\; 2 \int_{\ti\Sigma} \tilde u^*\omega 
&\;=\; 2 \int_S v^*\omega + \sum_{e\in \E, j=1,\ldots,N/N_e} 2 \int_D v_{e,j}^*\omega \\
&\;=\; 2\tau c_1(v^*TM) + \sum_{e\in \E, j=1,\ldots,N/N_e} \tau I(v_{e,j})  
\;=\; \tau I(\tilde u) \;=\; \tau N I(u) ,
\end{align*}
using properties of the Maslov index explained in \cite[Appendix~C]{ms:jh}. 
The first equality $N \int u^*\omega = \int \tilde u^*\omega$ can be confirmed by integrating over the complement of the branch points.
The last equality $I(\tilde u) = N I(u)$ holds since any trivialization $u^*TM\cong\Sigma\times\C^n$ induces a trivialization $\tilde u^*TM\cong\ti\Sigma\times\C^n$, which restricts to the $N$-fold covering $\rho\times{\rm Id}:\ti C_e \times\C^n \to C_e \times\C^n$ on each boundary component.

In the case of (b) we can take a multiple cover of the annulus such that one boundary loop is contractible in $M$, and hence the multiply covered annulus is homotopic to two disks to which we can apply monotonicity of the single Lagrangians.
\end{proof}

In the exact case, with $\omega=d\lambda$, any tuple of exact
Lagrangians $(L_e)_{ e \in \E}$, that is with $[\lambda|_{L_e}]=0\in
H^1(L_e)$, is automatically monotone.  Moreover, note that
monotonicity is invariant under Hamiltonian isotopies of one or several
Lagrangians.

\begin{remark}\label{BS}
Another situation in which one naturally has monotonicity is the
Bohr-Sommerfeld setting, as pointed out to us by P. Seidel.  Suppose
that the cohomology class $[\omega]$ is integral.  Let $(\cL,\nabla)
\to (M,\omega)$ be a unitary line-bundle-with-connection having
curvature $(2\pi /i) \omega$.  The restriction of $(\cL,\nabla)$ to
any Lagrangian $L \subset M$ is flat.  $L$ is {\em Bohr-Sommerfeld} if
the restriction of $(\cL,\nabla)$ to $L$ is trivial, that is, there
exists a non-zero horizontal section.  In that case, we choose a horizontal section
$\phi^\cL_L$ of unit length, which is unique up to a collection of phases 
$U(1)^{\pi_0(L)}$.
Suppose that $M$ is monotone, $[\omega]= \lambda c_1(M)$ for some
$\lambda > 0$.  Since $c_1(M)$ and $[\omega]$ are integral, we must have 
$\lambda = k/l$ for some integers $k,l> 0$.  Let $\cK^{-1} \to M$ denote the
anticanonical bundle, $\cK_m^{-1} = \Lambda^{\top}_\C (T_m^{0,1} M )$,
which satisfies $ k c_1(\cK^{-1}) = l \frac i{2\pi}[\curv(\nabla)] = l c_1(\cL) $. 
Hence there exists an isomorphism
$$ \Phi: (\cK^{-1})^{\otimes k} \to \cL^{\otimes l} . $$
Let $L\subset M$ be an oriented Lagrangian submanifold.  
The restriction of $\cK^{-1}$ to $L$ has a natural non-vanishing section 
$\phi_L^\cK$ given by the orientation and the isomorphisms
$$ \Lambda^{\top}_\R TL \to  \Lambda^{\top}_\C T^{0,1}M|_L, \ \ \
v_1 \wedge \ldots \wedge v_n \mapsto (v_1 + iJ v_1) \wedge \ldots
\wedge (v_n + iJ v_n) .$$
We say that $L$ is {\em Bohr-Sommerfeld monotone} with respect to
$(\cL,\nabla,\Phi)$ if it is Bohr-Sommerfeld and
the section $(\phi^{\cL}_L)^{\otimes l}$ is homotopic to 
$\Phi \circ (\phi^{\cK}_L)^{\otimes k}$, that is, there 
exists a function $\psi:L \to \R$ such that
\begin{equation} \label{stupid}
 (\exp( 2 \pi i \psi ) \phi^{\cL}_L)^{\otimes l} = \Phi \circ
(\phi^\cK_L)^{\otimes k} .\end{equation}
%
%K from referee notes
%Note that in case .........
%any Bohr-Sommerfeld Lagrangian will be monotone with respect to some connection 
%.................
%In the following applications, however, we will need the notion of Bohr-Sommerfeld monotone with respect to a fixed connection $\nabla$.
\end{remark}

\begin{lemma} \label{BS monotone} 
Let $(L_e)_{e \in \E}$ be a collection of Lagrangians such
that each is Bohr-Sommerfeld monotone with respect to
$(\cL,\nabla,\Phi)$.  Then $(L_e)_{e\in\cE}$ is monotone.
\end{lemma}
 
\begin{proof}  Let $\Sigma$ be a compact Riemann surface with boundary
components $(C_e)_{e \in \E}$.  Let $u: \Sigma \to M$ be a map with
boundary $u(C_e)\subset L_e$.  The index $I(u)$ is the sum of Maslov
indices of the bundles $ (u|_{C_e})^* TL_e$, with respect to some
fixed trivialization of $u^* TM$.  Equivalently, $I(u)$ is the sum of
winding numbers of the sections $\phi_{L_e}^{\cK}$ with respect to the
induced trivialization of $u^* \cK^{-1}$.  Since each $L_e$ is
Bohr-Sommerfeld, $k I(u)$ is the sum of the winding numbers of the
sections $(\phi_{L_e}^{\cL})^{\otimes l}$, with respect to the induced
trivialization of $u^* \cL^{\otimes l}$.  Write $u^* \nabla^{\otimes
l} = \d + \alpha$ for some $\alpha \in \Omega^1(\Sigma)$ in this
trivialization, so that $u^* \curv(\nabla^{\otimes l}) = \d \alpha$.
Since the sections are horizontal, we have
$$ k I(u) = (i/ 2\pi) \int_{\partial \Sigma} \alpha = (i /2\pi) \int_\Sigma
u^* \curv(\nabla^{\otimes l}) = l A(u) .$$
\end{proof}

\subsection{Graded Floer cohomology for pairs of Lagrangians}
\label{standard Floer}

Let $L_0,L_1 \subset M$ be compact Lagrangian submanifolds.  For a
time-dependent Hamiltonian $H \in C^\infty([0,1] \times M)$ let
$(X_t)_{t\in[0,1]}$ denote the family of Hamiltonian vector fields for
$(H_t)_{t\in[0,1]}$, and let $\phi_{t_0,t_1}: M \to M$ denote its
flow.  (That is, $\phi_{t_0,t_1}(y)=x(t_1)$, where $x:[0,1]\to M$
satisfies $\dot x=X_t(x)$, $x(t_0)=y$.)  We will abbreviate
$\phi_1:=\phi_{0,1}$ for the time 1 flow from $t_0=0$ to $t_1=1$.  Let
$\Ham(L_0,L_1)$ be the set of $H \in C^\infty([0,1] \times M)$ such
that $\phi_1(L_0)$ intersects $L_1$ transversally.  Then we have a
finite set of perturbed intersection points
$$
\cI(L_0,L_1) := \bigl\{ \gamma: [0,1] \to M \,\big|\, \gamma(t) 
= \phi_{0,t}(\gamma(0)), \ \gamma(0) \in L_0, \  \gamma(1) \in L_1 \bigr\} .
$$ 
It is isomorphic to the intersection $\phi_1(L_0)\pitchfork L_1$.
If we assume that $M$ and $L_0,L_1$ are graded as in (G1-2), then we obtain
a degree map from Section \ref{gradings},
$$  \cI(L_0,L_1) \to \Z_N ,
\qquad x\mapsto |x|=d(\sigma_{L_0}^N(x),\sigma_{L_1}^N(x)) . $$ 
Since $N$ is even the sign $(-1)^{|x|}$ is well-defined.  It agrees
with the usual sign in the intersection number, given by the
orientations of $\phi_1(L_0)$ and $L_1$, which also determine the mod
$2$ grading by assumption.

Next, we denote the space of time-dependent $\omega$-compatible 
almost complex structures by
$$
\J_t(M,\omega) := C^\infty([0,1], \J(M,\omega)) .
$$ 
For any $J\in\J_t(M,\omega)$ and $H\in\Ham(L_0,L_1)$ we say that a map 
$u: \R \times [0,1] \to M$ is {\em $(J,H)$-holomorphic} 
with Lagrangian boundary conditions if
\begin{equation} \label{Floer}
\overline{\partial}_{J,H} u :=
\partial_s u({s},t) + J_{t,u({s},t)}(\partial_tu({s},t) -
X_t(u({s},t))) = 0 ,\end{equation}
\begin{equation} \label{Lag bc}
u(\R,0) \subset L_0, \quad u(\R,1)\subset L_1 .
\end{equation}
The (perturbed) energy of a solution is
$$
E_H(u):=\int_{\R\times[0,1]} |\pd_s u|^2 = 
\int_{\R\times[0,1]} u^*\omega + \d (H(u) \d t) . 
$$
The following exponential decay lemma of Floer \cite{Floer:unregularized}
will be needed later and is part of the proof of Theorem \ref{traj} below.

\begin{lemma} \label{expconv} 
Let $H\in\Ham(L_0,L_1)$ and $J\in\J_t(M,\omega)$. 
Then for any $(J,H)$-holomorphic strip $u:\R \times [0,1] \to M$ 
with Lagrangian boundary conditions in $L_0,L_1$
the following are equivalent:
\begin{enumerate}
\item  $u$ has finite energy $E_H(u)=\int_{\R \times [0,1]} | \partial_s u |^2 < \infty $;
\item  There exist $x_\pm \in \cI(L_0,L_1)$ such that 
$u(s,\cdot)$ converges to $x_\pm$ exponentially in all derivatives 
as $ s \to \pm \infty$ .
\end{enumerate}
\end{lemma} 

For any $x_\pm \in \cI(L_0,L_1)$ we denote by 
$$
\M(x_-,x_+):= \bigl\{ u:\R \times [0,1] \to M \,\big|\, \eqref{Floer}, \eqref{Lag bc},
E_H(u)<\infty, \lim_{{s} \to \pm \infty} u(s,\cdot) = x_\pm
\bigr\} / \R
$$
the space of finite energy $(J,H)$-holomorphic maps modulo translation
in $s\in\R$.  It is isomorphic to the moduli space of finite energy
$J'$-holomorphic maps with boundary conditions in $\phi_1(L_0)$ and
$L_1$, and without Hamiltonian perturbation.  Here
$J'\in\J_t(M,\omega)$ arises from $J$ by pullback with $\phi_{t,1}$.

\begin{remark} \label{rmk:monotone}
Suppose that the pair $(L_0,L_1)$ is monotone, then for any $x_\pm \in
\cI(L_0,L_1)$ there exists a constant $c(x_-,x_+)$ such that for all
$u\in\M(x_-,x_+)$ the energy-index relation holds:
\begin{equation}\label{energy index}
2 E_H(u) =  \tau\cdot {\rm Ind}(D_u)
 + c(x_-,x_+) ,
\end{equation}
where $D_u$ denotes the linearized operator at $u$ of the Cauchy-Riemann equation
\eqref{Floer} on the space of sections of $u^*TM$ satisfying the linearized Lagrangian 
boundary conditions from \eqref{Lag bc}. 
Its Fredholm index is given by the Maslov-Viterbo index of $u$. 
This monotonicity ensures energy bounds for the moduli spaces of fixed index and thus compactness up to bubbling.
\end{remark}

\begin{theorem} \label{traj}  (Floer,Oh) 
Let $L_0,L_1 \subset M$ be a monotone pair\footnote{Throughout, we are working with monotone pairs of Lagrangians in the sense of Definition~\ref{tuple monotone}.} 
of Lagrangian submanifolds
satisfying (L1-2) and (M1-2).
For any $H \subset \Ham(L_0,L_1)$ there exists a dense subset
$\J_t^{\reg}(L_0,L_1;H) \subset \J_t(M,\omega)$ such that the following holds for all $x_\pm\in\cI(L_0,L_1)$.
\ben 
\item 
$\M(x_-,x_+)$ is a smooth manifold whose dimension near a nonconstant solution $u$
is given by the formal dimension ${\rm Ind}(D_u) -1$.
We denote $\M(x_-,x_+)_j:=\{{\rm Ind}(D_u) = j+1\}$; thus excluding the constant solution from $\M(x,x)_0$.
\item 
The component $\M(x_-,x_+)_0 \subset \M(x_-,x_+)$ 
of formal dimension zero is finite.
\item 
Suppose that $L_0$ and $L_1$ have minimal Maslov numbers $N_{L_k}\geq 3$.  
Then the one-dimensional component $\M(x_-,x_+)_1 \subset \M(x_-,x_+)$ has
a compactification as one-dimensional manifold with boundary
\begin{equation} \label{glue}
\partial \overline{\M(x_-,x_+)_{1}} \cong \bigcup_{x \in \cI(L_0,L_1)}
\M(x_-,x)_{0} \times \M(x,x_+)_{0} \end{equation}
\item If $(L_0,L_1)$ is relatively spin (as defined in e.g.\ \cite{orient}), 
then there exists a coherent set of orientations on
$\M(x_-,x_+)_0,\M(x_-,x_+)_1$ for all $x_\pm \in \cI(L_0,L_1)$, that is,
orientations compatible with \eqref{glue}.
\een
\end{theorem}

For the proofs of (a-c) we refer to Oh's paper \cite{oh:fl1} and the
clarifications \cite{oh:st}, \cite{ohkw:st}. For the exact case see
\cite{se:bo}. The proof of (d) is contained in 
\cite{orient} loosely following \cite{fooo}.  From (d) we obtain a map
$$ \eps: \M(x_-,x_+)_0 \to \{ \pm 1 \} $$
defined by comparing the given orientation to the canonical
orientation of a point. 

Now let $M$ be a monotone symplectic manifold satisfying (M1-2) and equipped
with an $N$-fold Maslov covering.
Let $L_0,L_1\subset M$ be a monotone, relative spin pair of
graded Lagrangian submanifolds satisfying (L1-3), and let $H \in \Ham(L_0,L_1)$.
The {\em Floer cochain group} is the $\Z_N$-graded group
$$ CF(L_0,L_1) = \bigoplus_{d \in \Z_N} CF^d(L_0,L_1), \qquad
CF^d(L_0,L_1) =\bigoplus_{x\in\cI(L_0,L_1), |x| = d} \Z \bra{x} , $$
and the {\em Floer coboundary operator} is the map of degree $1$,
$$ \partial^d : \ CF^d(L_0,L_1) \to CF^{d+1}(L_0,L_1) ,$$
defined by
$$
\partial^d\bra{x_-} := \sum_{x_+\in\cI(L_0,L_1)}
\Bigl( \sum_{u \in\M(x_-,x_+)_0} \eps(u)\Bigr) \bra{x_+} .$$
Here we choose some $J \in \J_t^\reg(L_0,L_1;H)$.  If an isolated
trajectory $u \in\M(x_-,x_+)_0$ exists, then the degree identity
$|x_+|=|x_-|+1$ can be seen by concatenating the paths
$\ti{\gamma_0},\ti{\gamma_1}$ of graded Lagrangians in the definition
of $|x_-|$ with the unique graded lifts of $u^*TL_0,u^*TL_1$ to obtain
paths of graded Lagrangians defining $|x_+|$ (using a trivialization
of $u^*TM$ over the strip, compactified to a disk).  By additivity of
the Maslov index this shows $|x_+|=|x_-|+I(u^*TL_0,u^*TL_1)=|x_-|+1$.
It follows from Theorem \ref{traj} that $\partial^2 = 0$.
Alternatively, following \cite{oh:fl1}, we could drop assumption (L3), then
$\partial^2=(w(L_0)-w(L_1)){\rm Id}$, and to obtain a well defined cohomology it
suffices to assume that this disk count vanishes.
%
% Note: homology has natural coproducts, cohomology products
% categories have products, cocategories have coproducts }
%
In either case, the {\em Floer cohomology}
\footnote{Note that our conventions differ from Seidel's definition of graded 
Floer cohomology in \cite{se:gr} in two points which cancel each other:
The roles of $x_-$ and $x_+$ are interchanged and we switched 
the sign of the Maslov index in the definition of the degree \eqref{the degree}.}
$$ HF(L_0,L_1) := \bigoplus_{d\in\Z_N} HF^d(L_0,L_1), \qquad
HF^d(L_0,L_1) := {\ker(\partial^d)}/{\on{im}}(\partial^{d - 1}) $$
is $\Z_N$-graded.
It is independent of the choice of $H$ and $J$; a generalization of this
fact is proved in Section~\ref{sec:invariance} below.
If the gradings moreover satisfy (G1-2), then we have a well defined 
splitting
$$ HF(L_0,L_1) = HF^{\even}(L_0,L_1) \oplus HF^{\odd}(L_0,L_1) ,
$$
which coincides with the splitting induced by the orientations of $L_0,L_1\subset M$.

%\begin{remark}  
%The Floer complex $CF(L_1,L_0)$ for the switched pair is 
%chain homotopic  to the dual complex $\Hom(CF(L_0,L_1),\Z)$.
%Hence the Floer cohomology $HF(L_1,L_0)$ is isomorphic to the direct sum of
%$\Hom(HF(L_0,L_1),\Z)$ with \marginpar{universal coefficient theorem???}
%, see \cite{fieldb} for details.
%\end{remark}
%
%  
%  \begin{remark} \label{flat} 
%  If $L_0,L_1$ admit flat vector bundles $E_0,E_1$, one can extend the
%  definition to {\em Floer cohomology with values in $E_0,E_1$} by
%  defining
%  %
%  $$ CF^*(L_0,L_1;E_0,E_1) := \bigoplus_{x\in\cI(L_0,L_1)}
%  E_{0,x}^\dual \otimes E_{1,x}.$$
%  %
% \comment{K: deleted bullets but might want to put degrees in}
%  The differential then takes into account the parallel transport maps
%  along the boundary of $u \in \M(x_-,x_+)_0$.
%  \end{remark} 
%  

\subsection{Floer cohomology for generalized Lagrangian correspondences}
\label{sequences}

The goal of this section is to define a first version of Floer cohomology for a cyclic generalized Lagrangian correspondence $\ul{L}$ as in Definition \ref{Lag cycle}.  
So we consider $\ul{L}=(L_{01},\ldots, L_{r(r+1)})$, a sequence of smooth Lagrangian correspondences $L_{(j-1)j} \subset M_{j-1}^-\times M_{j}$ between a sequence 
$M_0,M_1,\ldots,M_{r+1}=M_0$ of symplectic manifolds.
For example, we could consider a non-cyclic sequence of Lagrangians $L_{01}\subset M_1$, $\bigl(L_{(i-1) i} \subset M_{i-1}^-\times M_i\bigr)_{i=2,\ldots,r}$,
$L_{r0} \subset M_{r}^-$, which is a special case of the cyclic setup with $M_0=\{\pt\}$. 
The usual Floer cohomology for pairs of Lagrangians fits into this case with $r=1$ and Lagrangian submanifolds $L_{01}, L_{10}\subset M_1$.

We assume that $\ul{L}$ satisfies (M1-2,L1-3), i.e.\ each $M_j$ satisfies (M1-2) 
and each $L_{(j-1)j}$ satisfies (L1-3) with a fixed monotonicity constant $\tau\ge0$.
We moreover assume that $\ul{L}$ is graded in the sense of Definition~\ref{gen grad c}  and equipped with a relative spin structure in the following sense.
Alternatively, we may replace the minimal Maslov assumption (L3) by the assumption 
that the sum of disk numbers from \eqref{diskcount} vanishes,
\begin{equation}\label{w=0}
w(L_{01})+ \ldots + w(L_{r(r+1)}) = 0.
\end{equation}

\begin{definition} \label{gen spin c}
Let $\ul{L}=(L_{01},\ldots, L_{r(r+1)})$ be a cyclic generalized
Lagrangian correspondence (i.e.\ $L_{j(j+1)} \subset M_{j}^-\times
M_{j+1}$ for a cyclic sequence $M_0, M_1, \ldots, M_{r+1}=M_0$ of
symplectic manifolds).  A {\em relative spin structure} on $\ul{L}$
consists of a collection of background classes $b_j \in H^2(M_j,\Z_2)$
for $j=0,\ldots,r+1$ and relative spin structures on $L_{j(j+1)}$ with
background classes $- \pi_{j}^* b_{j} + \pi_{j+1}^* b_{j+1}$.  The
cyclic requirement on the background classes $b_0\in H^2(M_0,\Z_2)$
and $b_{r+1}\in H^2(M_{r+1},\Z_2)=H^2(M_0,\Z_2)$ is $b_{r+1}=b_0$ for
$r$ odd and $b_{r+1}=b_0+w_2(M_0)$ for $r$ even.\footnote{ This shift
is necessary in order to fit in the canonical relative spin structure
for the diagonal $\Delta_0$, see \cite{cat} for details.}
%
%
%The {\em background class} of the Lagrangian spin structure
%is the class of $(b_0,b_1\ldots,b_r)$ in $H^2(\prod M_j,\Z_2)$.
%
\end{definition}

Eventually, in Section \ref{quilted FH}, we will define the Floer cohomology $HF(\ul{L})$ directly,  using ``quilts of pseudoholomorphic strips''. In this section however we define $HF(\ul{L})$ as a special case of the Floer cohomology for pairs of Lagrangian submanifolds -- which are constructed from the sequence $\ul{L}$ as follows.
If $\ul{L}$ has even length $r+1\in 2\N$ we define a pair of graded Lagrangian submanifolds,
$$
\begin{aligned}
{L}_{(0)}&:= ( L_{01}\times L_{23} \times \ldots \times L_{(r-1)r} )  \\
{L}_{(1)}&:= ( L_{12} \times L_{34} \times \ldots \times L_{r(r+1)})^T  
\end{aligned}
\quad\subset\; M_0^-\times M_1 \times M_2^-\times\ldots \times M_{r} =: \widetilde{M} .
$$
Here we denote by $M_1^-\times \ldots \times M_{r}^-\times M_0 \to 
M_0^-\times M_1\times \ldots\times M_{r}$, $Z\mapsto Z^T$ 
the transposition of the last to the first factor, combined with an
overall sign change in the symplectic form.
If $\ul{L}$ has odd length $r+1\in 2\N +1$ we insert the diagonal  
$\Delta_0\subset M_{0}^-\times M_0 = M_{r+1}^-\times M_0$ (equipped with its
canonical grading) into $\ul{L}$
before arranging it into a pair of Lagrangian submanifolds as above, yielding
$$
\begin{aligned}
{L}_{(0)}&= ( L_{01}\times L_{23} \times \ldots \times L_{r(r+1)} )   \\
{L}_{(1)}&= ( L_{12} \times L_{34} \times \ldots \times L_{(r-1)r}\times \Delta_{0})^T  
\end{aligned}
\quad\subset\;
M_0^-\times M_1 \times\ldots \times M_r^-\times M_{r+1} = \widetilde{M}
$$
In the case of a noncyclic correspondence with $M_0=M_{r+1}=\{\pt\}$ 
the transposition as well as insertion of the diagonal are trivial operations.
Note that, beyond the grading, also the 
monotonicity, compactness, and orientation assumptions (L1-2) on $\ul{L}$ 
transfer directly to properties (L1-2) for ${L}_{(0)}$ and ${L}_{(1)}$.
Similarly, a relative spin structure on $\ul{L}$ induces compatible relative spin 
structures on ${L}_{(0)}$ and ${L}_{(1)}$, see \cite{orient}.
Moreover, we say that $\ul{L}$ is {\em monotone} if the pair of Lagrangians
$({L}_{(0)},{L}_{(1)})$ is monotone in the sense of Definition~\ref{tuple monotone}(b).
If this is the case, then a graded Floer cohomology for $\ul{L}$ can be defined by
$$ HF(\ul{L}) := HF({L}_{(0)} , {L}_{(1)} ) .
$$

\begin{remark} \label{rmk split J}
To see that $HF({L}_{(0)} , {L}_{(1)} )$ is well defined we need to make sure that $\partial^2=0$.
This holds immediately if ${L}_{(0)}$ and ${L}_{(1)}$ also satisfy (L3), if the bubbling of holomorphic discs is otherwise excluded, or if the effect of bubbling sums up to zero.
This can be achieved if all Lagrangians satisfy (L3) or, weaker, if the total disk count vanishes \eqref{w=0}.
\begin{enumerate}
\item 
Note that the assumption (L3) on the factors of $\ul{L}$ does not directly transfer to the product Lagrangians ${L}_{(0)}$ and ${L}_{(1)}$ since a difference of Maslov numbers greater than $3$ could give a total Maslov number less than $3$.
However, if we use a split almost complex structure $\widetilde{J}=J_0\oplus\ldots\oplus J_r$ on $\widetilde{M}$, induced from compatible almost complex structures $J_k$ on each $M_k$, then any nonconstant holomorphic disc in $\widetilde{M}$ with boundary on ${L}_{(0)}$ or ${L}_{(1)}$ will simply be a product of $J_k$-holomorphic discs. Pairs of these discs take boundary values in the Lagrangian correspondences $L_{(k-1)k}$ which satisfy the monotonicity assumptions as well as (L3). Hence each of these double discs must have nonnegative area and hence index, and at least one of them has positive area and hence Maslov index at least $3$. This excludes bubbling in moduli spaces of index $1$ or $2$, hence proves $\partial^2=0$.

The proof that transversality can be achieved with an almost complex structure (and also Hamiltonian perturbation) of split type can be found in Theorem \ref{quilt traj} and Proposition \ref{pert split} below. This excludes bubbling such that $\partial^2=0$ for this specific choice of perturbation data (and hence for any other choice of regular perturbation data). So the Floer cohomology $HF({L}_{(0)} , {L}_{(1)} )$ is indeed well defined.
\item
In the absence of (L3) we have $\partial^2=(w(L_{(0)})- w(L_{(1)})){\rm Id}$ by \cite{oh:fl1}.
Using a split almost complex structure we show in \cite{fieldb} (here stated in the case of even length $r+1$) that 
$w({L}_{(0)})= w(L_{01})+ w(L_{23}) + \ldots + w(L_{(r-1)r})$
and 
$w({L}_{(1)})= - w(L_{12}) - w(L_{34}) - \ldots - w(L_{r(r+1)})$,
hence
$\partial^2=(w(L_{01}) + \ldots + w(L_{r(r+1)})){\rm Id}$, which vanishes if we assume \eqref{w=0}.
The relative minus sign in $w({L}_{(1)})$ arises from the fact that e.g.\ $w(L_{12})$ is the disk count for the Lagrangian $L_{12}\subset M_1^-\times M_2$, whereas in the construction of $L_{(1)}$ we use $L_{12}^- \subset M_1\times M_2^-$, the same submanifold but viewed as Lagrangian with respect to the reversed symplectic structure.
The disk counts are related by $w(L_{12}^-) = - w(L_{12})$, since the $(J_1,-J_2)$-holomorphic discs with boundary on $L_{12}^-$ are identified with $(-J_1,J_2)$-holomorphic discs with boundary on $L_{12}$ via a reflection of the domain, which is orientation reversing for the moduli spaces.
\end{enumerate}
\end{remark}

In the case of a non-cyclic sequence the Floer cohomology $HF(\ul{L})$ specializes to
$$ 
HF(L_{1}, L_{12}, \ldots ,L_{(r-1)r}, L_r)
= HF(L_1 \times L_{23} \times \ldots , L_{12} \times L_{34} \times \ldots ) .$$
In particular we reproduce the definition of Floer cohomology for a pair of Lagrangians
$L_0,L_1\subset M$, viewed as cyclic correspondence 
$\{pt.\}\overset{L_0}{\longrightarrow} M \overset{L_1}{\longrightarrow} \{pt.\}$.
We moreover define a Floer cohomology for any Lagrangian $L\subset M^-\times M$, viewed as cyclic correspondence $M\overset{L}{\longrightarrow} M$, in particular for graphs 
$L={\rm graph}(\phi)$ of symplectomorphisms $\phi:M\to M$. By definition, this invariant is
$HF(L):= HF(L,\Delta_M)$, which reproduces the Floer cohomology 
$HF({\rm graph}(\phi))= HF({\rm graph}(\phi),\Delta_M)=HF(\phi)$ of a symplectomorphism.

\section{Quilted Floer cohomology}
\label{quilted FH}

The purpose of this section is to reformulate the definition of Floer cohomology for
generalized Lagrangian correspondences in terms of quilted surfaces (consisting of strips).
As in Section \ref{sequences} consider a cyclic generalized Lagrangian
correspondence $\ul{L}$, that is, a sequence of symplectic manifolds
$M_0,M_1,\ldots, M_{r},M_{r+1}$ with $M_0=M_{r+1}$
for $r \ge 0$, and a sequence of Lagrangian correspondences
$$L_{01} \subset M_0^- \times M_1, \ \ 
L_{12} \subset M_1^- \times M_2, \ \ \ldots, \ \
L_{r(r+1)} \subset M_{r}^- \times M_{r+1} .$$

\subsection{Unfolding of Floer cohomology in products} \label{unfold}

We defined the Floer cohomology $HF(\ul{L})$ as the standard Floer cohomology in the product manifold $\widetilde{M}=M_0^-\times M_1\times M_2^- \times \ldots$ of a pair of Lagrangians $L_{(0)}, L_{(1)}$ that is built from the cyclic sequence $\ul{L}$.
We will show how quilts arise naturally from ``unfolding" this construction and phrasing it in terms of tuples of holomorphic curves in the $M_j$.

Informally, $HF(\ul{L})$ can be viewed as the Morse homology on the path space
$$
\cP({L}_{(0)},{L}_{(1)}) = \bigl\{
{y}:[0,1]\to \widetilde{M} \,\big|\, y(0)\in{L}_{(0)}, y(1)\in {L}_{(1)} \bigr\}
$$
of the (potentially multivalued) symplectic action functional 
$$
\cA_H(y) = \int_{[0,1]\times[0,1]} v^*\omega_{\widetilde M} + \int_0^1 H(t,y(t)) {\rm d}t .
$$
Here $v:[0,1]\to\cP({L}_{(0)},{L}_{(1)})$ is a smooth homotopy from a fixed
$v(0)=y_0\in\cP({L}_{(0)},{L}_{(1)})$ (in a given connected component) to $v(1)=y$,
which can also be viewed as map $v:[0,1]\times[0,1]\to\widetilde{M}$ satisfying Lagrangian boundary conditions on $\{0\}\times[0,1]$ and $\{1\}\times[0,1]$.

Suppose for now that $r$ is odd, then the path space can be identified
with the set of tuples of paths in the manifolds $M_j$, 
connected via $L_{j(j+1)}$-matching conditions at the ends,
\begin{equation*} 
\cP(\ul{L}) = \bigl\{ \ul{x}=\bigl(x_j: [0,1] \to M_j\bigr)_{j=0,\ldots,r} \, \big|
(x_{j}(1),x_{j+1}(0)) \in L_{j(j+1)} \bigr\} .
\end{equation*}
Here and throughout we will use the index $j\in\{0,\ldots,r\}$ modulo $r+1$, so e.g.\ $x_{r+1}:=x_0$ and the matching condition for $j=r+1$ is $(x_{r}(1),x_{0}(0)) \in L_{r(r+1)}$.
We make the identification with $\cP({L}_{(0)},{L}_{(1)})$ by ${y}(t)=\bigl(x_0(1-t),x_1(t),x_2(1-t),\ldots,x_r(t)\bigr)$, then the unperturbed ($H=0$) symplectic action functional on $\cP(\ul{L})$ becomes
$$
\cA_0(\ul{x}) = \sum_{j=0}^r \int_{[0,1]\times[0,1]} v_j^*\omega_{M_j} .
$$
Here $v_j:[0,1]\times[0,1]\to M_j$ interpolate between fixed paths $v_j(0,\cdot)$ and
$v_j(1,\cdot)=x_j$, and satisfy what we will call ``seam conditions" 
$(v_{j}(s,1),v_{j+1}(s,0)) \in L_{j(j+1)}$ for all $s\in[0,1]$.
Next, assume that the almost complex structure on $\widetilde{M}$ is of time-independent split form ${J}=(-J_0)\oplus J_1 \oplus (-J_2)\oplus\ldots \oplus J_r$, given by a tuple
$J_j\in\J(M_j,\omega_j)$ of almost complex structures on the factors of $\widetilde{M}$. 
This defines a metric on the path space, and the gradient flow lines, viewed as solutions of PDE's are
the $J$-holomorphic strips $w:\R\times[0,1]\to\widetilde{M}$ with boundary values in $L_{(0)}$ and $L_{(1)}$. They are in one-to-one correspondence with $(r+1)$-tuples of 
$J_j$-holomorphic maps $u_j : \R \times [0,1] \to M_j$ satisfying the seam conditions
$$
(u_{j}({s},1),u_{j+1}({s},0)) \in L_{j(j+1)},\qquad\text{for all}\; 
j = 0,\ldots,r ,\ s\in\R.
$$
Here we again use cyclic notation $u_{r+1}:=u_0$, and the correspondence is given by
${w}(s,t)=\bigl(u_0(s,1-t),u_1(s,t),u_2(s,1-t),\ldots,u_r(s,t)\bigr)$.

For $r$ even there is a slight modification of the previous correspondence. The product manifold $\widetilde{M}$ has two factors $M_0$ and $M_{r+1}=M_0$ matched up via the diagonal. So the path space can be identified with the generalized path space $\cP(\ul{L})$ as above with the exception that the path $x_0:[0,2]\to M_0$ in $M_0=M_{r+1}$ is parametrized by an interval of length $2$ and satisfies the matching condition $(x_{0}(2),x_{1}(0) \in L_{01}$ at its end.
Similarly, a $J$-holomorphic strip $w:\R\times[0,1]\to\widetilde{M}$ corresponds via
${w}(s,t)=\bigl(u_0(s,2-t), u_1(s,t),u_2(s,1-t) \ldots, u_r(s,1-t),u_0(s,t)\bigr)$
to a tuple of $J_j$-holomorphic strips as above, with the exception that the strip 
$u_0:\R\times[0,2]\to M_0$ has width $2$.
This tuple $(u_j)_{j=0,\ldots,r}$ is the first instance of a nontrivial pseudoholomorphic quilt -- containing strips of different widths.

When $r$ is even, the Floer trajectories of the pair $L_{(0)},L_{(1)}$ in fact cannot be identified with an $(r+1)$-tuple of pseudoholomorphic maps, all defined on strips of width $1$, with seam conditions in $L_{j(j+1)}$.
Conformal rescaling $\tilde u_0(s,t):= u_0(2s,2t)$ would result in a ``time-shifted" matching condition $(\tilde u_{0}({s},1),u_{1}({2s},0)) \in L_{01}$ unless $u_1$ is rescaled, too, which would result in $\tilde u_0$ having width $1$ but all other strips having width $\frac 12$.
In fact, only simultaneous rescaling of all components in these pseudoholomorphic quilts preserves holomorphicity and seam conditions (unless the Lagrangian correspondences are of split type, e.g.\ $L_{01}=L_0\times L_1$ for Lagrangians $L_j\subset M_j$). It cannot change the relative widths of strips.

By a reparametrization of the path in $M_0$, one could identify $\cP({L}_{(0)},{L}_{(1)})$ and the action functional with the generalized path space $\cP(\ul{L})$ and a corresponding action functional, where all paths are parametrized by $[0,1]$. 
However, the reparametrized $\hat u_0(s,t):= u_0(s,2t)$ now satisfies $\partial_s \hat u_0 + \frac 12 J_0 \partial_t \hat u_0 =0$ with a no longer complex structure $\frac 12 J_0$ that squares to $-\frac 14$.
This is due to the fact that the pullback of the metric on $\cP({L}_{(0)},{L}_{(1)})$ to $\cP(\ul{L})$ is the $L^2$-metric on each factor with respect to $\omega_j(\cdot, J_j \cdot)$ for $j=1,\ldots,r$ but $\frac 12 \omega_0(\cdot, J_0 \cdot)$ on $M_0$. 
We could drop the factor $\frac 12$ in the metric on $M_0$ to obtain $J_j$-holomorphic strips of width $1$ in each factor as trajectories, however these would be the ``gradient flow lines" with respect to a different metric. In general, it is not known how Floer homology behaves under a change of metric. However, we will show that it is independent of the choice of weights 
$\delta_j^{-1}\omega_j(\cdot, J_j \cdot)$ in the $L^2$-metric on $\cP(\ul{L})$.
This setup is equivalent to defining the generalized path space with varying widths
$x_j:[0,\delta_j]\to M_j$ but fixing the standard $L^2$-metric induced by $\omega_j$ and $J_j$ on each factor.

\subsection{Construction of quilted Floer cohomology}
\label{cqf}

In the quilted setup for $HF(\ul{L})$ we fix widths $\ul{\delta}=(\delta_j>0)_{j=0,\ldots,r}$
and consider the generalized path space
\begin{equation*} \label{cP}
\cP(\ul{L}) := \bigl\{ \ul{x}=\bigl(x_j: [0,\delta_j] \to M_j\bigr)_{j=0,\ldots,r} \, \big|
(x_{j}(\delta_j),x_{j+1}(0)) \in L_{j(j+1)} \bigr\} .
\end{equation*}
We define a perturbed symplectic action functional on $\cP(\ul{L})$ by
picking a homotopy $\ul{v}=(v_j)_{j=0,\ldots,r}:[0,1]\to \cP(\ul{L})$ from a 
fixed $\ul{v}(0)$ to $\ul{v}(1)=\ul{x}$ and setting
$$
\cA_{\ul{H}}(\ul{x}) = \sum_{j=0}^r \biggl( \int_{[0,1]\times[0,\delta_j]} v_j^*\omega_{M_j} 
+ \int_0^{\delta_j} H_j(t,x_j(t)) {\rm d} t  \biggr),
$$
using a tuple of Hamiltonian functions
$$ 
\ul{H}=\bigl( H_j \in C^\infty([0,\delta_j] \times M_j) \bigr)_{j = 0,\ldots,r} .
$$
By folding and rescaling as in the previous section, this is equivalent to the path space 
$\cP(L_{(0)},L_{(1)})$ with symplectic action functional perturbed by a Hamiltonian of split type, e.g. $H=\sum_{j=0}^{r} (-1)^{j+1} \delta_j\ti{H}_{j}$ for $r$ odd,
% resp.\  $H= \frac 12 \ti{H}_0 + \sum_{j=1}^{r} \ti{H}_{j} + \frac 12 \ti{H}_{r+1}$ for $r$ even.
where $\ti{H_j}(t,x)= H_j(\delta_j t,x)$ for $j$ odd and
$\ti{H_j}(t,x)= H_j(\delta_j (1-t),x)$ for $j$ even.
Here the critical points correspond to the perturbed intersection points
$\phi_1^H(L_{(0)})\cap{L}_{(1)}$, where $\phi^H_1$ is the time-one flow of $H$.
In the quilted setup, the critical points of $\cA_{\ul{H}}$ are tuples of Hamiltonian chords,
\begin{equation*} \label{cI}
\cI(\ul{L}) := \left\{ \ul{x}=\bigl(x_j: [0,\delta_j] \to M_j\bigr)_{j=0,\ldots,r} \, \left|
\begin{aligned}
\dot x_j(t) = X_{H_j}(x_j(t)), \\
(x_{j}(\delta_j),x_{j+1}(0)) \in L_{j(j+1)} 
\end{aligned} \right.\right\} .
\end{equation*}
 $\cI(\ul{L})$ is canonically identified with $\times_{\phi_{\delta_0}^{H_0}}\bigl( L_{01}
\times_{\phi_{\delta_1}^{H_1}} L_{12} \ldots \times_{\phi_{\delta_r}^{H_r}}
L_{r(r+1)}\bigr)$,
the set of points 
$$ \bigl\{ (m_0,\ldots,m_r) \in M_0 \times \ldots \times M_r
\,\big|\, (\phi_{\delta_j}^{H_j}(m_{j}),m_{j+1}) \in L_{j(j+1)} \bigr\} , $$
where $\phi_{\delta_j}^{H_{j}}$ is the time $\delta_j$ flow of the Hamiltonian $H_{j}$.
In this setting we can check that Hamiltonians of split type suffice
to achieve transversality for the intersection points.

\begin{proposition}  \label{pert split}
There is a dense open subset 
$\Ham(\ul{L})\subset\oplus_{j=0}^r C^\infty([0,\delta_j] \times M_j)$
such that for every $(H_0,\ldots,H_r)\in\Ham(\ul{L})$
the set $ \times_{\phi_{\delta_0}^{H_0}}\bigl( L_{01} \times_{\phi_{\delta_1}^{H_1}} L_{12} \ldots \times_{\phi_{\delta_r}^{H_r}} L_{r(r+1)}\bigr)$ 
is smooth and finite, that is, the defining equations are transversal.
\end{proposition}  

\begin{proof}
The defining equations for 
$ \times_{\phi_{\delta_0}^{H_0}}\bigl( L_{01} \times_{\phi_{\delta_1}^{H_1}} L_{12} \ldots \times_{\phi_{\delta_r}^{H_r}} L_{r(r+1)}\bigr)$ 
are
\begin{equation} \label{kdefine}
m_j' = \phi_{\delta_j}^{H_j}(m_{j})
\qquad \text{for all}\;  j=0,\ldots,r
\end{equation}
for
$(m_0',m_1,m_1',m_2,\ldots,m_r',m_0) \in L_{01}\times L_{12}\ldots \times L_{r(r+1)}$.
Consider the universal moduli $\U$ space of data
$(H_0,\ldots,H_r,m_0',m_1,\ldots,m_r',m_0)$ satisfying \eqref{kdefine}, where now
each $H_j$ has class $C^\ell$ for some $\ell\geq1$.  
It is cut out by the diagonal values of the $C^\ell$-map
\begin{align*}
 L_{01}\times L_{12}\ldots \times L_{r(r+1)} \times 
 \bigoplus_{j=0,\ldots,r} C^\ell([0,1] \times M_{j})
&\longrightarrow \bigoplus_{j=0,\ldots,r} M_j\times M_j, \\ 
(m_j,m_j',H_j)_{j=0,\ldots,r}& \longmapsto ( \phi_{\delta_j}^{H_j}(m_{j}),m_j')_{j=0,\ldots,r} .
\end{align*}
The linearized equations for $\U$ are
\begin{equation} \label{surj} 
v_j' - D \phi_{\delta_j}^{H_j} (h_{j},v_{j}) = 0 \in TM_j
\qquad \text{for all}\; j = 0,\ldots, r .
\end{equation}
for $v_j,v_j' \in T_{m_j} M_j$  and $h_j \in C^\ell([0,1] \times M_j)$.
The map
$$ C^\ell([0,1] \times M_{j}) \to T_{\phi_{\delta_j}^{H_j}(m_{j})}
M_{j}, \quad h_{j} \mapsto D \phi_{\delta_j}^{H_j} (h_{j},0) $$
is surjective, which shows that the product of the operators on the
left-hand side of \eqref{surj} is also surjective.
So by the implicit function theorem $\U$ is a $C^\ell$ Banach manifold, 
and we consider its projection to $\oplus_{k=0}^r C^\ell([0,\delta_k] \times M_k)$.  
This is a Fredholm map of class $C^\ell$ and index $0$.
Hence, by the Sard-Smale theorem, the set of regular values (which coincides with the set of functions $H=(H_0,\ldots,H_r)$ such that the perturbed intersection is
transversal) is dense in $\oplus_{k=0}^r C^\ell([0,\delta_k] \times M_k)$.
Moreover, the set of regular values is open for each $\ell\geq 1$.
Indeed, by the compactness of $L_{01}\times L_{12}\ldots \times L_{r(r+1)}$, a $C^1$-small change in $H$ leads to a small change in perturbed intersection points, with small change in the linearized operators.

Now, by approximation of $C^\infty$-functions with $C^\ell$-functions, the set of regular values in $\oplus_{k=0}^r C^\infty([0,\delta_k] \times M_k)$ is dense in the $C^\ell$-topology for all $\ell\geq 1$, and hence dense in the $C^\infty$-topology.
Finally, the set of regular smooth $H$ is open in the $C^\infty$-topology as a special case of the $C^1$-openness.
\end{proof} 

In the proof of Theorem~\ref{main2} we will use the following special choice of corresponding regular Hamiltonian perturbations, which will provide a canonical identification of intersection points, as in Remark~\ref{rmk:canonical ident} for the unperturbed case.

\begin{remark} \label{rmk:canonical ident H}
Consider two cyclic generalized Lagrangian correspondences 
\begin{align*}
\ul{L}&=(L_{01},\ldots,L_{(j-1)j},L_{j(j+1)}, \ldots, L_{r(r+1)}), \\
\ul{L}'&=(L_{01},\ldots,L_{(j-1)j}\circ L_{j(j+1)}, \ldots, L_{r(r+1)})
\end{align*}
such that the composition $L_{(j-1)j}\circ L_{j(j+1)}$ is embedded
in the sense of Definition~\ref{embedded}.
Pick $\ul{H}'=(\ldots,H_{j-1},H_{j+1},\ldots )\in\Ham(\ul{L}')$
such that $ \times_{\phi_{\delta_0}^{H_0}}\bigl( L_{01} \ldots \times_{\phi_{\delta_{j-1}}^{H_{j-1}}} L_{(j-1)j}\circ L_{j(j+1)} \times_{\phi_{\delta_{j+1}}^{H_{j+1}}} \ldots L_{r(r+1)}\bigr)$ 
is transverse.
Then we have $\ul{H}:=(\ldots,H_{j-1},H_j\equiv 0,H_{j+1},\ldots )\in\Ham(\ul{L})$, that is
$ \times_{\phi_{\delta_0}^{H_0}}\bigl( L_{01} \ldots \times_{\phi_{\delta_{j-1}}^{H_{j-1}}} L_{(j-1)j} \times_{\rm Id} L_{j(j+1)} \times_{\phi_{\delta_{j+1}}^{H_{j+1}}} \ldots L_{r(r+1)}\bigr)$ 
is transverse, since by assumption $L_{(j-1)j}\times L_{(j+1)j}$ is transverse
to the diagonal $M_{j-1}\times \Delta_{M_j}\times M_{j+1}$.
Moreover, the generalized intersection points
\begin{align*}
\cI(\ul{L},\ul{H})&=
\bigl\{(\ldots,m_{j-1},m_j,m_{j+1},\ldots) \in \ldots \times M_{j-1}\times M_j \times M_{j+1} \ldots  
\,\big|\,  \\
&\qquad\qquad\qquad\qquad
\ldots,(\phi^{H_{j-1}}_{\delta_{j-1}}(m_{j-1}),m_j)\in L_{(j-1)j}, (m_{j},m_{j+1})\in L_{j(j+1)}, \ldots  \bigr\}\\
&=
\bigl\{ (\ldots,m_{j-1},m_{j+1},\ldots) \in \ldots \times M_{j-1} \times M_{j+1} \ldots  \,\big|\, \\
&\qquad\qquad\qquad\qquad
 \ldots , (\phi^{H_{j-1}}_{\delta_{j-1}}(m_{j-1}),m_{j+1})\in L_{(j-1)j}\circ L_{j(j+1)}, \ldots  \bigr\}
=\cI(\ul{L}',\ul{H}')
\end{align*}
are canonically identified, since the intermediate point 
$m_{j}\in M_j$ is uniquely determined by the pair $(\phi^{H_{j-1}}_{\delta_{j-1}}(m_{j-1}),m_{j+1})\in L_{(j-1)j}\circ L_{j(j+1)}$.
\end{remark}

With this  split Hamiltonian perturbation we have a canonical bijection
of critical points $\phi^H_1(L_{(0)})\cap L_{(1)}\cong\cI(\ul{L})$, and hence
the (graded) Floer chain group $CF({L}_{(0)},{L}_{(1)})$ is identified with
$$ 
CF(\ul{L}) := \bigoplus_{d\in\Z_N} CF^d(\ul{L}), \qquad 
CF^d(\ul{L}) := \bigoplus_{\ul{x} \in \cI(\ul{L}), |\ul{x}|=d} \Z \bra{\ul{x}} .$$
The grading is defined as in Section \ref{graded cor},
$$  \cI(\ul{L})\cong \phi_1^H({L}_{(0)})\cap{L}_{(1)} \to \Z_N ,
\qquad \ul{x}\cong {y}\mapsto |{y}|=|\ul{x}| . $$ 
Next, fix a tuple of almost complex structures
$$ 
\ul{J}= (J_j)_{j=0,\ldots,r} \in
\oplus_{j=0}^r C^\infty([0,\delta_j],\J(M_j,\omega_j)) =:  \J_t(\ul{L})
$$
and equip $\cP(\ul{L})$ with the $L^2$-metric induced by the $t$-dependent metric $\omega_j(\cdot,J_j\cdot)$ on each factor $M_j$. Then the Floer trajectories (obtained by reformulating the gradient flow as PDE) are $(r+1)$-tuples of maps
$u_j : \R \times [0,\delta_j] \to M_j$ that are $(J_j,H_j)$-holomorphic,
\beq \label{Jjhol}
\overline\partial_{J_j,H_j} u_j = \partial_s u_j + J_j \bigl(\partial_t u_j - X_{H_j}(u_j)  \bigr) = 0 
\qquad\forall j = 0,\ldots r ,
\eeq
and satisfy the seam conditions
\beq\label{ubc}
(u_{j}({s},\delta_j),u_{j+1}({s},0)) \in L_{j(j+1)}\qquad\forall j = 0,\ldots r ,\ s\in\R.
\eeq
For a Floer trajectory to be counted towards the differential between critical points $\ul{x}^\pm\in\cI(\ul{L})$ we moreover require finite energy and limits 
\beq \label{ulim}
E(\ul{u}) := \sum_{j=0}^r \int_{\R \times [0,\delta_j]} u_j^*\omega_j + {\rm d}(H_j(u_j) {\rm d}t) <\infty,
\qquad\lim_{s\to\pm\infty} u_j(s,\cdot) = x^\pm_j \quad\forall j = 0,\ldots,r  .
\eeq

\begin{figure}[ht]
\begin{picture}(0,0)%
\includegraphics{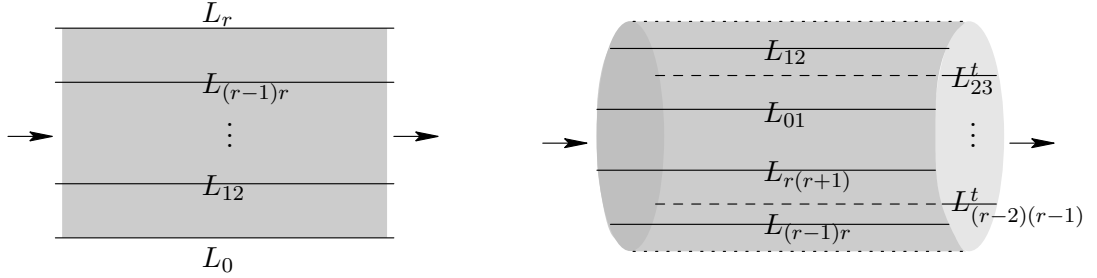}%
\end{picture}%
\setlength{\unitlength}{3729sp}%
\begingroup\makeatletter\ifx\SetFigFont\undefined%
\gdef\SetFigFont#1#2#3#4#5{%
  \reset@font\fontsize{#1}{#2pt}%
  \fontfamily{#3}\fontseries{#4}\fontshape{#5}%
  \selectfont}%
\fi\endgroup%
\begin{picture}(6310,1890)(3139,-1231)
\put(4441,-116){\makebox(0,0)[lb]{$L_{(r-1) r}$}}
\put(4441,404){\makebox(0,0)[lb]{$L_{r}$}}
\put(4441,-1231){\makebox(0,0)[lb]{$L_{0}$}}
\put(4601,-401){\makebox(0,0)[lb]{$\vdots$}}
\put(4441,-751){\makebox(0,0)[lb]{$L_{12}$}}
\put(8176,150){\makebox(0,0)[lb]{$L_{12}$}}
\put(8176,-266){\makebox(0,0)[lb]{$L_{01}$}}
\put(8176,-701){\makebox(0,0)[lb]{$L_{r(r+1)}$}}
\put(8176,-1051){\makebox(0,0)[lb]{$L_{(r-1)r}$}}
\put(9541,-401){\makebox(0,0)[lb]{$\vdots$}}
\put(9421,-51){\makebox(0,0)[lb]{$L_{23}^t$}}
\put(9421,-931){\makebox(0,0)[lb]{$L_{(r-2)(r-1)}^t$}}
\end{picture}%
\caption{Quilted Floer trajectories for $M_0=\{\pt\}$ and in general}
\label{qFH trajectories}
\end{figure}

As in standard Floer theory, the moduli spaces of ``quilted holomorphic strips"
$$
\M(\ul{x}^-,\ul{x}^+) :=
\bigl\{ \ul{u}=\bigl( u_j : \R \times [0,\delta_j] \to M_j \bigr)_{j=0,\ldots,r} \,\big|\, 
\eqref{Jjhol} , \eqref{ubc}, \eqref{ulim} \bigr\} / \R
$$
arise from quotienting out by simultaneous $\R$-shift in all components $u_j$.
(Separate shifts will not preserve the seam condition unless the correspondences are of split type.)
We will see that they have the same Fredholm, exponential decay, and compactness properties as usual for Floer trajectories. For that purpose we restrict ourselves to the monotone case.

\begin{remark} \label{L monotone}
The ``monotonicity for Floer theory" assumption for the pair $({L}_{(0)},{L}_{(1)})$ in Definition \ref{tuple monotone}
can be phrased directly for $\ul{L}$ in the language of \cite{quilts}:
``$\ul{L}$ is a monotone boundary condition for the quilted cylinder". 
That is, the action-index relation
$$ 
2\sum_{j=0}^r \int u_j^*\omega_j = \tau \cdot 
I\bigl((u_j^*TM_j)_{j=0,\ldots,r}, ( s_{j(j+1)}^*TL_{j(j+1)})_{j=0,\ldots,r} \bigr) $$
holds for each tuple of maps $u_j:S^1\times[0,\delta_j]\to M_j$ that satisfies the seam conditions
$s_{j(j+1)}(s):=(u_{j}({s},\delta_j),u_{j+1}({s},0)) \in L_{j(j+1)}$ for $j=0,\ldots, r$.
Here the topological index $I$ is defined by choosing a trivialization for each $u_j^*TM_j$ and then summing over the Maslov indices of the loops $s_{j(j+1)}^*TL_{j(j+1)}$ of Lagrangian subspaces with respect to these trivializations.

Note that the monotonicity condition for $\ul{L}$ is independent of the width $\delta_j$ of the annuli that parametrize the maps $u_j$.
Moreover, it implies monotonicity for the sequence $\ul{L}'=(L_{01},\ldots, L_{(j-1)j} \circ L_{j(j+1)}, \ldots,L_{r(r+1)})$ obtained from an embedded composition
$L_{(j-1)j} \circ L_{j(j+1)}$. 
To see the latter note that any seam condition $s_{(j-1)(j+1)}:S^1\to L_{(j-1)j} \circ L_{j(j+1)}$ induces a smooth map $u_j:S^1\times [0,1]\to M_j$ that is constant in $[0,1]$, fits the seam conditions for $\ul{L}$, but contributes zero to both energy and Maslov index. 
Hence the action index relation for $\ul{L}$ implies the same relation for $\ul{L}'$.

Indeed, to identify the Maslov indices pick the trivialization of $u_j^*TM_j$ constant across $[0,1]$.
Then the Maslov index for $\ul{L}$ has a contribution 
$I(\Lambda_{(j-1)j})+I(\Lambda_{j(j+1)})= I(\Lambda_{(j-1)j}\times\Lambda_{j(j+1)})$ 
from the trivializations $\Lambda_{i(i+1)}:S^1\to\C^{n_i+n_{i+1}}$ of $s_{i(i+1)}^*L_{i(i+1)}$, 
and the contribution to the Maslov index for $\ul{L}'$ is the index
$I(\Lambda_{(j-1)(j+1)}:=\Lambda_{(j-1)j}\circ\Lambda_{j(j+1)})$ of the geometric composition of the trivializations.
By Lemma~\ref{the one} we can homotope $\Lambda_{(j-1)j}\times\Lambda_{j(j+1)}$ to $(\Lambda_{(j-1)(j+1)} \times \Lambda_{11})^T$, where $\Lambda_{11}\subset\C^{n_1}\times \C^{n_1}$ is a fixed complement of the diagonal and $(\cdot)^T$ is the exchange of factors as in Section~\ref{graded cor}.
Then the crossing form definition of the Maslov index proves
$$
I(\Lambda_{(j-1)j}\times\Lambda_{j(j+1)})
= I\bigl((\Lambda_{(j-1)(j+1)} \times \Lambda_{11})^T\bigr)
= I(\Lambda_{(j-1)(j+1)}) + I(\Lambda_{11}) 
= I(\Lambda_{(j-1)(j+1)}).
$$
\end{remark}

\begin{theorem} \label{quilt traj} 
Suppose that the symplectic manifolds $M_j$ satisfy (M1-2) with the same value of the
monotonicity constant $\tau$, the Lagrangian correspondences $L_{j(j+1)}$ satisfy (L1-2), 
and $\ul{L}$ satisfies the monotonicity assumption of Remark~\ref{L monotone}.

For any  choice of widths $\ul{\delta}$ and regular Hamiltonian perturbations 
$\ul{H} \subset \Ham(\ul{L})$ there exists a dense subset
$\J_t^{\reg}(\ul{L};\ul{H}) \subset \J_t(\ul{L})$  such that the following holds for all $\ul{x}_\pm\in\cI(\ul{L})$.
\ben 
\item 
$\M(\ul{x}_-,\ul{x}_+)$ is a smooth manifold whose dimension near a nonconstant solution $\ul{u}$
is given by the formal dimension, equal to $\Ind(D_{\ul{u}})-1$.
Here $D_{\ul{u}}$ is the linearized operator at $\ul{u}$ of \eqref{Jjhol} on the space of sections satisfying the linearized boundary- and seam conditions of \eqref{ubc}.
We denote the component of nonconstant solutions of formal dimension $j$ by  $\M(\ul{x}_-,\ul{x}_+)_j := \{ \Ind(D_{\ul{u}})-1 = j \}$.
\item 
The component $\M(\ul{x}_-,\ul{x}_+)_0 \subset \M(\ul{x}_-,\ul{x}_+)$ 
of formal dimension zero is finite.
\item 
Suppose that each $L_{j(j+1)}$ has minimal Maslov number $N_{L_{j(j+1)}}\geq 3$.  
Then the one-dimensional component $\M(\ul{x}_-,\ul{x}_+)_1 \subset \M(\ul{x}_-,\ul{x}_+)$ has
a compactification as one-dimensional manifold with boundary
\begin{equation} \label{qglue}
\partial \overline{\M(\ul{x}_-,\ul{x}_+)_{1}} \cong \bigcup_{\ul{x} \in \cI(\ul{L})}
\M(\ul{x}_-,\ul{x})_{0} \times \M(\ul{x},\ul{x}_+)_{0} .
\end{equation}
\item If $\ul{L}$ is relatively spin (as defined in Definition~\ref{gen spin c}), 
then there exists a coherent set of orientations on
$\M(\ul{x}_-,\ul{x}_+)_0,\M(\ul{x}_-,\ul{x}_+)_1$ for all $\ul{x}_\pm \in \cI(\ul{L})$, that is,
orientations compatible with \eqref{qglue}.
\een
\end{theorem}

\begin{proof}
Suppose for simplicity that $r$ is odd. (For even $r$ we can insert a diagonal into the sequence $\ul{L}$, then the quilted holomorphic strips of widths $\ul{\delta}$ can be identified with quilted holomorphic strips for the new sequence with widths $(\frac{\delta_0}2, \delta_1,\ldots,\delta_r,\frac{\delta_0}2)$.) Then the quilted moduli space $\M(\ul{x}_-,\ul{x}_+)$ is canonically identified with the moduli space of $(J_{\ul{\delta}},H)$-holomorphic maps ${w} : \R \times [0,1] \to \widetilde{M}$ 
with boundary conditions ${w}(\R,0)\subset{L}_{(0)}$, ${w}(\R,1)\subset{L}_{(1)}$, 
finite energy $E_H(w)<\infty$, and limits $\lim_{s\to\pm\infty}w(s,\cdot) = \ul{x}_\pm \in \cI({L}_{(0)},{L}_{(1)})$.

The correspondence is by 
$$
{w}(s,t)=\bigl(u_0(s,\delta_0(1-t)),u_1(s,\delta_1 t),u_2(s,\delta_2(1-t)),\ldots,u_r(s,\delta_r t)\bigr),
$$
where $H=\sum_{j=0}^{r} (-1)^{j+1} \delta_j\ti{H}_{j}$ as above and
$$
J_{\ul{\delta}}:=\bigl(-\delta_0^{-1}J_0(\delta_0(1-t)),\delta_1^{-1}J_1(\delta_1 t),\ldots,\delta_r^{-1}J_r(\delta_r t)\bigr)
$$
satisfies all properties of a $t$-dependent $\omega_{\widetilde{M}}$-compatible almost complex structure except that it squares to the negative definite diagonal matrix 
$J_{\ul{\delta}}^2= - (\delta_0^{-2}{\rm Id}_{TM_0}\oplus \ldots \oplus \delta_r^{-2}{\rm Id}_{TM_r})$ instead of $-{\rm Id}$. Let us call it a ``scaled almost complex structure".
Most analytic properties of pseudo-holomorphic curves carry over directly to $J_{\ul{\delta}}$-holomorphic curves.
Indeed, $\overline{\partial}_{J_{\ul{\delta}},H}$ still presents a partial differential operator of the form $\partial_s + \D$, where the linearizations of $\D$ are self-adjoint operators 
on $L^2([0,1],T\widetilde{M})$ with boundary conditions in $TL_{(0)}, TL_{(1)}$. 
Moreover, in local coordinates 
$(\partial_s - \D) (\partial_s + \D) = \partial_s^2 - J_{\ul{\delta}}^2 \partial_t^2 + \text{(lower order terms)}$ is an elliptic operator (i.e.\ has an elliptic symbol), and in the splitting $T\widetilde{M}= TL_{(0)} \oplus J_{\ul{\delta}}^{-1}TL_{(0)}$ the Lagrangian boundary conditions induce Neumann conditions (from\ $\partial_t w - X_H(w) |_{t=0}= - J_{\ul{\delta}}^{-1}\partial_s w |_{t=0} \in  J_{\ul{\delta}}^{-1}  TL_{(0)}$) resp.\ Dirichlet conditions (from $w|_{t=0}\in L_{(0)}$) on the two factors.

With these remarks in mind, we can follow the standard construction of Floer theory (which is currently probably best outlined in \cite{salamon:floer lectures} for the case of holomorphic cylinders, fully executed in \cite{Donaldson book} for a gauge theoretic setting, and hopefully soon available in \cite{oh book} for holomorphic strips).
The moduli space (before quotienting by the $\R$-action) is described as the zero set of the scaled Cauchy-Riemann operator $\overline{\partial}_{J_{\ul{\delta}},H}$, which is
a Fredholm section of a Banach bundle over the usual Banach manifold of maps
$w : \R \times [0,1] \to \widetilde{M}$ satisfying boundary conditions in $L_{(0)}, L_{(1)}$ and converging uniformly to $\ul{x}_\pm$ for $s\to\pm\infty$.
(The Banach manifold is modeled as usual for some $p>2$ on the Sobolev space of maps in $W^{1,p}(\R\times[0,1],\C^N)$ which take boundary values in $\R^N$ resp.\ $i\R^N$.)

The Fredholm property of the linearized operators follows as in e.g.\ \cite{Floer:unregularized} from the general theory of \cite{lock mc} 
(also see \cite{Donaldson book}) for operators of the form
$\partial_s + D_s$, where the operators $D_s$ converge to invertible operators as $s\to\pm\infty$. 
After a Hamiltonian transformation (moving the perturbation onto the Lagrangian $L_{(0)}$ and replacing $J_{\ul{\delta}}$ with $\phi^H_*J_{\ul{\delta}}$, which retains the same properties) these operators take the form $(\phi^H_*J_{\ul{\delta}})\partial_t$ on $L^p([0,1],T\widetilde{M})$ with domain given by $W^{1,p}$-paths satisfying boundary conditions in $T\phi^H_1(L_{(0)}), TL_{(1)}$. 
Note that $\phi^H_*J_{\ul{\delta}}$ is an invertible operator on $L^p([0,1],T\widetilde{M})$,
and invertibility of $\partial_t$ follows as usual from the fact that the boundary conditions are transverse on the ends $s\to\pm\infty$.
Similar considerations (for $p=2$ showing that each $D_s$ is self-adjoint and for sufficiently large $|s|$ has a spectral gap -- eigenvalues bounded uniformly away from $0$) are the crucial ingredient in proving that solutions of the nonlinear equation on long strips of small energy converge exponentially to intersection points. Details can be found in \cite{Floer:unregularized}, \cite{salamon:floer lectures}, or  \cite[Lemma 3.2.3]{isom}
(where we prove the presently irrelevant fact that the exponential decay rate is in fact uniform for certain families of widths $\ul{\delta}$). 

To calculate the index of the linearization of $\overline{\partial}_{J_{\ul{\delta}},H}$, 
one can deform $J_{\ul{\delta}}$ through the endomorphisms 
$\bigl(\exp(\tau \ln\delta_0) {\rm Id}_{TM_0},\ldots,\exp(\tau\ln\delta_r) {\rm Id}_{TM_r}\bigr) 
\circ J_{\ul{\delta}}$ to a true almost complex structure at ${\tau=1}$.
This provides a continuous family of Fredholm operators, along which the index is constant, and ending at a traditional Cauchy-Riemann operator whose index is given by a Maslov index. (see e.g.\ \cite{Floer:Morse index,se:bo}).
In particular, we have ${\rm ind}(D_{\ul{u}})=0$ for the constant solution in case $\ul{x}_-=\ul{x}_+$.
This identification of the index with a Maslov index together with the monotonicity assumption implies an energy-index relation as in Remark~\ref{rmk:monotone} for solutions (for another proof see e.g.\ \cite{oh:fl1}). Hence all solutions with a fixed index satisfy a uniform $L^2$-bound on the gradient. But before proceeding to compactness, let us assume that the section $\overline{\partial}_{J_{\ul{\delta}},H}$ is transverse to the zero section, i.e.\ its linearization at any zero is surjective. (This will be achieved by an appropriate choice of $\ul{J}$ -- see below.)
Then the implicit function theorem for Banach bundles (see e.g.\ \cite{ms:jh}) implies that the space of solutions ($\M(\ul{x}_-,\ul{x}_+)$ before quotienting by $\R$) is a smooth manifold, whose dimension near $\ul{u}$ is the Fredholm index ${\rm ind}(D_{\ul{u}})$.
Now, except for at a constant solution, the $\R$-action on these finite dimensional manifolds is smooth, proper, and free 
(whereas on the Banach manifold it is not even differentiable), inducing a smooth structure on the moduli spaces $\M(\ul{x}_-,\ul{x}_+)$ of dimension ${\rm ind}(D_{\ul{u}})-1$ at nonconstant solutions. This proves (a).

To prove (b) and (c) one next proves compactness properties of the moduli spaces of fixed index, as in \cite{Floer:unregularized, salamon:floer lectures}. 
By monotonicity, the fixed index provides a uniform $L^2$-bound on the gradients of solutions.
If, on the other hand, one had an $L^\infty$ gradient bound for a sequence of solutions, then elliptic estimates would imply that a subsequence converges with all derivatives on all compact subsets of $\R\times[0,1]$. (These estimates work exactly as in \cite[Appendix~B]{ms:jh}, using the splitting into Dirichlet and Neumann problem described above.)
For the moduli spaces of index $1$ and $2$ we can ensure $L^\infty$-bounds as follows:
We analyze any blow-up point of the gradient in the formulation as a tuple of maps $\ul{u}$. Then the usual rescaling analysis (see e.g.\ \cite{Floer:unregularized}) is local, in the interior of one component $u_j$ (leading to a $J_j$-holomorphic sphere in $M_j$) or near a seam, where we can consider $u_j(s,-t)\times u_{j+1}(s,t)$ as $(-J_j)\oplus J_{j+1}$-holomorphic curve with boundary condition in $L_{j(j+1)}$. The latter type of bubbling hence leads to a holomorphic disc in $M_j^-\times M_{j+1}$ with boundary on $L_{j(j+1)}$. Away from these blow-up points, the solution converges $C^\infty_{\rm loc}$ to a punctured solution, to which we can apply the usual removable singularity theorems
(e.g.\ \cite{ms:jh,oh:remsing}). In the limit we obtain a new solution (possibly with different end points) of nonnegative index and a number of holomorphic spheres and disks, each of which is nonconstant, so by monotonicity has positive index. By assumption (L3) they in fact must have Maslov index at least $3$ (i.e.\ Chern number $2$ for spheres). Since the Maslov index of all components adds up to $1$ or $2$, any bubbling is excluded.
(See \cite{oh:fl1} for the analogous argument in the standard theory.)
This discussion ensures that the moduli spaces of index $1$ and $2$ are compact up to the breaking of trajectories as in Morse theory. 
This is proven by combining the local elliptic estimates with the exponential decay on long strips, see \cite{Floer:unregularized, Donaldson book}.
Finally, part (c) requires a gluing theorem identifying the ends of the moduli space with broken trajectories. Again the proof in \cite{Floer:lag} (or \cite[Section 3.2]{sch:coh} for the closed case) carries over directly. The crucial ingredient is a uniformly bounded family of right inverses for the linearized operator as in \cite[Proposition~3.9]{salamon:floer lectures}, which is established by combining the already established exponential decay and Fredholm estimates.

Orientations can also be defined as in the standard Floer theory \cite{se:bo} since the linearized operator canonically deforms (as above) through Fredholm operators to a standard Cauchy-Riemann operator.

The only part of the standard construction of Floer theory on $\widetilde M$ that has to be adapted substantially is transversality: The scaled almost complex structures $J_{\ul{\delta}}$ obtained from tuples $\ul{J}$ of almost complex structures on each factor of 
$\widetilde M$ are highly non-generic as scaled almost complex structures on $\widetilde M$ (which generically do not respect the splitting into factors).
Nevertheless, we proceed as usual and define the set $\J_t^{\reg}(\ul{L};\ul{H})$ to be those tuples of complex structures $\ul{J}$, such that the corresponding scaled almost complex structure $J_{\ul{\delta}}\in\J_t({L}_{(0)},{L}_{(1)};{H})$ is regular in the sense that the linearized operator $D_u\overline{\partial}_{J_{\ul{\delta}},H}$ is surjective for every solution of $\overline{\partial}_{J_{\ul{\delta}},H}u=0$ with $L_{(0)},L_{(1)}$ boundary values.
(This is equivalent to the surjectivity at every solution $\ul{u}$ of the linearized operator $D_{\ul{u}}$ of \eqref{Jjhol} on the space of sections satisfying the linearized boundary- and seam conditions of \eqref{ubc}.)
In order to find a dense set of regular $\ul{J}$ we note that the unique continuation theorem \cite[Theorem 4.3]{fhs:tr} applies to the
interior of every single nonconstant strip $u_j:\R\times(0,\delta_j)\to M_j$.  
It implies that the set of regular points, $(s_0,t_0)\in\R\times(0,\delta_j)$ with
$\pd_s u_j(s_0,t_0)\neq 0$ and $u_j^{-1}(u_j(\R\cup\{\pm\infty\}) , t_0
) =\{(s_0,t_0)\}$, is open and dense. These points can be used to
prove surjectivity of the linearized operator for a universal moduli
space of solutions with respect to almost complex structures of class $C^k$. 
(The constant solutions are automatically transverse due to the previously
ensured transversality of the intersection points $\phi_1^H(L_{(0)})\pitchfork L_{(1)}$.)
The existence of a $C^k$-dense set of regular $\ul{J}$ then follows from the usual Sard-Smale argument as in \cite{ms:jh} for each $k\in\N$.
Finally, an intersection argument by Taubes (which here works exactly as in \cite{fhs:tr}) proves that the regular smooth almost complex structures are a 
comeagre\footnote{A subset of a topological space is {\it comeagre} if it is the intersection of countably many open dense subsets. Many authors in symplectic topology would use the term ``Baire second category'', which however in classical Baire theory \cite[Chapter 7.8]{royden} denotes more generally subsets that are not meagre, i.e.\ not the complement of a comeagre subset.
} subset of $\J_t(\ul{L})$ in the $C^\infty$-topology.
Since $\J_t(\ul{L})$ is a complete metric space\footnote{
It can be expressed as a closed subspace of a function space $C^\infty([0,1]\times\widetilde{M},\R^N)$. The  latter carries a metric $d(f,g) = \sum_{k=0}^\infty 2^{-k}\|\nabla^k(f-g)\|_\infty (1 + \|\nabla^k(f-g)\|_\infty )^{-1}$, which induces the $C^\infty$-topology.
}, 
and hence a Baire space, any such comeagre subset is dense as claimed, see e.g.\ \cite[Thm.27]{royden}.
\end{proof}

Now, assuming monotonicity and choosing regular $\ul{H}$ and $\ul{J}$ 
we can define the Floer homology $HF(\ul{L})$ just as in the standard case:
The Floer coboundary operator $ \partial^d : \ CF^d(\ul{L}) \to CF^{d+1}(\ul{L}) $ is 
defined by
$$
\partial^d\bra{\ul{x}_-} := \sum_{\ul{x}_+\in\cI(\ul{L})}
\Bigl( \sum_{\ul{u} \in\M(\ul{x}_-,\ul{x}_+)_0} \eps(\ul{u})\Bigr) \bra{\ul{x}_+} ,
$$
where the signs $ \eps: \M(\ul{x}_-,\ul{x}_+)_0 \to \{ \pm 1 \} $
are defined by comparing the given orientation to the canonical orientation of a point. 
It follows from Theorem \ref{quilt traj} (c) that $\partial^2 = 0$, 
and $\partial$ is a map of degree $1$ by index calculations as in the standard case.
This defines the {\em quilted Floer cohomology}
$$ HF(\ul{L}) := \bigoplus_{d\in\Z_N} HF^d(\ul{L}), \qquad
HF^d(\ul{L}) := {\ker(\partial^d)}/{\on{im}}(\partial^{d - 1}) $$
as $\Z_N$-graded group.
It is independent of the choice of $\ul{H}$ and $\ul{J}$ by a standard construction of continuation maps. The same construction also allows for a deformation of the widths $\ul{\delta}$, in the folded setup of the above proof, where the $\delta_j$ are merely scale factors in the endomorphism $J_{\ul{\delta}}$. For a more conceptual proof based on quilts interpolating between strips of different widths see Section~\ref{sec:invariance} below.

\begin{remark}  \label{empty} One can also allow the sequence $\ul{L}$ to have length
zero (that is, the empty sequence) as a generalized correspondence
from $M$ to $M$; this is the case $r = -1$ in the previous notation.
In this case we define $HF(\ul{L}) = HF(\Id_M)$, the cylindrical Floer
homology.  This would be the case without seams in Figure \ref{qFH trajectories}.
\end{remark} 

In the special case of a cyclic generalized Lagrangian correspondence $\ul{L}$ with
$M_0=\{\pt\}$ and some intermediate correspondence $L_{j(j+1)}$ of split form, we have the following
K\"unneth Theorem.

\begin{theorem} \label{kuenneth}
Suppose that $\ul{L}=(L_{01},\ldots, L_{r(r+1)})$ is a cyclic generalized Lagrangian correspondence as in Theorem~\ref{quilt traj} such that $M_0=\{\pt\}$ and $L_{j(j+1)}= L_j \times L_{j+1}$ for some $1\leq j < r$ and Lagrangian submanifolds $L_j\subset M_j$ and $L_{j+1}\subset M_{j+1}$.
Then the quilted Floer complex for $\ul{L}$ is canonically isomorphic to the tensor product 
of the two Floer complexes for the cyclic correspondences 
$(L_{01},\ldots,L_j )$ and $(L_{j+1},\ldots, L_{r(r+1)})$,
$$
 CF( L_{01},\ldots,L_j \times L_{j+1},\ldots, L_{r(r+1)})
\cong 
CF( L_{01},\ldots,L_j)
\otimes CF(L_{j+1},\ldots, L_{r(r+1)}).
$$
In particular, if either $HF( L_{01},\ldots,L_j)$
or $HF(L_{j+1},\ldots, L_{r(r+1)})$ is torsion-free, then there is a canonical isomorphism
$$
HF( L_{01},\ldots,L_j \times L_{j+1},\ldots, L_{r(r+1)})
\cong 
HF( L_{01},\ldots,L_j)
\otimes HF(L_{j+1},\ldots, L_{r(r+1)}).
$$
The latter isomorphism also holds if we work with coefficients in a field, e.g.\ with $\Z_2$ coefficients.
\end{theorem}
\begin{proof}
The generators and boundary operators of the Floer complex
$CF( L_{01},\ldots,L_j \times L_{j+1},\ldots, L_{r(r+1)})$
are trivially identified with the generators and differential of the Floer complex
$CF( L_{01},\ldots,L_j, L_{j+1},\ldots, L_{r(r+1)})$, where we view the two factors of
$L_j \times L_{j+1}$ as correspondences $L_j\subset M_j^-\times\{\pt\}$ 
and $L_{j+1}\subset \{\pt\}^-\times M_{j+1}$ to and from the point.
This is the first trivial case of Theorem~\ref{main2} below.
At the same time, the intersection points of $\ul{L}$ are trivially identified with the 
pairs of intersection points of $(L_{01},\ldots,L_j)$ and $(L_{j+1},\ldots, L_{r(r+1)})$,
and the Floer trajectories for $\ul{L}$ are either 
a pair of a Floer trajectory for $( L_{01},\ldots,L_j)$ 
and an intersection point of $(L_{j+1},\ldots, L_{r(r+1)})$
or a pair of an intersection point for $( L_{01},\ldots,L_j)$ 
and a Floer trajectory for $(L_{j+1},\ldots, L_{r(r+1)})$.
This canonically identifies the Floer complex $CF(\ul{L})$ with the tensor product of complexes $CF( L_{01},\ldots,L_j) \otimes CF(L_{j+1},\ldots, L_{r(r+1)})$.

Now the isomorphism of homologies follows from the general K\"unneth formula \cite[Lemma~5.3.1]{spanier} and the fact that the torsion product vanishes if it has one torsion-free factor \cite[Lemma 5.2.5]{spanier}.
Modules over fields are generally torsion-free.
\end{proof}

\subsection{Invariance of quilted Floer cohomology and relative quilt invariants}
\label{FH quilt inv} \label{sec:invariance}

The purpose of this section is to prove the independence of quilted Floer cohomology from the choice of perturbation data, in particular the choice of widths.

Consider a cyclic generalized Lagrangian
correspondence $\ul{L}=(L_{k(k+1)})_{k=0,\ldots,r}$ satisfying the monotonicity
conditions of Theorem~\ref{quilt traj}.
Fix a tuple of widths $\ul{\delta}=(\delta_k)_{k=0,\ldots,r}$. 
Then  Proposition \ref{pert split} and Theorem \ref{quilt traj} provide tuples of Hamiltonians
$\ul{H}=(H_k)_{k=0\ldots,r}$ and almost complex structures 
$\ul{J}=(J_k)_{k=0,\ldots,r}$ such that the Floer homology $HF(\ul{L})$ can be defined by counting quilted Floer trajectories $\ul{u}\in\M(\ul{x}^-,\ul{x}^+)$ between generalized intersection points 
$\ul{x}^\pm\in\cI(\ul{L})$.

The proof of independence of Floer cohomology from the choice of perturbations and particularly the widths goes somewhat beyond the proof for standard Floer theory. It is best formulated by using quilted surfaces that are not obtained by ``unfolding of strips".
With Proposition \ref{width indep} below in place we can in particular identify the two definitions of Floer cohomology $HF(\ul{L}) \cong HF({L}_{(0)},{L}_{(1)})$ for a cyclic sequence in Sections~\ref{sequences} and \ref{cqf}.
For that purpose one chooses special widths in the quilted setup of Section \ref{cqf}, namely those that correspond by the discussion in Section \ref{unfold} to the ``folded" Floer trajectories of $HF({L}_{(0)},{L}_{(1)})$.
The proof of the Proposition however uses the notation and construction of relative quilt invariants in \cite{quilts}. For readers familiar with this notation, the following Remark describes the quilted Floer trajectories as quilted surfaces. For readers unfamiliar with the quilt notation we will summarize the definition of relative quilt invariants in the special case used in the proof of Proposition~\ref{width indep} below.

\begin{remark} \label{quilt lingo}
In the language of quilted surfaces developed in \cite{quilts} the Floer trajectories correspond to the {\em holomorphic quilted cylinders}
$\ul{u}\in \M_{\ul{Z}}(\ul{x}^-,\ul{x}^+)$
with $\ul{K}=(H_k \d t)_{k=0,\ldots,r}$ and $\ul{J}=(J_k)_{k=0,\ldots,r}$.
Here the quilted surface is the quilted cylinder $\ul{Z}$ whose patches are strips  
$(S_k=\R\times[0,\delta_k])_{k=0,\ldots,r}$ of the given widths with the canonical complex structure and the obvious (up to a shift chosen as $\pm 1$) ends $\eps_{k,e_\pm}:\R^\pm\times[0,\delta_k]\to S_k, (s,t)\mapsto(s,\pm1+t)$.
The seams are $\sigma_k=\{(k,\R\times\{\delta_k\}),(k+1,\R\times\{0\})\}$ for $k=0,\ldots,r$ modulo $(r+1)$, with seam maps $\phi_{\sigma_k}: \partial S_k \supset(s,\delta_k)\mapsto (s,0)\subset\partial S_{k+1}$. 
This quilted surface is shown on the right in Figure \ref{qFH trajectories}.
There are no remaining boundary components except for in the special case of a noncyclic sequence with $M_0=\{\pt\}$, which is indicated on the left in Figure \ref{qFH trajectories}. In that case $\ul{Z}$ has no seam $\sigma_r$ between $S_r$ and $S_0$ but true boundary components $(0,\R\times\{0\})$ and $(r,\R\times\{\delta_r\})$.
The ends of the quilted surface are the incoming $\ul{e}_-=\bigl((0,e_-),(1,e_-),\ldots,(r,e_-)\bigr)$ and
the outgoing $\ul{e}_+=\bigl((0,e_+),(1,e_+),\ldots,(r,e_+)\bigr)$.
Note however that the perturbation data $(\ul{J},\ul{K})$ is $\R$-invariant and the count for the Floer differential is modulo simultaneous $\R$-shift of all maps $u_k$. That is, unlike in the definition of relative quilt invariants in \cite{quilts}, where no symmetries are divided out and index $0$ solutions are counted, we here count the isolated solutions $\M(\ul{x}^-,\ul{x}^+)_0 = \M_{\ul{Z}}(\ul{x}^-,\ul{x}^+)_1/\R $, which are pseudoholomorphic quilts of index $1$.
\end{remark}

\begin{proposition}  \label{width indep}
$HF(\ul{L})$ is independent, up to isomorphism of $\Z_N$-graded
groups, of the choice of perturbation data $(\ul{H},\ul{J})$ and
widths $\ul{\delta}$ of the strips.
\end{proposition}

\begin{proof}  Suppose that $(\ul{H}^i, \ul{J}^i, \ul{\delta}^i)$
are two different choices for $i=0,1$.  
Let $\ul{Z}_{01}$ resp.\ $\ul{Z}_{10}$ be the quilted cylinder as before, 
but with complex structures $j_k$ on each strip $S_k\cong\R\times[0,1]$
that interpolate between the two widths $\delta_k^0$ and $\delta_k^1$ at the ends $(k,e_-)$ and $(k,e_+)$, 
in this order for $\ul{Z}_{01}$ and in reversed order for $\ul{Z}_{10}$.
In order for the seams to be real analytic we pick the standard complex structure near the boundary components $\R\times\{0,1\}\subset\partial S_k$ and only in the interior of $S_k$ scale to the appropriate width and interpolate.
Figure \ref{interp} shows the example for $r=3$ and $M_0=M_4=\{\pt\}$.
For readers unfamiliar with \cite{quilts} we specify, as an example, that the 
{\em quilted surface} $\ul{Z}_{01}$ consists of the following data:
\ben 
\item A collection $\ul{S} = (S_k)_{k=0,\ldots,r}$ of {\em patches}.
Here these are the twice punctured disks $S_k=\R\times[0,1]\cong D^2\setminus\{-1,1\}$ 
with a complex structure $j_k$ and strip-like ends as follows:
We fix embeddings of half-strips 
$\eps_\pm : \R^\pm \times [0,\delta_\pm] \to S_k , (s,t) \mapsto (\pm 1 + \delta_\pm^{-1} s, \delta_\pm^{-1} t)$.
of width $\delta_-=\delta^0_k$ resp.\ $\delta_+=\delta^1_k$.
These are incoming resp.\ outgoing strip-like ends with disjoint images.
To construct the complex structure we moreover 
let $\Delta=\tfrac13 \min(\delta_+,\delta_-)$ and fix disjoint embeddings
$\eps_0: \R \times [0,\Delta] \to S_k$ and $\eps_1: \R \times [1-\Delta,\Delta] \to S_k$ 
such that $\eps_0(s,0)=(s,0)$ resp.\ $\eps_1(s,1)=(s,1)$ and $\eps_{0/1}^{-1}\circ\eps_\pm$
are four biholomorphisms with respect to the canonical complex structure on the half-strip and strip.
Now we choose a complex structure $j_k$ on $S_k$ such that it pulls back to the canonical complex structure under each of $\eps_-,\eps_+,\eps_0,\eps_1$.
\item 
A collection $\S$ of {\em seams}. Here these are the set $\S=(\sigma_k)_{k=0,\ldots,r}$ of seam labels 
$\sigma_k=\{(k,\R\times\{1\}\subset\partial S_k),(k+1,\R\times\{0\}\subset\partial S_{k+1})\}$, where we use the label $k$ modulo $(r+1)$, 
with seam maps $\phi_{\sigma_k}: \partial S_k \ni (s,1)\mapsto (s,0)\in\partial S_{k+1}$
providing diffeomorphisms of boundary components.
By our construction these seams are {\it real analytic} and {\it compatible with strip-like ends}.
\item
Orderings of patches, boundary components, and quilted ends. Here they are given by the enumeration with $k=0,\ldots,r$, and there is only one incoming and one outgoing quilted end consisting of all incoming resp.\ outgoing strip-like ends on the patches.
\een
\begin{figure}[ht] 
\begin{picture}(0,0)%
\includegraphics{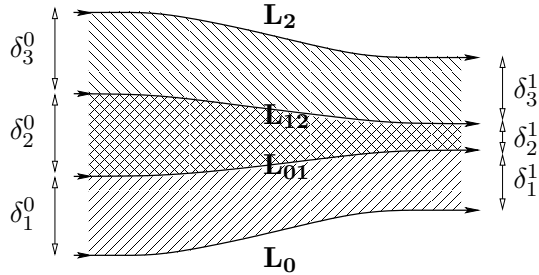}%
\end{picture}%
\setlength{\unitlength}{2072sp}%
\begingroup\makeatletter\ifx\SetFigFont\undefined%
\gdef\SetFigFont#1#2#3#4#5{%
  \reset@font\fontsize{#1}{#2pt}%
  \fontfamily{#3}\fontseries{#4}\fontshape{#5}%
  \selectfont}%
\fi\endgroup%
\begin{picture}(6030,3024)(-134,-2758)
\put(-134,-431){\makebox(0,0)[lb]{{$\delta^0_3$}}}
\put(-134,-1421){\makebox(0,0)[lb]{{$\delta^0_2$}}}
\put(-134,-2411){\makebox(0,0)[lb]{{$\delta^0_1$}}}
\put(5896,-926){\makebox(0,0)[lb]{{$\delta^1_3$}}}
\put(5896,-1511){\makebox(0,0)[lb]{{$\delta^1_2$}}}
\put(5896,-2006){\makebox(0,0)[lb]{{$\delta^1_1$}}}
\put(2900,0){\makebox(0,0)[lb]{{$\mathbf{L_2}$}}}
\put(2900,-1191){\makebox(0,0)[lb]{{$\mathbf{L_{12}}$}}}
\put(2900,-1800){\makebox(0,0)[lb]{{$\mathbf{L_{01}}$}}}
\put(2900,-2911){\makebox(0,0)[lb]{{$\mathbf{L_0}$}}}
\end{picture}
\caption{Interpolating between two widths}
\label{interp}
\end{figure}
We moreover interpolate the perturbation data on the two ends by some 
regular $(\ul{K}_{01},\ul{J}_{01})$ on $\ul{Z}_{01}$ and similar for $\ul{Z}_{10}$.
Summarizing from \cite{quilts} this means the following.
\ben
\item 
We pick any tuple of function valued one-forms
$\ul{K}_{01}= \bigl( K_k \in \Omega^1(S_k,C^\infty(M_k)) \bigr)_{k=0,\ldots,r} $
such that $K_k|_{\pd S_k}=0$ and on each end 
$\eps_{-}^*K_k = H_{k}^0 \d t$ resp.\ $\eps_{+}^*K_k = H_{k}^1 \d t$.
We denote the corresponding Hamiltonian vector field
valued one-forms by $\ul{Y}_{01}=\bigl( Y_k \in \Omega^1(S_k,\Vect(M_k))$. 
These satisfy $\eps_{\pm}^* Y_k = X_{H_{k}^{0/1}} \d t$ on each strip-like end.
\item
Let $\J$ denote the set of collections
$$ \ul{J} = \bigl( J_k \in \Map(S_k,\J(M_k,\omega_k)) \bigr)_{k =0,\dots,r } $$
of compatible almost complex structures
agreeing with the almost complex structures $J_k^0$ resp.\ $J_k^1$ on the incoming ends $\im\eps_-$
resp.\ the outgoing ends $\im\eps_+$.
We prove in \cite{quilts} that $\J$ contains an open dense subset of regular $\ul{J}$, of which we pick one 
$\ul{J}_{01}$.
\een
Given this perturbation data and intersection points $\ul{x}^-\in\cI(\ul{L})^0$,
$\ul{x}^+\in\cI(\ul{L})^1$
we construct in \cite{quilts} the moduli spaces of pseudoholomorphic quilts 
$\M_{\ul{Z}_{01}}(\ul{x}^-,\ul{x}^+)$ and $\M_{\ul{Z}_{10}}(\ul{x}^-,\ul{x}^+)$. 
For example,
$$
\M_{\ul{Z}_{01}}(\ul{x}^-,\ul{x}^+) := \bigl\{ \ul{u}=\bigl( u_k :S_k \to M_k \bigr)_{k=0,\ldots,r} \,\big|\, (a)-(d) \bigr\}
$$ 
is the space of collections of $(\ul{J}_{01},\ul{Y}_{01})$-holomorphic maps with 
Lagrangian seam conditions, finite energy, and fixed limits, that is
\ben
\item
$J_k(u_k) \circ (\d u_k - Y_k(u_k)) -
(\d u_k - Y_k(u_k)) \circ j_k = 0$ for $k=0,\ldots,r$,
\item  
$(u_{k},u_{k+1}\circ\varphi_{\sigma_k}) (\R\times\{0\}) \subset L_{k(k+1)}$
for all $k=0,\ldots,r$,
\item $\sum_{k=0}^r \int_{S_k} \bigl( u_k^*\omega_k + \d (K_k\circ u_k ) \bigr) <\infty$,
\item $\lim_{s \to \pm \infty} u_k(\eps_{\pm}(s,t)) = x^\pm_{k}(t) $ 
for all $k=0,\ldots,r$.
\een
The relative invariants, constructed in \cite{quilts} 
from the zero-dimensional moduli spaces of pseudoholomorphic quilts
then provide maps between the Floer cohomology groups
$$
\Phi_{\ul{Z}_{01}} : HF(\ul{L})^0 \to HF(\ul{L})^1, \qquad\quad
\Phi_{\ul{Z}_{10}} : HF(\ul{L})^1 \to HF(\ul{L})^0 .
$$ 
We briefly review the construction: On chain level a map $ C\Phi_{\ul{Z}_{01}} : CF(\ul{L})^0 \to  CF(\ul{L})^1$
can be defined by 
\begin{equation*} \label{defrel}
 C\Phi_{\ul{Z}_{01}}\bra{\ul{x}^-} := \sum_{\ul{x}^+\in\cI(\ul{L})^1} 
 \biggl( \sum_{\ul{u} \in \M_{\ul{Z}_{01}}(\ul{x}^-,\ul{x}^+)_0} \eps(u) \biggr) 
 \bra{\ul{x}^+} ,
\end{equation*}
where $\eps:\M_{\ul{Z}_{01}}(\ul{x}^-,\ul{x}^+)_0 \to \{-1,+1\}$
is defined by comparing the orientation constructed in \cite{orient} on the zero dimensional component of the moduli space to the canonical orientation of a point. 
We prove in \cite{quilts} that the map $C\Phi_{\ul{Z}_{01}}$ is a chain map and so descends 
to a map of Floer cohomologies. In fact, the map on cohomology level is independent of the choice
of perturbation data $(\ul{K}_{01},\ul{J}_{01})$.

Next, the quilted surface $\ul{Z}_{01}\#\ul{Z}_{10}$ that is obtained by gluing the incoming ends of
$\ul{Z}_{01}$ to the outgoing ends of $\ul{Z}_{10}$ 
can be deformed with fixed ends to the infinite strip with translationally invariant complex structures
(reflected in the widths $\ul{\delta}^1$) and perturbation data $(\ul{H}^1,\ul{J}^1)$. 
The relative quilt invariant defined by the latter is the identity on $HF(\ul{L})^1$ since only constant strips can contribute (all nonconstant solutions lie in at least $1$-dimensional moduli spaces due to the nontriviality of the $\R$-action), see \cite{quilts}.
Since the relative quilt invariants are independent of the above choices, the relative invariant $\Phi_{\ul{Z}_{01}\#\ul{Z}_{10}}$ is the identity on $HF(\ul{L})^1$ (and similarly $\Phi_{\ul{Z}_{10}\#\ul{Z}_{01}}$ is the identity on $HF(\ul{L})^0$).
Then by the gluing theorem for relative quilt invariants \cite{quilts} (where the sign is positive) we have
$$
\Phi_{\ul{Z}_{01}}\circ \Phi_{\ul{Z}_{10}} = \Phi_{\ul{Z}_{01}\#\ul{Z}_{10}}
={\rm Id} , \ \ \ \
\Phi_{\ul{Z}_{10}}\circ \Phi_{\ul{Z}_{01}} = \Phi_{\ul{Z}_{10}\#\ul{Z}_{01}}
={\rm Id}.
$$ 
This proves that the Floer cohomology groups $HF(\ul{L})^0$ and
$HF(\ul{L})^1$ arising from the different choices of data are
isomorphic.\end{proof}

\subsection{Geometric composition and quilted Floer cohomology}
\label{shrink}

In this section we prove and discuss the isomorphism \eqref{eq:iso}, more precisely stated as follows.

\begin{theorem} \label{main2}  
Let $\ul{L}= (L_{01},\ldots,L_{r(r+1)})$ be a cyclic sequence of 
Lagrangian correspondences between
symplectic manifolds $M_0,\ldots,M_{r+1}=M_0$ as in
Definition~\ref{Lag cycle}.  Suppose
\begin{enumerate} 
\item the symplectic manifolds all satisfy (M1-2) with 
the same monotonicity constant $\tau$,
\item  the Lagrangian correspondences all satisfy (L1-2) and at least one of the following: 
\begin{enumerate}
\item Each Lagrangian correspondence satisfies (L3).
\item The sum ${w(L_{01})+\ldots+w(L_{(j-1)j})+w(L_{j(j+1)})+\ldots +w(L_{r(r+1)})=0}$ of holomorphic disk counts  \eqref{diskcount} vanishes.
\item The sum ${w(L_{01})+\ldots+w(L_{(j-1)j}\circ L_{j(j+1)})+\ldots +w(L_{r(r+1)})=0}$ of holomorphic disk counts  \eqref{diskcount} vanishes.
\end{enumerate}
\item  the sequence $\ul{L}$ is monotone, relatively spin, and graded in
the sense of Section~\ref{sequences},
\item the composition $L_{(j-1)j} \circ L_{j(j+1)}$ is embedded in the sense of Definition \ref{embedded}.
\end{enumerate}
Then with respect to the induced relative spin structure, orientation,
and grading \footnote{ The grading of $L_{(j-1)(j+1)}$ is given by \eqref{dfn comp grad}, the orientation is given by Remark~\ref{rmk:embedded}(b), and for the relative spin structure see \cite{orient}.  } 
on the modified sequence $\ul{L}'=(L_{01},\ldots, L_{(j-1)j} \circ L_{j(j+1)}, \ldots,L_{r(r+1)})$
the following two Floer homologies are well defined and canonically isomorphic as graded groups:
\begin{equation*}
HF(\ul{L})= HF(\ldots L_{(j-1)j} , L_{j(j+1)} \ldots)
\overset{\sim}{\longrightarrow} HF(\ldots L_{(j-1)j} \circ L_{j(j+1)}
\ldots)=HF(\ul{L}') .
\end{equation*}
The isomorphism is canonical in the following sense: The Floer cohomologies defined by any two choices of perturbation data and widths are canonically isomorphic. For sufficently small width $\delta_j>0$ and corresponding perturbation data, the isomorphism $HF(\ul{L})\cong HF(\ul{L}')$ is given by the identity map on the generators $\cI(\ul{L})=\cI(\ul{L}')$, which are canonically identified by Remark~\ref{rmk:canonical ident H}.
\end{theorem} 

Before summarizing the proof let us mention the (im)possibility of various generalizations.

\begin{remark} \label{main rmk}
\begin{enumerate}
\item
Note the orientation conventions when comparing with the Floer cohomology $HF(L_0,L_1)$ 
for a pair of Lagrangians $L_0,L_1\subset M$ in  \cite{oh:fl1}.
For quilted Floer cohomology, this pair is viewed as cyclic sequence $\{pt.\}\to M \to \{pt.\}$, that is the Lagrangian correspondences are $L_0\subset \{pt.\}^-\times M\cong M$ and $L_1^- \subset M^-\times \{pt.\} \cong M^-$, the same submanifold but viewed as Lagrangian with respect to $-\omega$.
Thus we obtain $w(L_0) + w(L_1^-) = w(L_0) - w(L_1)$, since the $-J_1$-holomorphic discs with boundary on 
$L_1^- \subset M^-$ are identified with $J_1$-holomorphic discs with boundary on $L_1\subset M$ via a reflection of the domain, which is orientation reversing for the moduli spaces.
\item
If in Theorem~\ref{main2}~(b) we only assume (L1-2) then Floer cohomology may not be well-defined due to $\partial^2=w{\rm Id}$, where 
$w = w(L_{01})+\ldots+w(L_{r(r+1)})$, see \cite{oh:fl1,fieldb}.
In case $w\neq 0$ the chain homotopy equivalence
$(CF(\ul{L}),\partial) \cong (CF(\ul{L}'),\partial')$ 
continues to hold in the derived category of matrix factorizations\footnote{
Explicitly, there exist chain maps
$f: CF(\ul{L}) \to CF(\ul{L}')$ and $g: CF(\ul{L}')\to CF(\ul{L})$ such that both 
$f\circ g$ and $g\circ f$ are homotopy equivalent to the identity in the sense that
e.g.\  $ {\rm Id} - g\circ f = k\partial - \partial k$ for some $k: CF(\ul{L}) \to CF(\ul{L})$.
}, see \cite{fieldb}.

We will show that either of the extra assumptions (i), (ii), or (iii) implies that $\partial^2=0$ on both Floer complexes, and hence implies both (ii) and (iii) (but not (i)).
In fact, (i) directly implies (ii) (but not (iii)) since all minimal Maslov numbers being at least $3$ implies $w(L_{(i-1)i})=0$ for each $i$. The other implications $(i)\Rightarrow (iii)$,
$(ii)\Rightarrow (iii)$, $(iii)\Rightarrow (ii)$ also require the monotonicity and embeddedness assumptions since they follow indirectly from the above isomorphism in the derived category.
\item
The relative spin structures are only needed to define the Floer
cohomology groups with $\Z$ coefficients.  Here we only prove the
isomorphism with $\Z_2$ coefficients.  The full result then follows
from a comparison of signs in \cite{orient}.
\item
There should also be versions of this result for Floer cohomology with
coefficients in flat vector bundles, and Novikov rings, using an
understanding of their behaviour under geometric composition, similar
to the theory presented here for gradings.  The gradings on the
Lagrangians can be dropped if one wants only an isomorphism of
ungraded groups.
\item Note that the geometric composition $L_{(j-1)j} \circ L_{j(j+1)}$
could be a smooth Lagrangian despite the composition not being
embedded. If this failure is in the transversality, then our approach
does not apply (as e.g. for a $G$-invariant Lagrangian $L\subset
\mu^{-1}(0)$ in the zero set of the moment map, whose composition with
$\Sigma_{\mu}$ is the smooth projection $L\circ\Sigma_{\mu}=\pi(L)$
despite $L$ not being transverse to $\mu^{-1}(0)$). (For such
Lagrangians one would expect a correspondence between holomorphic
curves in $M\qu G$ and symplectic vortices in $M$, in the spirit of
\cite{ga:gw} and the Lagrangian version of the Atiyah-Floer conjecture
\cite{we:su}.)  However, when $L_{(j-1)j} \times_{M_j} L_{j(j+1)}$
is transverse but a $k$-fold cover of $L_{(j-1)j} \circ L_{j(j+1)}$,
then the map of intersection points $\cI(\ul{L}) \to \cI(\ul{L}')$ is
a $k$-to-$1$ map as well. In this case our analysis still applies and
gives a $k$-to-$1$ map of moduli spaces as long as bubbling is excluded. 
This leads to further calculation tools for Floer cohomology but 
needs to be investigated on a case-by-case basis.
\item
The monotonicity assumptions (M1) and (L1) cannot simply be replaced by other tools which allow the definition of Floer cohomology (such as Novikov rings, twisted coefficients, obstructions, or deformations).
This is since a new type of bubbling can occur in the strip shrinking that we use to prove the isomorphism. We have called it the ``figure eight bubble" and describe it in \cite{isom}. However, we are lacking the construction of a moduli space of figure eight bubbles. Our present method for excluding these bubbles hinges on strict monotonicity with nonnegative constant $\tau\geq 0$ as well as the $2$-grading assumption implied by orientations. 
In general, we expect figure eight bubbles to be a codimension $1$ phenomenon in a $1$-parameter family of strip widths approaching zero. 
We hence expect the isomorphism to fail in more general settings, except for special topological assumptions restricting the expected dimension of the moduli space of figure eight bubbles.
Eventually, we expect to construct obstruction classes and an $A_\infty$-type structure from moduli spaces of figure eight bubbles, and to replace the isomorphism by a morphism of $A_\infty$-modules. However, all of this depends on a basic removable singularity result for figure eight bubbles, which has not yet been accomplished.
\end{enumerate}
\end{remark}

Theorem~\ref{main2} is fairly obvious if one of the composed Lagrangians correspondences is the graph of a symplectomorphism. 
It suffices to observe that symplectomorphisms map pseudoholomorphic curves to pseudoholomorphic curves. However, there is no corresponding effect for more general Lagrangian correspondences. Here the natural approach to a proof is to degenerate the holomorphic curve equation in $M_j$ until solutions become constant across the strip (or, equivalently, shrink the width of that strip to zero). This limit corresponds to geometric composition of the two Lagrangian correspondences attached to the strip.
Clearly, most difficulties in this proof are localized near the degenerating strip.
We thus banished the analysis to \cite{isom}, where we prove the special case
$HF(L_0,L_{01},L_{12},L_2)\overset{\sim}{\to} HF(L_0,L_{01} \circ L_{12},L_2)$ of Theorem~\ref{main2} by establishing a bijection between the
Floer trajectories for $(L_0,L_{02},L_2)$ on strips of width $(1,1)$
and those for $(L_0,L_{01},L_{12},L_2)$ on strips of width $(1,\delta,1)$
for sufficiently small width $\delta$ of the middle strip. 
These quilted Floer trajectories are shown in Figure \ref{shrinkfig}.
\begin{figure}[ht]
\begin{picture}(0,0)%
\includegraphics{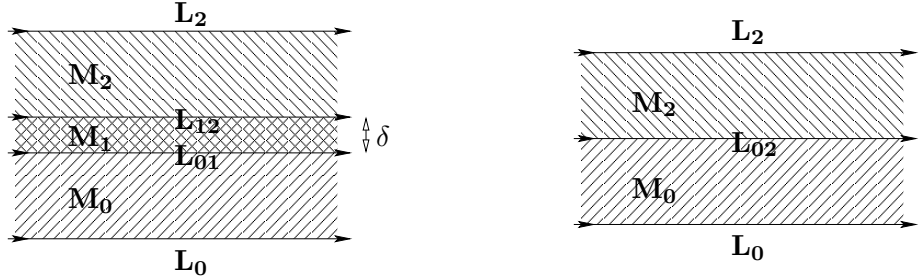}%
\end{picture}%
\setlength{\unitlength}{2200sp}%
\begingroup\makeatletter\ifx\SetFigFont\undefined%
\gdef\SetFigFont#1#2#3#4#5{%
  \reset@font\fontsize{#1}{#2pt}%
  \fontfamily{#3}\fontseries{#4}\fontshape{#5}%
  \selectfont}%
\fi\endgroup%
\begin{picture}(9549,3135)(589,-2461)
\put(2476,504){\makebox(0,0)[lb]{{$\mathbf{L_2}$}}}
\put(2476,-700){\makebox(0,0)[lb]{{$\mathbf{L_{12}}$}}}
\put(2476,-1111){\makebox(0,0)[lb]{{$\mathbf{L_{01}}$}}}
\put(2476,-2300){\makebox(0,0)[lb]{{$\mathbf{L_0}$}}}
\put(8751,-2130){\makebox(0,0)[lb]{{$\mathbf{L_0}$}}}
\put(8751,-961){\makebox(0,0)[lb]{{$\mathbf{L_{02}}$}}}
\put(8751,300){\makebox(0,0)[lb]{{$\mathbf{L_2}$}}}
\put(7626,-511){\makebox(0,0)[lb]{{$\mathbf{M_2}$}}}
\put(7626,-1486){\makebox(0,0)[lb]{{$\mathbf{M_0}$}}}
\put(1276,-1600){\makebox(0,0)[lb]{{$\mathbf{M_0}$}}}
\put(1276,-900){\makebox(0,0)[lb]{{$\mathbf{M_1}$}}}
\put(1276,-200){\makebox(0,0)[lb]{{$\mathbf{M_2}$}}}
\put(4750,-850){\makebox(0,0)[lb]{{$\delta$}}}
\end{picture}%
\caption{Shrinking the middle strip} \label{shrinkfig}
\end{figure}
The missing piece of proof in \cite{isom} is the independence of the Floer cohomology from the choice of $\delta>0$, which we here established in Proposition~\ref{width indep}.

\begin{proof}[Summary of proof of Theorem~\ref{main2}.]
We first consider the case of assumption (i) or (ii) holding in (b). Then the assumptions of the Theorem guarantee that $HF(\ul{L})$ is well-defined. By Remark~\ref{L monotone} the monotonicity of $\ul{L}$ also implies monotonicity of $\ul{L}'$ and hence monotonicity in the sense of (L1) for $L_{(j-1)j} \circ L_{j(j+1)}$. (Assuming the symplectic manifolds to be connected, any disk can be extended to a quilted cylinder.) Compactness and orientation (L2) also holds for the composed correspondence, but the minimal Maslov index condition (L3) may not transfer.
However, this only affects the question whether $\partial^2=w{\rm Id}=0$ on the Floer chain group for $\ul{L}'$. In fact, we just assumed that $w= w(L_{01})+\ldots+w(L_{r(r+1)}) =0$ on $CF(\ul{L})$, and a priori it is not clear that one should have $w(L_{(j-1)j} \circ L_{j(j+1)})=w(L_{(j-1)j}) + w(L_{j(j+1)})$.
If, on the other hand, assumption (iii) holds in (b), then we are guaranteed that $HF(\ul{L}')$ is well-defined, but it is not a priori clear that $\partial^2=0$ on $CF(\ul{L})$.

In either case, to define the differential on $CF(\ul{L})$ we choose some widths $\ul{\delta}'$, Hamiltonian perturbations
$\ul{H}'$ to make the intersection $\cI(\ul{L}')$ transverse, and almost complex structures $\ul{J}'$ to make the moduli spaces of Floer trajectories for $\ul{L}'$ regular.
Thanks to Proposition~\ref{width indep} we may then choose the same widths $\ul{\delta}$ except for some small $\delta_j>0$, the same Hamiltonian perturbations $\ul{H}$ except for the additional $H_j\equiv 0$, and the same almost complex structures $\ul{J}$ except for some additional time-independent $J_j\in\J(M_j,\omega_j)$, 
to define $HF(\ul{L})$. We only need to make sure that this choice makes the intersection points $\cI(\ul{L})$ and the moduli spaces of Floer trajectories for $\ul{L}$ regular. The first is automatically the case by the transversality assumption for $L_{(j-1)j} \times_{M_j} L_{j(j+1)}$, the latter is true for $\delta_j>0$ sufficiently small and is proven as part of the adiabatic limit analysis \cite{isom}.
(Actually, precisely following the constructions of \cite{isom}, we can achieve transversality for $\ul{L}'$ with $J'_{j-1}$ and $J'_{j+1}$ being time-independent near the seam; then $J_{j-1}$ and $J_{j+1}$ are obtained by a slight linear dilation and constant extension near the new seams.)
With these choices, our assumption (b) on holomorphic disk counts (or (L3) on the minimal Maslov index) implies $\partial^2=0$ on $CF(\ul{L})$. 

Next, the injectivity assumption for the composition $L_{(j-1)j} \times_{M_j} L_{j(j+1)}$ provides a canonical bijection of generalized intersection points $\cI(\ul{L}')\cong\cI(\ul{L})$ as in Remark~\ref{rmk:canonical ident}. In \cite{isom} we establish bijections between the corresponding moduli spaces of Floer trajectories for $\delta_j>0$ sufficiently small. This means that the Floer differentials on $CF(\ul{L}')$ and  $CF(\ul{L})$ agree under the canonical identification of generators. In particular that implies $\partial^2=0$ on both complexes as soon as it is true on one (which is ensured by each version of assumption (b)). Hence both Floer cohomologies are well-defined and isomorphic as claimed.
(In fact, we deduce a posteriori that $w(L_{(j-1)j} \circ L_{j(j+1)})=w(L_{(j-1)j}) + w(L_{j(j+1)})$.)
\end{proof}

\begin{remark}  To see that the assumption that the composition $L_{(j-1)j}\circ L_{j(j+1)}$ is embedded is necessary, consider the case that $r=2$ and $M_0,M_2$ are points. In this case, if $v:\R\times[0,1]\to M_1$ is a Floer trajectory of index one
with limits $x^+\neq x^-$, we can consider the rescaled maps
$v_\delta:\R\times[0,\delta]\to M_1$.  In this case a figure eight
bubble always develops in the limit $\delta \to 0$. 
This shows that the bijection between trajectories fails in this case.
\end{remark}

\section{Applications} \label{sec:app}

Quilted Floer homology was originally designed to construct symplectic versions of
gauge theoretic invariants, in particular symplectic versions of
Donaldson invariants, which we develop in later papers \cite{field,
  fieldb}, Seiberg-Witten invariants as in Perutz \cite{per:lag} and
Lekili \cite{yanki}, and Khovanov invariants as in Rezazadegan \cite{reza}. Applications to symplectic topology are given for moduli spaces of
flat bundles in \cite{fieldb}, and to classification of Lagrangians in
tori in Abouzaid-Smith \cite{as:hms}.  In this section,
we give a few brief applications of the main result to symplectic topology, which show how the results work in practice. 
Most of the concrete examples can be achieved with
other, less sophisticated methods. Our point in giving them is to show how many Floer homology calculations can be obtained from a single principle: {\em Floer homology is well defined under embedded geometric composition}. 
The reader looking for more sophisticated applications is encouraged to look at the sequel papers and the references above.

\subsection{Direct computation of Floer cohomology}

\begin{theorem} \label{thm1}
Let $L_{01}\subset M_0^-\times M_1$ be a Lagrangian correspondence and
suppose that the Lagrangian submanifolds $L_0\subset M_0$ and
$L_1\subset M_1$ are such that both $L_0\circ L_{01}$ and $L_{01}\circ
L_1$ are embedded compositions. 
Assume that $M_0,M_1$ satisfy (M1-2), $L_0,L_1,L_{01}$ satisfy (L1-2), 
and $(L_0\times L_1, L_{01})$ is a monotone pair in the sense of Definition~\ref{tuple monotone}~(b).
Then there exists a canonical isomorphism
$$
CF(L_{0}\circ L_{01},L_{1})\overset{\sim}{\longrightarrow} CF(L_{0},L_{01}\circ L_{1}) 
$$
in the category of derived matrix factorizations. It is an isomorphism of Floer cohomologies if one of the following holds: All Lagrangians satisfy (L3), or $w(L_0)+w(L_{01})+w(L_1)=0$, or
$w(L_0\circ L_{01})+w(L_1)=0$, or $w(L_0)+w(L_{01}\circ L_1)=0$.
\end{theorem}
\begin{proof}
By Theorem~\ref{main2} both Floer cohomologies are isomorphic to the
quilted Floer cohomology $HF(L_0,L_{01},L_1)=HF(L_0\times L_1,L_{01})$.
In case $\partial^2\neq 0$, the isomorphism of homologies is replaced by a chain homotopy equivalence in the derived category, see Remark~\ref{main rmk}.
\end{proof}

\begin{example} 
We begin with a ``warm-up'' example.  Let $N$ be a compact,
simply-connected, monotone symplectic manifold. The submanifold
$\Delta_{ij} := \{ (x_1,x_2,x_3,x_4) | x_i=x_j\} \subset N^-\times N
\times N^-\times N$ is coisotropic for appropriate choices of $1\leq
i<j\leq 4$.  Then we can identify
\begin{equation} \label{identify}
 HF( \Delta_{14}\cap\Delta_{23}, \Delta_{12}\cap\Delta_{34} ) \cong
 HF(\Delta_N, \Delta_N) \cong H(N) \end{equation}
with the homology of $N$.  This follows from Theorem~\ref{thm1}
applied to $L_0 = \Delta_{N} \subset N^-\times N = M_0$, $L_1 =
\Delta_{12}\cap\Delta_{34} \subset N^-\times N \times N^-\times N =
M_1$, and $L_{01} = \{ (w,z;w,x,x,z) | w,x,z\in N \} \subset
M_0^-\times M_1$.  Then the compositions $L_0\circ L_{01} =
\Delta_{14}\cap\Delta_{23}$ and $L_{01}\circ L_1 = \Delta_{N} $ are
clearly embedded.  Monotonicity together with simply connectedness
ensures the monotonicity of all the Lagrangians and pairs of
Lagrangians. Since $N$ is orientable, all minimal Maslov indices are
at least $2$.  The reader can easily verify the identification
\eqref{identify} using the fact that the components of a holomorphic
trajectory for $(\Delta_{14} \cap \Delta_{23},\Delta_{12} \cap
\Delta_{34})$ fit together to a holomorphic cylinder $v: S^1 \times
      [0,1] \to N$.
\end{example}

The following is a more non-trivial example of Theorem~\ref{thm1}.

\begin{example}   \label{mn}
Let $M_n$ be the moduli space of Euclidean $n$-gons of edge length
$1$, as in for example Kirwan \cite{ki:coh}:
$$M_n = (S^2)^n \qu SO(3) = \{ (v_1 ,\ldots, v_n ) \in (S^2)^n |
v_1 + \ldots + v_n = 0 \} / SO(3) .$$ 
We take on $S^2$ the standard symplectic form $\omega$ with volume $ 4
\pi$ so that $c_1(S^2) = [\omega]$.  For $n \ge 5$ odd $M_n$ is a
monotone symplectic manifold with minimal Chern number $1$ and monotonicity
constant $1$.
For example, $M_3$ is a point and $M_5$ is diffeomorphic to the fourth
del Pezzo surface, given by blowing-up of $\P^2$ at four points
\cite[Example~1.10]{se:dehn}.  For $i\neq j$ the submanifold $\Delta_{ij} = \{
     [v_1,\ldots,v_n] \in M_n \,|\, v_i = -v_j \}$ is a coisotropic,
     spherically fibered over $M_{n-2}$ by the map that forgets
     $v_i,v_j$.  The image of $\Delta_{ij}$ in $M_{n-2}^- \times M_n$
     is a Lagrangian correspondence, also denoted $\Delta_{ij}$.
For $i,j,k$ distinct the composition $\Delta_{ij} \circ \Delta_{jk}^t
$ is embedded and yields the graph of a permutation on $M_{n-2}$.  For
$k=i\pm 1$ or for $i,k=j\pm 1$ this permutation is trivial, so we have
$$\Delta_{ij} \circ \Delta_{j(i\pm1)}^t =\Delta_{M_{n-2}}, \quad \Delta_{(j\pm1)j} \circ \Delta_{j(j\mp1)}^t =\Delta_{M_{n-2}} .$$
Now let $L \subset M_{n-2}$ be a compact, oriented, monotone
Lagrangian, and $L_{ij} = L \circ
\Delta_{ij}$ be its inverse image in $M_n$.  This composition is
embedded and we can also identify it with $L_{ij} = \Delta_{ij}^t
\circ L$.  The latter allows to calculate 
$$\Delta_{ij} \circ L_{jk} = \Delta_{ij}\circ\Delta_{jk}^t \circ L .
$$ 
For $i,j,k$ distinct it is an embedded composition, which yields the
image of $L$ under permutation.  Suppose that the pair $(L,L)$ is
monotone, so that $HF(L,L)$ is well-defined.  Using Theorem~\ref{thm1}
we can also calculate
\begin{eqnarray*} 
HF(L_{ij}, L_{j(i\pm 1)}) 
= HF(L \circ \Delta_{ij}, L_{j(i\pm 1)}) 
= HF(L, \Delta_{ij} \circ L_{j(i-1)})
= HF(L,L) .
\end{eqnarray*}
Here all Floer homologies are defined since $\partial^2=0$ on $CF(L,L)$.
Similarly, we obtain $HF(L_{(j\pm1)j}, L_{j(j\mp 1)})= HF(L,L)$ .
\end{example}  

The symplectic manifolds in the above example are certain moduli spaces of parabolic bundles on punctured spheres, which can also be seen as $SU(2)$-representation spaces of punctured spheres with fixed conjugacy class for each puncture \cite{MeWo}.
The Lagrangian submanifolds $\Delta_{ij}$ arise from elementary tangles connecting two punctures, by taking parabolic bundles on the corresponding cobordism of punctured spheres, and restricting them to the boundary. The fact that their composition is the graph of a permutation reflects the fact that all braid moves on the punctures can be decomposed into elementary tangles. Thus one may hope to define knot or tangle invariants by decomposition into elementary tangles and representation in the symplectic category as above. The invariant should then be the quilted Floer cohomology of the sequence of Lagrangian correspondences arising from the sequence of elementary tangles.
In order to prove invariance, one needs to check that moves between different decompositions are reflected by isomorphisms of Floer cohomology. The above example shows how this follows from our main Theorem for the cancellation of two elementary tangles.

This general approach to defining topological invariants is also described in the introduction. The moduli spaces in Example~\ref{mn} actually yield the trivial invariant since they arise from conjugacy classes close to the center of $SU(2)$.
In \cite{field} and \cite{fieldb} we employ similar moduli spaces with different conjugacy classes, giving rise to nontrivial invariants of tangles and $3$-manifolds.

\subsection{Computations in $\CP^n$}
\label{cpn}

In this section we demonstrate, at the example of $\CP^n$, how some Floer cohomologies in toric symplectic varieties can be calculated by reduction.
We equip $\CP^n=\{[\ul{z}]=[z_0:z_1:\ldots:z_n]\}$ with Fubini-Study symplectic form and moment maps $\mu_j([\ul{z}])=\pi |z_j|^2/|\ul{z}|^2$ for $j=1,\ldots, n$. 
We denote by 
$$
\Sigma_{j}:=\bigl\{ \bigl( [\ldots z_{j-1}:z_{j+1}:\ldots], [z_0:\ldots:z_n]\bigr) \,\big|\, \mu_j( [z_0:\ldots:z_n])=\tfrac \pi {n+1} \bigr\} 
\subset (\CP^{n-1})^-\times \CP^n
$$ 
the Lagrangian sphere arising from reduction at the level set 
$$
\mu_j^{-1}(\tfrac \pi {n+1})=\bigl\{ [z_0:\ldots:z_n] \,\big|\, z_j=\tfrac 1{\sqrt{n+1}}, {\textstyle \sum_{i\neq j}} |z_i|^2 = \tfrac{n}{n+1} \bigr\}  .
$$
Note that the reduced space, e.g.\ $\mu_n^{-1}(\tfrac \pi {n+1})/S^1 = \{ [z_0:\ldots z_{n-1}] \,\big|\, {\textstyle \sum_i} |z_i|^2 = \tfrac{n}{n+1} \}/S^1$ for $j=n$,
is $\CP^{n-1}$ with Fubini-Study form scaled by $\tfrac{n}{n+1}$, hence has the same monotonicity constant $\tau = n^{-1} \frac{n \pi}{n+1} = (n+1)^{-1} \pi$ as $\CP^n$.
(Recall that the generator $\CP^1\subset\CP^n$ of $\pi_2(\CP^n)$ has Fubini-Study symplectic area $\pi$ and Chern number $n+1$.)

More generally, for each $1<k\leq n$ a Lagrangian correspondence 
$$
\Sigma_{(k,\ldots,n)}:=\bigl\{ \bigl( [z_0:\ldots:z_{k-1}], [z_0:\ldots:z_n]\bigr) \,\big|\, \mu_j( [\ul{z}])=\tfrac \pi {n+1} \;\forall j\geq k \bigr\} 
\subset (\CP^{k-1})^-\times \CP^n
$$ 
arises from reduction at the level set 
$$
(\mu_k\times\ldots\times\mu_n)^{-1}(\tfrac \pi {n+1}, \ldots, \tfrac \pi {n+1})
=\Bigl\{ [\ul{z}] \,\Big|\, z_k=|z_{k+1}|=\ldots=|z_n| = \tfrac 1{\sqrt{n+1}} , \sum_{i=0}^{k-1} |z_i|^2 = \tfrac{k}{n+1} \Bigr\}.
$$
Here again the reduced spaces $\CP^{k-1}
=\{ [z_0:\ldots:z_{k-1}] \,|\, \sum_{i} |z_i|^2 = \tfrac{k}{n+1} \}/S^1$ 
carry scaled Fubini-Study forms with monotonicity constant $\tau=\frac{\pi}{n+1}$.
Moreover, note that $\Sigma_{(k,\ldots,n)}$, diffeomorphic to the product 
of an $(n-k)$-torus with a $(2k-1)$-sphere, can be viewed as Lagrangian embedded 
in $(\CP^{k-1})^-\times \CP^n$ and also as coisotropic submanifold of $\CP^n$.
One can check explicitly that the Lagrangians $\Sigma_{(k,\ldots,n)}$ are monotone, and we will see in Corollary~\ref{cor:sphere} below that they are nondisplaceable by Hamiltonian diffeomorphisms. The reason is that as coisotropic they contain the nondisplaceable Clifford torus
$$
T^n_{\rm Cl}=(\mu_1\times\ldots\mu_n)^{-1}(\tfrac \pi {n+1}, \ldots, \tfrac \pi {n+1})
=\bigl\{ [\ul{z}] \,\big|\, z_0=|z_{1}|=\ldots=|z_n| = \tfrac 1{\sqrt{n+1}} \bigr\} \subset \CP^n.
$$ 
That $T^n_{\rm Cl}$ is the only nondisplaceable fibre of the torus fibration is known by e.g.\ \cite{biran-entov-polterovich}. 
Its Floer cohomology was calculated by Cho \cite{cho:hol} with all possible spin structures. Here we reproduce this calculation for the standard spin structure, employing the above Lagrangian correspondences and the isomorphism of Floer cohomology under embedded geometric composition (Theorem~\ref{main2}). 
This approach also allows for a direct computation of Floer cohomology for any pair of nonstandard spin structures on $T^n_{\rm Cl}$, which we will discuss in \cite{orient}.

\begin{theorem}\cite{cho:hol} \label{clifford}
For any $n\in\N$ with the standard spin structure (given by \cite[Prp.8.1]{cho:hol})
$$
HF(T^n_{\rm Cl}, T^n_{\rm Cl})\cong H_*(T^n)\cong \Z^{2^n} .
$$
\end{theorem}
\begin{proof}
The isomorphism between the Floer cohomology and the homology of the Clifford $n$-torus follows inductively from the following chain of isomorphisms:
\begin{align} \label{calc}
HF(T^n_{\rm Cl}, T^n_{\rm Cl})
&= HF(T^1_{\rm Cl}\circ \Sigma_{(2,\ldots,n)}, \Sigma_1^t \circ T^{n-1}_{\rm Cl}) \nonumber \\
&\cong HF(T^1_{\rm Cl}, \Sigma_{(2,\ldots,n)}, \Sigma_1^t,T^{n-1}_{\rm Cl}) \nonumber\\
&\cong HF(T^1_{\rm Cl}, \Sigma_{(2,\ldots,n)}\circ\Sigma_1^t,T^{n-1}_{\rm Cl}) \\
&\cong HF(T^1_{\rm Cl}, T^1_{\rm Cl} \times T^{n-1}_{\rm Cl} ,T^{n-1}_{\rm Cl}) \nonumber\\
&\cong HF(T^1_{\rm Cl}, T^1_{\rm Cl}) \otimes HF(T^{n-1}_{\rm Cl} ,T^{n-1}_{\rm Cl}) .
\nonumber
\end{align}
Let us go through this step by step:
The geometric composition
$T^1_{\rm Cl}\circ \Sigma_{(2,\ldots,n)}= T^n_{\rm Cl}$ 
is the preimage of $T^1_{\rm Cl}$ under the projection $
(\mu_2\times\ldots\times\mu_n)^{-1}(\tfrac \pi {n+1}, \ldots, \tfrac \pi {n+1})
\to\CP^1$,
hence automatically embedded in the sense of Definition~\ref{embedded}.
Similarly, $T^{n-1}_{\rm Cl} \circ \Sigma_1 = T^n_{\rm Cl}$
is the preimage of $T^{n-1}_{\rm Cl}$ under the projection 
$\mu_1^{-1}(\frac\pi{n+1})\to\CP^{n-1}$,
and by transposition we obtain the embedded composition
$\Sigma_1^t \circ T^{n-1}_{\rm Cl} = T^n_{\rm Cl}$.
Next, the intersection 
$$
\Sigma_{(2,\ldots,n)} \times_{\CP^n} \Sigma_1^t
\; \cong \; (\mu_2\times\ldots\times\mu_n)^{-1}(\tfrac \pi {n+1}, \ldots, \tfrac \pi {n+1})
\cap \mu_1^{-1}(\tfrac \pi {n+1}) \; = \; T^n_{\rm Cl} \;\subset\; \CP^n
$$
is transverse and embeds to 
$$
\Sigma_{(2,\ldots,n)} \circ \Sigma_1^t
= \bigl\{ \bigl( [z_0:z_1],[z_0:z_2:\ldots:z_n] \bigr) \,\big| \, [\ul{z}]\in T^n_{\rm Cl} \bigr\}
= T^1_{\rm Cl} \times T^{n-1}_{\rm Cl} \subset \CP^1 \times \CP^{n-1} .
$$
To make sure that Theorem~\ref{main2} indeed implies all the above isomorphisms of Floer cohomology, it remains to ensure that the maximally decomposed tuple
$(T^1_{\rm Cl}, \Sigma_{(2,\ldots,n)}, \Sigma_j^t,T^{n-1}_{\rm Cl})$ is monotone.
That follows from the monotonicity of all factors together with the torsion fundamental groups of the symplectic manifolds involved. Moreover, it turns out that we need not worry about the minimal Maslov indices $2$. This is since we have $\partial^2=0$ on the first chain group $CF(T_{\rm Cl},T_{\rm Cl})$, i.e.\ assumption (b)(iii) is satisfied (see \cite{oh:fl1} and note Remark~\ref{main rmk}~(a)).
Now Theorem~\ref{main2} implies $\partial^2=0$ (i.e.\ (b)(ii) resp.\ (iii)) for each of the other chain groups in \eqref{calc}.

Moreover, we need to fix spin structures on $T^{n-1}_{\rm Cl}$ and $\Sigma_1$ as well as on $T^1_{\rm Cl}$ and $\Sigma_{(2,\ldots,n)}$ such that the induced spin structure on the composition, $T^{n}_{\rm Cl}$ is the standard one.
For the Clifford tori we pick the standard spin structure given by the trivialization of $TT^k_{\rm Cl}\subset \C^k$ in the coordinate chart $\C^k\cong\{z_0= \frac 1{\sqrt{n+1}}\}\subset\CP^k$. On the sphere $\Sigma_1\subset\CP^n$ we fix the spin structure given by the standard orientation in the chart $\{z_1=\frac 1{\sqrt{n+1}}\}$. (The orientation provides a trivialization over the $0$-skeleton, which coincides with the $1$- and $2$-skeleton of this sphere of dimension $\geq 3$;
see \cite{cho:hol} or \cite{orient} for more details on spin structures.) We can read off the standard spin structure induced on $T^{n}_{\rm Cl}$ from the identification
$$
{\rm T} T^{n}_{\rm Cl} \cong {\rm pr}^* {\rm T} T^{n-1}_{\rm Cl} \oplus E, \qquad
E =  ({\rm pr}^* {\rm T}\CP^{n-1})^\perp \subset {\rm T}\Sigma_1|_{\rm T^{n}_{\rm Cl}}.
$$
Here ${\rm T}\CP^{n-1}|_{T^{n-1}_{\rm Cl}}=  {\rm T} T^{n-1}_{\rm Cl} \oplus i  {\rm T} T^{n-1}_{\rm Cl}$ inherits a trivialization from $T^{n-1}_{\rm Cl}$, so the orientation of $\Sigma_1$ induces a trivialization of the line bundle $E$ (given by the linearized action of $\mu_1$).

For the spin structure on $\Sigma_{(2,\ldots,n)}\subset\CP^n$we identify $\Sigma_{(2,\ldots,n)}\cong T^{n-2}\cdot S^3$ with the orbit of the sphere $S^3=\bigl\{[z_0:z_1:\frac 1{\sqrt{n+1}}:\ldots:\frac 1{\sqrt{n+1}}] \,\big|\, |z_0|^2+|z_1|^2=\frac 2{n+1} \bigr\}\subset\CP^n$ under the action of the torus $T^{n-2}\subset\C^{n-2}$ in the $z_3,\ldots,z_n$-coordinates. If we pick the standard trivialization of $T^{n-2}$ and the standard orientation of $S^3\subset\C^2$ in the above chart, then again the standard spin structure is induced on $T^{n}_{\rm Cl}$ by the identification
$$
{\rm T}_{\ul{z}} T^{n}_{\rm Cl} \cong {\rm T}_{{\rm pr}(\ul{z})} T^1_{\rm Cl} \oplus {\rm T}_{\ul{z}} (T^{n-2} \ul{z}) \oplus F_{(z_0,z_1)}, \qquad
F =  ({\rm pr}^*{\rm T}\CP^1)^\perp \subset {\rm T}S^3|_{T^1_{\rm Cl}}.
$$
Here ${\rm T}\CP^1|_{T^1_{\rm Cl}}=  {\rm T} T^1_{\rm Cl} \oplus i  {\rm T} T^1_{\rm Cl}$ inherits a trivialization from $T^1_{\rm Cl}$, so the orientation of $S^3$ induces a trivialization of the line bundle $F$.

The last isomorphism 
$HF(T^1_{\rm Cl}, T^1_{\rm Cl} \times T^{n-1}_{\rm Cl} ,T^{n-1}_{\rm Cl}) \cong HF(T^1_{\rm Cl}, T^1_{\rm Cl}) \otimes HF(T^{n-1}_{\rm Cl} ,T^{n-1}_{\rm Cl})$
in \eqref{calc} follows from the K\"unneth Theorem~\ref{kuenneth}. 
Here the right hand side is indeed the tensor product of homologies since the first factor is torsion-free.
Indeed, we know from elementary curve counts (see e.g.\ \cite{cho:hol}) that, with the standard spin structure on both factors, $HF(T^1_{\rm Cl},T^1_{\rm Cl})\cong \Z\oplus\Z \cong H_*(S^1=T^1)$. 
Finally, the homology $H_*(T^n)$ satisfies the same inductive relation \eqref{calc} as the Floer cohomology. This proves the theorem.
\end{proof}

This  Floer cohomology calculation directly generalizes when replacing 
$T^1_{\rm Cl}\subset\CP^1$ with another Lagrangian submanifold in a possibly higher dimensional complex projective space.

\begin{theorem} \label{L cl}
Let $1\leq k <n$ and let $L\subset\CP^k$ and $L'\subset\CP^{n-k}$ be oriented, monotone Lagrangian submanifolds. Denote by ${\rm pr}:(\mu_{k+1}\times\ldots\times\mu_n)^{-1}(\tfrac \pi {n+1}, \ldots, \tfrac \pi {n+1})\to\CP^k$ and
${\rm pr'}:(\mu_{1}\times\ldots\times\mu_k)^{-1}(\tfrac \pi {n+1}, \ldots, \tfrac \pi {n+1})\to\CP^{n-k}$ the 
reductions of $\CP^n$ at complementary monotone level sets.
Then ${\rm pr}^{-1}(L)\subset\CP^n$ and ${\rm pr'}^{-1}(L')\subset\CP^n$ are monotone Lagrangian submanifolds and there exists a canonical chain homotopy equivalence
$$
CF({\rm pr}^{-1}(L),{\rm pr'}^{-1}(L'))\cong
CF(L, T^k_{\rm Cl}) \otimes CF(T^{n-k}_{\rm Cl} ,L').
$$
Here we assume that $\partial^2=0$ on either the left hand side or the ride hand side complex; otherwise this is an equivalence in the category of derived matrix factorizations.
\end{theorem}
\begin{proof}
Denote by $\Sigma_{(1,\ldots,k)}\subset(\CP^{n-k})^-\times\CP^n$ the Lagrangian correspondence arising from reduction at the level set 
$(\mu_{1}\times\ldots\times\mu_k)^{-1}(\tfrac \pi {n+1}, \ldots, \tfrac \pi {n+1}) \subset\CP^n$.
Then 
$$
\Sigma_{(k+1,\ldots,n)} \times_{\CP^n} \Sigma_{(1,\ldots,k)}^t
\cong (\mu_{k+1}\times\ldots\mu_n)^{-1}(\tfrac \pi {n+1}, \ldots )
\cap (\mu_1\times\ldots\mu_k)^{-1}(\tfrac \pi {n+1},\ldots) = T^n_{\rm Cl}
$$
is transverse and embeds to 
$\Sigma_{(k+1,\ldots,n)} \circ \Sigma_{(1,\ldots,k)}^t
= T^{k}_{\rm Cl} \times T^{n-k}_{\rm Cl} \subset \CP^k \times \CP^{n-k}$.
Now in complete analogy to the proof of Theorem~\ref{clifford} above, we have 
a sequence of chain homotopy equivalences
\begin{align*}
CF({\rm pr}^{-1}(L),{\rm pr'}^{-1}(L'))
&= CF(L\circ \Sigma_{(k+1,\ldots,n)}, \Sigma_{(1,\ldots,k)}^t \circ L') \nonumber \\
&\cong CF(L, \Sigma_{(k+1,\ldots,n)}, \Sigma_{(1,\ldots,k)}^t,L) \nonumber\\
&\cong CF(L, \Sigma_{(k+1,\ldots,n)}\circ\Sigma_{(1,\ldots,k)}^t,L') \\
&\cong CF(L, T^k_{\rm Cl} \times T^{n-k}_{\rm Cl} ,L') \nonumber\\
&\cong CF(L,T^{k}_{\rm Cl}) \otimes CF(T^{n-k}_{\rm Cl} , L') .
\nonumber
\end{align*}
\end{proof}

A special case of Theorem~\ref{L cl} arises from taking $L'$ to be a Clifford torus, then
$$
HF({\rm pr}^{-1}(L),T_{\rm Cl}^n)
\cong HF(L, T_{\rm Cl}^k) \otimes H_*(T^{n-k}),
$$
where we used $HF(T^{n-k}_{\rm Cl} ,T^{n-k}_{\rm Cl})\cong H_*(T^{n-k})$
by \cite{cho:hol} or Theorem~\ref{clifford}.
Since the latter is torsion-free, the K\"unneth Theorem~\ref{kuenneth} indeed implies the above isomorphism of homologies.
This applies, for example, to $L=\RP^{1}\subset\CP^{1}$ and yields
another Lagrangian torus ${\rm pr}^{-1}(\RP^{1})\subset\CP^n$.
Although it is Hamiltonian isotopic to the Clifford torus, we need not check this but can
calculate directly the relative Floer cohomology
\begin{align*}
HF({\rm pr}^{-1}(\RP^{1}), T^n_{\rm Cl})
\cong HF(\RP^{1}, T^1_{\rm Cl}) \otimes H_*(T^{n-1}) 
\cong H_*(T^n)
\end{align*}
as well as (for $n=2$)
\begin{align*}
HF({\rm pr}^{-1}(\RP^{1}),{\rm pr'}^{-1}(\RP^{1}))
\cong HF(\RP^{1}, T^1_{\rm Cl}) \otimes 
HF(T^1_{\rm Cl},\RP^{1})
\cong H_*(T^1)\otimes H_*(T^1) \cong H_*(T^2) .
\end{align*}
More generally, we can apply Theorem~\ref{L cl} to odd real projective spaces $L=\RP^{k}\subset\CP^{k}$ for $k=2\ell-1\geq 3$ with $\Z_2$-coefficients.\footnote{
The number of holomorphic discs through a generic point is $0$ for $\RP^k$ (which has minimal Maslov number $k+1 \geq 3$ for $k\geq 2$) and it is $k+1$ for $T^{k}_{\rm Cl}$ by \cite{cho:hol}, hence $\partial^2=0$ on
$CF(\RP^{k}, T^{k}_{\rm Cl})$ only holds for odd $k$ and with $\Z_2$ coefficients.
}
By explicit calculation due to Alston \cite{alston} the underlying Floer cohomology is 
$$
HF(\RP^{2\ell-1}, T^{2\ell-1}_{\rm Cl};\Z_2) \cong
 \Z_2^{2^\ell} .
$$
Now our calculations in Theorem~\ref{L cl} provide with $\Z_2$-coefficients
\begin{align*}
HF({\rm pr}^{-1}(\RP^{2\ell-1}), T^{n}_{\rm Cl};\Z_2)
\cong HF(\RP^{2\ell-1}, T^{2\ell-1}_{\rm Cl};\Z_2) \otimes H_*(T^{n-{2\ell+1}}_{\rm Cl};\Z_2) \cong \Z_2^{2^\ell + 2(n-2\ell + 1)}
\end{align*}
as well as for even $n=2m$
\begin{align*}
&HF({\rm pr}^{-1}(\RP^{2\ell-1}), {\rm pr'}^{-1}(\RP^{2(m-\ell) +1});\Z_2) \\
&\cong HF(\RP^{2\ell-1}, T^{2\ell-1}_{\rm Cl};\Z_2) \otimes 
HF(T^{2(m-\ell)+1}_{\rm Cl},\RP^{2(m-\ell)+1};\Z_2) 
\cong \Z_2^{2^\ell + 2^{m-\ell + 1}}.
\end{align*}

\subsection{Detecting nontrivial Floer cohomology of a Lagrangian correspondence}

In this section we provide a tool for deducing nontriviality of Floer cohomology and hence nondisplaceability of a Lagrangian correspondence itself (as Lagrangian submanifold).

\begin{theorem} \label{thm5}
Let $L_{01}\subset M_0^-\times M_1$ be a Lagrangian correspondence.
Suppose that there exists a Lagrangian submanifold $L_1\subset M_1$ such that $L_0:=L_{01}\circ L_1$ is an embedded composition and $HF(L_0,L_0)\neq 0$. 
Assume that $M_0,M_1$ satisfy (M1-2), $L_0,L_1,L_{01}$ satisfy (L1-2), and
$(L_0\times L_1, L_{01})$ is a monotone pair in the sense of Definition~\ref{tuple monotone}~(b).
Then the Lagrangian $L_{01}\subset M_0^-\times M_1$ has nonzero Floer cohomology $HF(L_{01},L_{01})\neq 0$.
\end{theorem}
\begin{proof}
The assumptions guarantee that $HF(L_0,L_0)\cong HF(L_0, L_{01}, L_1) = HF(L_0\times L_1 , L_{01})$ are all well-defined and isomorphic, hence nonzero. (The differential on all three squares to zero since the total disk count in the sense of \eqref{w=0} is $w(L_0)-w(L_0)=0$.)
Now $HF(L_0\times L_1 , L_{01})$ is a module over $HF(L_{01} , L_{01})$, where the multiplication is defined by counting pseudoholomorphic $3$-gons, see  e.g.\ \cite{se:bo} or \cite{cat}. The unit $1_{L_{01}}\in HF(L_{01} , L_{01})$ is defined by counting pseudoholomorphic $1$-gons; it is nontrivial since it acts as identity on a nontrivial group. Hence $HF(L_{01} , L_{01})$ contains a nonzero element, as claimed.
\end{proof}

\begin{corollary} \label{cor of thm5}
Let $\Sigma\subset M$ be the level set of the moment map of a Hamiltonian $G$-action.  Suppose that
$\Sigma$ contains a $G$-invariant Lagrangian submanifold $L\subset M$ such that
$HF(L,L)\neq 0$ and ${\rm pr}(L)\subset \Sigma/G = M\qu G$ is smooth.  
Assume that $M,M\qu G$ satisfy (M1-2), $\Sigma,L,{\rm pr}(L)$ satisfy (L1-2), and
$({\rm pr}(L)\times L, \Sigma)$ is a monotone pair in the sense of Definition~\ref{tuple monotone}~(b).
Then $HF(\Sigma,\Sigma)\neq 0$.
\end{corollary}
\begin{proof}
This is a case of Theorem~\ref{thm5}, where $L_1={\rm pr}(L)$, and the composition $\Sigma\circ {\rm pr}(L)={\rm pr}^{-1}({\rm pr}(L))=L$ is automatically embedded.
\end{proof}

The following example in case $k=n=2$ was initially pointed out to us in 2006 by Paul Seidel; we since learned of alternative proof methods by Biran-Cornea and Fukaya-Oh-Ono-Ohta. We use the notation of Section~\ref{cpn}; in particular 
$\Sigma_{(n)}\subset (\CP^{n-1})^-\times\CP^n$ is a Lagrangian $2n-1$-sphere arising from reduction at the level set $\mu_n^{-1}(\frac\pi{n+1})$.

\begin{corollary} \label{cor:sphere}
For every $2\leq k \leq n$ the Lagrangian embedding 
$\Sigma_{(k,\ldots,n)} \subset (\CP^{k-1})^-\times \CP^n$ 
of $(S^1)^{n-k} \times S^{2k-1}$ is Hamiltonian non-displaceable.
\end{corollary}
\begin{proof}
By construction $\Sigma_{(k,\ldots,n)}$ is the correspondence arising from the level set of
$\mu_k\times\ldots\times\mu_n$ at the level $(\frac\pi{n+1},\ldots,\frac\pi{n+1})$ which contains the nondisplaceable Clifford torus $T^n_{\rm Cl}\subset\CP^n$.
The projection ${\rm pr}(T^n_{\rm Cl})=T^n_{\rm Cl}\circ\Sigma_{(k,\ldots,n)}$ is the Clifford torus $T_{\rm Cl}^{k-1}\subset\CP^{k-1}$. The Clifford tori as well as $\Sigma_{(k,\ldots,n)}$ are monotone with minimal Maslov number $2$ (but this is immaterial here since $\partial^2=0$ on the relevant Floer complex), and the monotonicity of the pair $(T_{\rm Cl}^{k-1}\times T_{\rm Cl}^{n},\Sigma_{(k,\ldots,n)})$ follows directly from the monotonicity of the factors and the torsion fundamental groups of complex projective space.
Now Corollary~\ref{cor of thm5} says that the nonvanishing of $HF(T^n_{\rm Cl}, T^n_{\rm Cl})\neq 0$ (by \cite{cho:hol} or Section~\ref{cpn}) directly implies nonvanishing Floer cohomology $HF(\Sigma_{(k,\ldots,n)},\Sigma_{(k,\ldots,n)})\neq 0$, and hence nondisplaceability.

To spell things out in this example, Theorem~\ref{main2} provides an isomorphism
$$
HF(T^{k-1}_{\rm Cl}\times T^n_{\rm Cl},\Sigma_{(k,\ldots,n)}) = 
HF(T^{k-1}_{\rm Cl},\Sigma_{(k,\ldots,n)}, T^n_{\rm Cl})\cong HF(T^n_{\rm Cl}, T^n_{\rm Cl}) \neq 0.
$$ 
\end{proof}

\subsection{Gysin sequence for spherically fibered Lagrangian correspondence}

In this section, we give a conjectural relation between Floer
cohomology $HF(L,L')$ for $L, L' \subset M_0$ and the Floer cohomology
$HF(L_{01} \circ L, L_{01} \circ L')$ for the images in $M_1$ under a
correspondence $L_{01}\subset M_0^-\times M_1$.  Results of this type
can be viewed as {\em transfer of non-displaceability} results, in the
sense that non-triviality of $HF(L,L')$ implies non-triviality of
$HF(L_{01} \circ L , L_{01} \circ L')$ and hence non-displaceability
of $L_{01} \circ L$ from $L_{01} \circ L'$ by Hamiltonian
perturbation.

In our example, the Lagrangian correspondence arises from a
spherically fibered coisotropic $\iota: C \to M$ with projection $\pi:
C \to B$.  The image of $C$ under $\iota \times \pi$ is a Lagrangian
correspondence from $M$ to $B$, also denoted $C$.  Our standing
assumptions are compactness, orientability, and monotonicity,
i.e.\ $M, B$, and $C$ satisfy (M1-2) and (L1-2) with a fixed $\tau\geq
0$.
Perutz \cite{per:gys} proved the following analogue of the Gysin sequence.

\begin{theorem} 
Suppose that the minimal Maslov number of $C$ is at least $\codim(C
\subset M) +2$.  Then there exists a long exact sequence
$$ \ldots \to HF(C,C) \to HF(\Id) \to HF(\Id) \to HF(C,C) \to \ldots $$
where the map $HF(\Id) \to HF(\Id)$ is quantum multiplication by the
Euler class of $\pi$.
\end{theorem}   

One naturally conjectures the following relative version (for example,
compare the Seidel triangle with the relative version in \cite{se:lo}).

\begin{conjecture} \label{con1}
Let $L_0,L_1 \subset B$ be a monotone pair of Lagrangian submanifolds satisfying (L1-3). 
(Or replace (L3) by $w(L_0)=w(L_1)$.)
Suppose that the minimal Maslov number of $C$ is at least $\codim(C\subset M)+2$
Then there exists a long exact sequence
$$ \ldots \to HF(L_0,C^t,C,L_1) \to HF(L_0,L_1) \to HF(L_0,L_1) \to
HF(L_0,C^t,C,L_1) \to \ldots $$
where the middle map is Floer theoretic multiplication
\footnote{
This product is defined by counting pseudoholomorphic strips with boundary on $(L_0,L_1)$ and an internal puncture with asymptotics fixed by a given class in $HF(\Delta_B)\cong H(B)$.} 
by the Euler class of~$\pi$.
\end{conjecture}   

The compositions $C \circ L_1$ and $L_0\circ C^t = (C \circ L_0)^t$ are clearly embedded. Hence Conjecture~\ref{con1} together with Theorem~\ref{main2} implies the following.

\begin{corollary} 
Under the same assumptions as in Conjecture~\ref{con1} there exists a long exact
sequence
$$ \ldots \to HF(C \circ L_0,C \circ L_1) \to HF(L_0,L_1) \to
HF(L_0,L_1) \to HF(C \circ L_0,C \circ L_1) \to \ldots $$
\end{corollary}

In particular, we obtain a 'transfer of non-displaceability' result if the Euler class vanishes.

\begin{corollary} With the same assumptions as in Corollary~\ref{con1}, 
if the Euler class of $\pi : C \to B$ is zero, then 
$ HF(C \circ L_0, C \circ L_1)$ is isomorphic to two copies
of $HF(L_0,L_1)$. 
\end{corollary}

\begin{example} 
Suppose that $M$ is a monotone Hamiltonian $G=SU(2)$ manifold, with
moment map $\Phi$, and $\Phinv(0)$ is an $SU(2)$-bundle over the
symplectic quotient $M \qu G$.

Let $(L_0,L_1)$ be a monotone pair of $G$-invariant Lagrangians contained in the zero level set and with minimal Maslov number at least three.
Necessarily each $L_j$ is a principal $SU(2)$ bundle over $L_j / G
\subset M \qu G$.  Suppose that the minimal Maslov number of
$\Phinv(0)$, considered as a Lagrangian in $M^- \times M \qu G$, is at
least $5$. Then there is a long exact sequence
$$ \ldots \to HF(L_0,L_1) \to HF(L_0 / G,L_1/G) \to HF(L_0/G,L_1/G)
\to HF(L_0,L_1) \to \ldots .$$
In particular, if $M \to M \qu G$ is a trivial $G$-bundle, then
$HF(L_0,L_1)$ is isomorphic to two copies of $HF(L_0/G,L_1/G)$.
\end{example}

\def\cprime{$'$} \def\cprime{$'$} \def\cprime{$'$} \def\cprime{$'$}
  \def\cprime{$'$} \def\cprime{$'$}
  \def\polhk#1{\setbox0=\hbox{#1}{\ooalign{\hidewidth
  \lower1.5ex\hbox{`}\hidewidth\crcr\unhbox0}}} \def\cprime{$'$}
  \def\cprime{$'$}

\end{document}